\documentclass[12pt,francais]{article}
\usepackage[latin1]{inputenc}
\usepackage{babel}
\usepackage{amsmath,amsthm,amsfonts,amssymb,amscd}
\usepackage{epsf,graphicx}
\usepackage{color}
\evensidemargin = \oddsidemargin
\advance\textwidth by -2in
\makeatletter
\gdef\thmhead@plain#1#2#3{%
  \thmname{#1}\thmnumber{\@ifnotempty{#1}{ }#2}%
  \thmnote{ {\mdseries#3}}}
\let\thmhead\thmhead@plain
\makeatother
\theoremstyle{plain}
\newtheorem*{theoreme*}{Théorème}
\newtheorem{theoreme}{Théorème}[section]
\newtheorem{proposition}[theoreme]{Proposition}
\newtheorem{lemme}[theoreme]{Lemme}
\newtheorem{lemme-d'elagage}[theoreme]{Lemme d'élagage}
\newtheorem{corollaire}[theoreme]{Corollaire}
\newtheorem{assertion}[theoreme]{Assertion}
\theoremstyle{definition}
\newtheorem{definition}[theoreme]{Définition}
\newtheorem{question}[theoreme]{Question}
\theoremstyle{remark}
\newtheorem{remarque}[theoreme]{Remarque}

\newtheorem*{remarque*}{Remarque}
\newtheorem*{remarques*}{Remarques}
\newtheorem{exemple}[theoreme]{Exemple}

\newtheorem*{resume}{Résumé}
\newtheorem*{remerciements}{Remerciements}
\def\ly{\fontencoding{U}\fontfamily{lasy}\fontseries{m}\fontshape{n}\selectfont}
\def\guil#1{\leavevmode\hbox{{\ly(\kern-0.20em(\kern+0.20em}}\nobreak{}\,#1\,%
  \nobreak\hbox{{\ly\kern+0.20em)\kern-0.20em)}}}
\def\Alinea#1{\hfill\break%
  \hbox to \parindent{\hss{\textup{#1}}\enspace}\ignorespaces}
\def\alinea#1{\par\noindent%
  \hbox to \parindent{\hss{\textup{#1}}\enspace}\ignorespaces}
\def\up{\textup}
\def\from{\colon}
\def\res#1{\,\vert\,{}_{#1}}
\def\abs#1{\lvert #1 \rvert}

\def\classe#1{\mathcal C_{}^{#1}}
\def\eps{\varepsilon}
\def\Chi{\setbox0=\hbox{$\chi$} \mathord{\raise\dp0\hbox{$\chi$}}}

\def\D{{\mathbf D}}

\def\KK{{\mathcal K}}

\def\N{{\mathbb N}}
\def\P{{\mathbf P}}
\def\S{{\mathbf S}}

\def\R{{\mathbb R}}
\def\SS{{\mathcal S}}
\def\T{{\mathbf T}}
\def\TT{{\mathcal T}}
\def\F{{\mathcal F}}
\def\X{{\mathbf X}}
\def\Z{{\mathbb Z}}
\def\CP^#1{\mathbf P^{#1}(\mathbf C)}
\def\RP^#1{\mathbf P^{#1}(\mathbf R)}
\def\XX{{\mathcal X}}
\def\FF{{\mathcal F}}

\def\SCT{\mathop{\mathcal{SCT}}}
\def\id{\mathrm{id}}
\def\Adh{\mathop{\mathrm{Adh}}\nolimits}
\def\Int{\mathop{\mathrm{Int}}\nolimits}
\def\Reg{\mathop{\mathrm{Reg}}\nolimits}

\def\Card{\mathop{\mathrm{Card}}}

\def\tb{\mathop{\mathrm{tb}}}
\def\TB{\mathop{\mathrm{TB}}}
\def\arcinf{\partial_-}
\def\arcsup{\partial_\cap}
\def\grad{\nabla}

\def\adresse#1{\def\nl{\egroup\egroup\hbox\bgroup\itshape\bgroup}
  \par \noindent
  \hbox to \textwidth{\hskip .5\textwidth \vbox{\small
  \hbox \bgroup\itshape\bgroup #1 \egroup\egroup}\hfill}}

\title{Finitude homotopique et isotopique\\ des structures de contact tendues}
\author{Vincent Colin, Emmanuel Giroux et Ko Honda}
\date{}
\begin{document}

\maketitle

\begin{resume}
Soit $V$ une variété close de dimension~$3$. Dans cet article, on montre que les
classes d'homotopie de champs de plans sur~$V$ qui contiennent des structures de
contact tendues sont en nombre fini et que, si $V$ est atoroïdale, les classes
d'isotopie des structures de contact tendues sur $V$ sont elles aussi en nombre 
fini.

\smallskip

\noindent
\textit{Mots clés}~:
structure de contact tendue, champ de plans, homotopie, isotopie, conjugaison,
triangulation, surface branchée.

\smallskip

\noindent
\textit{Codes AMS}~: 57M50, 53C15.
\end{resume}

\bigskip

Le but de cet article est de démontrer les deux théorèmes suivants ainsi que 
quelques uns de leurs avatars~:

\begin{theoreme} \label{t:homotopie}
Sur toute variété close de dimension~$3$, les classes d'homotopie de champs de 
plans qui contiennent des structures de contact tendues sont en nombre fini.
\end{theoreme}

\begin{theoreme} \label{t:isotopie}
Sur toute variété close et atoroïdale de dimension~$3$, les classes d'isotopie 
des structures de contact tendues sont en nombre fini.
\end{theoreme}

Ces résultats s'inscrivent dans un contexte qu'on rappelle brièvement. D'abord, 
les structures de contact tendues sont officiellement nées dans~\cite{El2} où 
Y.~Eliashberg établit à leur sujet le théorème de finitude homologique suivant~:
sur toute variété close~$V$ de dimension~$3$, seul un nombre fini de classes de 
cohomologie entières de degré $2$ sont réalisables comme les classes d'Euler de
structures de contact tendues (orientées). Dans ce même article, le théorème de
finitude homotopique~\ref{t:homotopie} apparaît aussi explicitement --~Theorem 
2.2.2~-- mais l'argument proposé pour sa démonstration est incomplet. Dans~\cite
{KM}, par l'étude des monopoles de Seiberg-Witten, P.~Kronheimer et T.~Mrowka 
obtiennent l'analogue du théorème~\ref{t:homotopie} pour les structures de 
contact remplissables (au sens le plus faible qui soit). D'après un théorème de 
M.~Gromov~\cite{Gr} et Y.~Eliashberg~\cite{El3}, toute structure de contact 
remplissable est tendue mais bien des structures de contact tendues ne sont, en 
revanche, pas remplissables~\cite{EH2, LS}, de sorte que le théorème~\ref
{t:homotopie} ne découle pas des résultats de~\cite{KM}.

Un autre théorème de \cite{El2} prouve que la sphère $\S^3$ porte une seule 
structure de contact tendue à isotopie (et changement d'orientation) près. Ceci 
suggère que les classes d'isotopie des structures de contact tendues pourraient 
être en nombre fini sur toute variété close de dimension~$3$ \cite[conjecture 
8.6.1]{El2} mais les résultats de \cite{Gi3,Ka} réfutent cette hypothèse. Ils
montrent en effet qu'une modification de Lutz répétée sur un tore incompressible
peut produire une infinité de structures de contact tendues deux à deux non 
conjuguées, d'où la conjecture suivante de~\cite{Gi4}~: les structures de 
contact tendues sur une variété close de dimension~$3$ forment un nombre fini de
classes d'isotopie si et seulement si la variété est atoroïdale, \emph{i.e.} ne 
contient aucun tore plongé incompressible. Cette conjecture est à présent avérée
grâce au théorème~\ref{t:isotopie} ci-dessus et au résultat suivant démontré
indépendamment dans \cite{Co3} et dans \cite{HKM} (travail commun avec W.~Kazez 
et G.~Mati\'c)~:

\begin{theoreme}
Sur toute variété de dimension~$3$ close, orientée, toroïdale et
irréductible\footnote
{\upshape{Sans cette hypothèse, il se peut que la variété ne porte aucune 
structure de contact tendue (voir \cite{EH1}.}},
les structures de contact tendues forment une infinité de classes de
conjugaison.
\end{theoreme}

\medskip

Les méthodes développées dans le présent article permettent par ailleurs de 
préciser ce résultat en tenant compte de la \emph{torsion}, invariant des 
structures de contact introduit dans \cite{Gi4} et défini comme suit~:

\begin{definition} \label{d:torsion}
On appelle \emph{torsion} d'une structure de contact~$\xi$ sur une variété close
$V$ de dimension~$3$ la borne supérieure des entiers $n \in \N$ pour lesquels on
peut plonger dans $(V,\xi)$ le produit $\T^2 \times [0, 2n\pi]$ muni de la 
structure de contact d'équation $\cos \theta \, dx_1 - \sin \theta \, dx_2 = 0$,
où $(x, \theta) \in \R^2 \!/ \Z^2 \times [0, 2n\pi]$.
\end{definition}

En vertu du théorème de Y.~Eliashberg qui les classifie, les structures de 
contact non tendues --~aussi dites \emph{vrillées}~-- ont une torsion infinie. 
En outre, sur une variété atoroïdale et irréductible, les structures de contact 
tendues ont une torsion nulle. Pour ce qui concerne les variétés toroïdales, les
démonstrations du théorème ci-dessus établissent en fait que, lorsqu'elles sont 
irréductibles, elles admettent des structures de contact tendues de torsion 
finie arbitrairement grande. Cela dit, la réponse à la question suivante reste
inconnue~: 

\begin{question}
Les structures de contact tendues ont-elles toutes une torsion finie\,?
\end{question}

Le théorème qu'on obtient ici est le suivant~:

\begin{theoreme} \label{t:torsion}
Sur toute variété close de dimension~$3$, les structures de contact tendues et 
de torsion égale à un entier fixé ne forment qu'un nombre fini d'orbites sous
l'action du groupe engendré par les difféomorphismes isotopes à l'identité et 
les twists de Dehn le long de tores.
\end{theoreme}

Ce théorème entraîne le théorème \ref{t:isotopie} puisque, sur une variété~$V$ 
atoroïdale, tout twist de Dehn le long d'un tore est isotope à l'identité. Sur 
une variété toroïdale en revanche, les structures de contact tendues de torsion 
(finie) fixée peuvent très bien former une infinité de classes d'isotopie, comme
le montre l'exemple des fibrés en cercles au-dessus des surfaces~\cite{Gi5,Ho2}.

\medskip

La clé des théorèmes \ref{t:homotopie}, \ref{t:isotopie} et \ref{t:torsion} est 
le résultat suivant qui montre que les structures de contact tendues sur une 
variété close sont engendrées, à partir d'un nombre fini d'entre elles, par des 
modifications de Lutz (voir la partie 1.4 pour la définition de cette 
opération)~:

\begin{theoreme} \label{t:generation}
Sur toute variété close~$V$ de dimension~$3$, il existe un nombre fini de 
structures de contact $\xi_1, \dots, \xi_n$ et, pour chaque entier $i \in \{1, 
\dots, n\}$, un nombre fini de tores $T_1^i, \dots, T_{k_i}^i$ transversaux à 
$\xi_i$ tels que toute structure de contact tendue $\xi$ sur~$V$, à isotopie 
près, s'obtienne à partir d'une des structures~$\xi_i$ par une modification de 
Lutz de coefficients $n_j^i(\xi) \in \N$ le long des tores $T_j^i$, $1 \le j \le
k_i$.
\end{theoreme}

Un point crucial néanmoins --~qui empêche de déduire formellement le théorème
\ref{t:torsion} du théorème \ref{t:generation}~-- est ici que les tores $T_j^i$ 
fournis par ce dernier ne sont \emph{a priori} pas disjoints --~même pour~$i$ 
fixé. Du coup, l'ordre dans lequel on effectue les modifications de Lutz a une 
incidence.

\medskip

Tous ces résultats fournissent une classification géométrique grossière des 
structures de contact tendues en dimension~$3$. Il serait intéressant de la 
préciser en donnant par exemple des estimations effectives pour les nombres dont 
on se borne ici à prouver qu'ils sont finis.

Les preuves qu'on donne dans ce texte s'\'etendent 
au  cas des vari\'et\'es compactes \`a bord, pour peu qu'on impose
un germe de structure de contact le long du bord.
Par exemple, le th\'eor\`eme~\ref{t:torsion} devient~:

\begin{theoreme} \label{t:torsion-bord}
Soit $V$ une  variété compacte \`a bord de dimension~$3$ et $\zeta$ un
germe de structure de contact le long de $\partial V$. 
Les structures de contact tendues \'egales \`a $\zeta$ pr\`es de $\partial V$ et 
de torsion égale à un entier fixé ne forment qu'un nombre fini d'orbites sous
l'action du groupe engendré par les difféomorphismes isotopes à l'identité 
relativement au bord et 
les twists de Dehn le long de tores inclus dans $\Int (V)$.
\end{theoreme}

\medskip

Grâce au travail de Y.~Eliashberg et W.~Thurston~\cite{ET} qui procure, en 
dimension~$3$, une méthode générale pour déformer un feuilletage par surfaces en
une structure de contact, le résultat suivant de D.~Gabai~\cite{Ga} se déduit du
théorème~\ref{t:homotopie}~:

\begin{corollaire} \label{c:feuilletages}
Sur toute variété de dimension~$3$ close et orientable, les classes d'homotopie 
de champs de plans qui contiennent un feuilletage sans composantes de Reeb sont 
en nombre fini.
\end{corollaire}

Avant D.~Gabai, P.~Kronheimer et T.~Mrowka ont démontré ce théorème dans \cite
{KM} pour une classe plus restreinte de feuilletages --~les feuilletages 
tendus~-- comme une conséquence de leur travail sur les structures de contact 
remplissables. En fait, un des résultats principaux de~\cite{ET} dit que tout 
feuilletage tendu sur une variété close orientable peut être déformé en une 
structure de contact remplissable. Par suite, la finitude homotopique des
structures de contact remplissables implique celle des feuilletages tendus. Le 
corollaire~\ref{c:feuilletages} découle pareillement du théorème~\ref
{t:homotopie} et d'un résultat de~\cite{Co4} qui montre que, toujours sur une 
variété close et orientable, tout feuilletage sans composantes de Reeb se laisse
déformer en une structure de contact tendue.

Dans \cite{Ga}, D.~Gabai démontre le corollaire~\ref{c:feuilletages} par une 
méthode géométrique fondée sur les idées classiques de H.~Kneser~\cite{Kn}, 
W.~Haken~\cite{Ha} et F.~Waldhausen~\cite{Wa} mais utilisant de plus un 
processus infini de mise sous forme normale dû à M.~Brittenham~\cite{Br}. Cette 
méthode lui permet d'atteindre aussi un analogue du théorème~\ref{t:isotopie} 
qui, pour le coup, semble loin d'en être un corollaire.

Enfin, la version  relative~\ref{t:torsion-bord} 
du théorème~\ref{t:torsion} donne une classification 
grossière des n\oe uds legendriens dans la sphère $\S^3$ 
munie de sa structure de contact standard~: 

\begin{theoreme} \label{t:noeuds}
Dans la sphère de contact standard de dimension~$3$, les n\oe uds legendriens
ayant une classe d'isotopie lisse et un invariant de Thurston-Bennequin imposés 
--~arbitrairement~-- forment un nombre fini de classes d'isotopie legendrienne.
\end{theoreme}

\medskip

Pour démontrer les résultats énoncés ci-dessus, la méthode qu'on adopte consiste 
à mettre chaque structure de contact tendue sur~$V$, par isotopie, en position 
\guil{adaptée} et \guil{minimale} par rapport à une triangulation fixe $\Delta$
de~$V$ (sections 2.1 et 2.2.). Ensuite, grâce à la classification des structures
de contact tendues sur la boule~\cite{El2} et à des manipulations de surfaces 
convexes~\cite{Gi1}, on observe que chaque structure de contact tendue obtenue 
est entièrement décrite par une donnée combinatoire~: le découpage des faces de 
$\Delta$. Il s'agit, sur chaque $2$-simplexe~$F$, d'une collection finie d'arcs 
simples, propres, disjoints et évitant les sommets, leur union n'étant bien 
définie qu'à une isotopie près de~$F$ relative aux sommets. Ces arcs (à cause
des propriétés de la triangulation et des structures de contact), s'organisent 
en familles d'arcs parallèles, elles-mêmes contenues dans un nombre fini de 
\guil{domaines fibrés} (sections 2.3 et 2.4). En dehors de ces domaines, les
structures sont sous contrôle (théorème~\ref{t:domaines}). On établit alors une 
correspondance entre les structures de contact tendues \guil{ajustées} à un même
domaine fibré $M$  et les surfaces \guil{portées} par $M$. C'est là 
qu'apparaissent les tores générateurs du théorème~\ref{t:generation} (section
3). Pour démontrer la finitude isotopique, il faut encore simplifier le domaine 
$M$ pour le faire ressembler à un fibré en cercles au-dessus d'une surface 
(section 4). On applique alors un théorème classifiant les structures tendues 
portées par un fibré en cercles~\cite{Gi4, Gi5, Ho2} pour conclure. Enfin, on 
donne la preuve du théorème~\ref{t:noeuds} dans la section~5. La section~1 est 
consacrée à la présentation  des outils nécéssaires à la réalisation de ce 
programme.

\begin{remerciements}
Le travail présenté ici est, pour l'essentiel, le fruit de discussions entamées 
à l'American Institute of Mathematics --~Palo Alto, USA~-- à l'automne 2000 et 
poursuivies à l'Institut des Hautes Études Scientifiques --~Bures-sur-Yvette,
France~-- au printemps 2001. Les auteurs remercient vivement ces deux organismes
pour leur accueil et leur aide matérielle. Ils remercient également le CNRS, la 
NSF, la Fondation Alfred Sloan et l'Institut Universitaire de France pour avoir 
en partie financé leurs recherches. Ils tiennent enfin à adresser un grand merci
à François Laudenbach pour sa relecture attentive d'une ébauche du texte.
\end{remerciements}

\section{Outils de géométrie de contact et de topologie}

On présente ici succintement les outils de géométrie de contact et
de topologie qui serviront pour démontrer les théorèmes annoncés
dans l'introduction.

\subsection{Structures de contact tendues\,/\,vrillées}\label{subsection:tendues}

Dans ce texte, $V$ désigne toujours une variété compacte et orientée
de dimension~$3$ et chaque structure de contact~$\xi$ qu'on
considère sur~$V$ est directe et (co-)\,orientée. En d'autres
termes, il existe sur~$V$ une $1$-forme $\alpha$ --~unique à
multiplication près par une fonction positive~-- dont $\xi$ est le
noyau coorienté et dont le produit avec $d\alpha$ est partout un
élément de volume positif relativement à l'orientation de~$V$.

Sur une variété close, les structures de contact forment au plus une
infinité dénombra\-ble de classes d'isotopie. Cela découle du
théorème de Gray qui montre qu'elles sont stables dans la topologie
$\classe1$~: si $\xi_s$, $s \in [0,1]$, est un chemin de structures
de contact sur une variété close~$V$, il existe une isotopie
$\phi_s$ de $V$ partant de l'identité ($\phi_0 = \id$) et vérifiant
$\phi_{s*}\xi_0 = \xi_s$ pour tout $s$ --~en dimension~$3$, on
dispose même d'un résultat beaucoup plus fort de stabilité
$\classe0$~\cite{Co2}. Un autre fait important est qu'une variété de
contact close $(V,\xi)$ possède un grand nombre d'isotopies de
contact. Plus explicitement, la projection $TV \to TV / \xi$ met en
bijection les champs de vecteurs sur~$V$ qui préservent~$\xi$ et les
sections du fibré quotient $TV / \xi \to V$. Autrement dit, si $\xi$
est le noyau d'une forme de contact $\alpha$, toute fonction lisse
$u \from V \to \R$ détermine un unique champ de vecteurs
$\grad_\alpha u$ qui conserve~$\xi$ et sur lequel $\alpha$ vaut~$u$.
Le problème du prolongement des isotopies --~ou des champs de
vecteurs~-- de contact se ramène en particulier à un problème de
prolongement de fonctions.

\medskip

Un \emph{disque vrillé}~$D$ dans une variété de contact $(V, \xi)$
est un disque plongé dans~$V$ et dont le plan tangent $T_pD$, en
chaque point~$p$ du bord $\partial D$, coïncide avec le plan de
contact~$\xi_p$. La structure de contact $\xi$ est dite
\emph{vrillée} ou \emph{tendue} selon qu'un disque vrillé existe ou
non à l'intérieur de~$V$. Cette dichotomie a deux raisons d'être
majeures. En premier lieu, un théorème de Y.~Eliashberg~\cite{El1}
établit la classification isotopique des structures de contact
vrillées sur toute variété close~: chaque classe d'homotopie de
champs de plans tangents contient une et une seule classe d'isotopie
de structures de contact vrillées. D'autre part, les disques vrillés
s'interprètent comme des cycles évanescents. Plus précisément,
Y.~Eliashberg et W.~Thurston ont proposé dans~\cite{ET} la
définition suivante~: un \emph{cycle évanescent} pour un champ de
plans quelconque $\zeta$ est un disque plongé~$D$ qui est partout
tangent à~$\zeta$ le long de son bord mais n'est pas isotope à un
disque intégral de~$\zeta$ relativement à $\partial D$. Avec ce
langage, une structure de contact (resp. un feuilletage) sans cycles
évanescents est une structure de contact tendue (resp. un
feuilletage sans composantes de Reeb~: théorème de Novikov).

La classification des structures de contact tendues n'est connue que
sur peu de variétés, essentiellement celles que des découpages
successifs le long de tores et d'anneaux ramènent à des boules. Sur
la boule, un théorème de Y.~Eliashberg \cite{El2} affirme qu'il n'y
a essentiellement qu'une structure tendue.

\begin{theoreme}\label{t:boule} Des structures de contact tendues
sur la boule qui impriment sur le bord le même champ de droites
singulier (ou, autrement dit, le même feuilletage caractéristique~:
voir plus loin) sont isotopes relativement au bord.
\end{theoreme}

Dans la suite, une variété close et orientée $V$ de dimension~$3$
étant donnée, $\SCT(V)$ désigne l'espace de ses structures de
contact tendues et on dit qu'un sous-ensemble $\XX$ de $\SCT(V)$ est
\emph{complet} s'il représente toutes les classes d'isotopie,
c'est-à-dire s'il possède (au moins) un point dans chaque composante
connexe de $\SCT(V)$.

\subsection{Surfaces convexes}

La théorie des surfaces convexes~\cite{Gi1} fournit, pour les
variétés de contact, des techniques de découpage et de collage le
long de surfaces. On peut brièvement la résumer ainsi~: au voisinage
d'une surface générique, la géométrie d'une structure de contact est
entièrement inscrite dans une multi-courbe tracée sur la surface.

\smallskip

Toute surface $S$ plongée dans une variété de contact $(V, \xi)$
hérite d'un \emph{feuilletage caractéristique} noté $\xi S$ et
engendré par le champ singulier de droites $\xi \cap TS$. Ce
feuilletage détermine totalement le germe de $\xi$ le long de~$S$,
à isotopie de $V$ préservant $S$ et $\xi S$ près. Lorsque
$\xi$ et~$S$ sont orientées, $\xi S$ l'est aussi. Concrètement, si
$\xi$ est le noyau de~$\alpha$ et si $\omega$ est une forme d'aire
sur~$S$, le feuilletage~$\xi S$ est constitué des courbes intégrales
orientées du champ de vecteurs~$\eta$ dont le produit intérieur
avec~$\omega$ est égal à la $1$-forme induite par $\alpha$ sur~$S$.
Une singularité de $\xi S$ est un point de tangence entre $\xi$
et~$S$ et a donc un signe~: elle est positive ou négative selon que
les orientations de $\xi$ et~$S$ coïncident ou non. Ce signe est
également celui de la divergence de~$\eta$. Les singularités de $\xi
S$ sont aussi appelées \emph{points complexes} de~$S$. Un point
complexe est dit \emph{elliptique} (resp. \emph {hyperbolique}) si
c'est une singularité isolée de $\xi S$ qui est de type foyer (resp.
selle), c'est-à-dire d'indice~$1$ (resp.~$-1$).

\medskip

Soit $S$ une surface compacte, orientée, plongée dans $(V,\xi)$ et à
bord vide ou legendrien --~une courbe legendrienne est une courbe
partout tangente à~$\xi$. On dit que $S$ est \emph{$\xi$-convexe} si
elle a un épaississement
$$ U = S \times \R \supset S = S \times \{0\} $$
dans lequel l'action de~$\R$ par translations (action engendrée par
le champ de vecteurs $\partial_t$, $t \in \R$) préserve~$\xi$. Un
tel épaississement $U$ est dit \emph{homogène}.

Une surface $\xi$-convexe possède non seulement un épaississement
homogène mais tout un \guil{système fondamental} de tels voisinages.
En d'autres termes, tout voisinage de $S$ contient un voisinage
homogène de $S$. En effet, si $\nu$ est la section de $TV/\xi$
associée au champ de vecteurs $\partial_t$ dans un épaississement
homogène $U = S \times \R$, le fait que $\partial S$ soit legendrien
assure que le champ de vecteurs de contact associé à toute section
du type $f \nu$, où $f$ est une fonction qui ne dépend que de~$t$,
reste tangent à~$\partial U = \partial S \times \R$. En effet, pour
$\xi = \ker \alpha$ et $v = \alpha (\nu)$, si on cherche  le champ
de contact $\grad_\alpha (fv)$ (défini dans la section~1.1) sous
la forme $\grad_\alpha (fv )=f\frac{\partial}{\partial t} +X_0$,
avec $X_0 \in \xi$, cela résulte de la formule
$$i_{X_0} d\alpha \res \xi = -vdf \res \xi.$$

La $\xi$-convexité est ainsi une propriété locale qui, pour peu que
$S$ ne touche pas $\partial V$, se lit directement sur le
feuilletage caractéristique $\xi S$ comme expliqué ci-dessous.

\medskip

On considère sur~$S$ un feuilletage singulier~$\sigma$ --~famille
des courbes intégrales orien\-tées d'un champ de vecteurs~-- et
une multi-courbe~$\Gamma$ --~union finie de courbes fermées, simples
et disjointes. On dit que $\Gamma$ est une courbe de
\emph{découpage} de $\sigma$ (où \emph{découpe} $\sigma$) si
$\sigma$ est transversal à~$\Gamma$ et est porté par un champ de
vecteurs qui, sur chaque composante de $S \setminus \Gamma$, ou bien
dilate l'aire et sort le long de~$\Gamma$, ou bien contracte l'aire
et rentre le long de~$\Gamma$. Une telle courbe, si elle existe et
si $\sigma$ est tangent à $\partial S$, est unique à isotopie près
parmi les courbes de découpage.

Le lemme qui suit est établi dans~\cite{Gi1} pour les surfaces
closes mais sa démonstration s'étend au cas des surfaces à bord
legendrien~:

\begin{lemme} \label{l:convexite}
Soit $\xi$ une structure de contact sur $V$ et $S \subset \Int (V)$
une surface compacte, orientée et à bord legendrien. Pour que $S$
soit $\xi$-convexe, il faut et il suffit que son feuilletage
caractéristique $\xi S$ soit découpé par une multi-courbe. De
plus, cette condition est remplie dès que les propriétés suivantes
sont satisfaites\up{~:}
\begin{enumerate}
\item
chaque demi-orbite de $\xi S$ a pour ensemble limite une singularité
ou une orbite fermée\up{\,;}
\item
les orbites fermées de $\xi S$ sont toutes hyperboliques\up{\,;}
\item
aucune orbite \up(orientée\up) de $\xi S$ ne va d'une singularité
négative à une singularité positive.
\end{enumerate}
\end{lemme}
\begin{remarque} \label{r:demitour}
La dernière des trois propriétés ci-dessus mérite un commentaire~:
une orbite de $\xi S$ qui va d'une singularité négative à une
singularité positive est un arc legendrien sur $S$ le long duquel
$\xi$ effectue globalement autour de $S$ un demi-tour dans le sens
contraire des aiguilles d'une montre.
\end{remarque}

Soit $S$ une surface compacte orientée plongée dans $(V,\xi)$.
Compte tenu du théorème de M.~Peixoto qui affirme qu'un feuilletage
singulier de~$S$ est génériquement de type Morse-Smale (et donc
possède les propriétés 1, 2 et 3), le lemme \ref{l:convexite}
ci-dessus montre~\cite{Gi1} que, si $S$ est close, on peut la rendre
$\xi$-convexe en la perturbant par une isotopie arbitrairement
$\classe\infty$-petite. De même, si $S$ est à bord legendrien et si
aucune orbite de $\xi S$ sur $\partial S$ ne va d'une singularité
négative vers une singularité positive (\emph{i.e.} si, entre deux
singularités consécutives de signes contraires sur $\partial S$, le
demi-tour que $\xi$ fait autour de~$S$ s'effectue globalement dans
le sens des aiguilles d'une montre), on peut aussi la rendre
$\xi$-convexe en la dépla\c{c}ant par une isotopie relative au bord
et arbitrairement $\classe\infty$-petite.

Soit maintenant $L$ une courbe legendrienne tracée sur~$S$. On
appelle \emph{nombre de Thurston-Bennequin} de $L$ relatif à~$S$ le
nombre d'enroulement $\tb(L,S) \in \tfrac12 \Z$ de $\xi$ autour de
$TS$ le long de~$L$ (si $L$ est fermée, $\tb(L,S) \in \Z$), qui
décompte algébriquement les points complexes de $S$ le long de
$L$. On vérifie facilement que, quand $S$ est $\xi$-convexe et que
les points extrêmes de $L$ sont complexes --~si $L$ n'est pas
fermée~--,
$$ \tb(L,S) = - \tfrac12 \Card (L \cap \Gamma) $$
où $\Gamma$ est une multi-courbe qui découpe $\xi S$. En
particulier, chaque composante de $\partial S$ a un invariant de
Thurston-Bennequin négatif ou nul.  Inversement, si $S$ est une
surface compacte à bord legendrien et si chaque composante de
$\partial S$ a un nombre de Thurston-Bennequin négatif ou nul, une
isotopie arbitrairement $\classe0$-petite, relative au bord et à
support dans un voisinage aussi petit qu'on veut de $\partial S$,
permet de perturber~$S$ en une surface dont aucune orbite du
feuilletage caractéristique sur le bord ne lie une singularité
négative à une singularité positive. Une isotopie arbitrairement
$\classe\infty$-petite suffit ensuite pour obtenir une surface
$\xi$-convexe.
La condition portant sur le signe des nombres de Thurston-Bennequin des composantes
de $\partial S$ est en outre ais\'ee \`a satisfaire au moyen d'une isotopie, dite
de {\it stabilisation}~:  
pour toute courbe legendrienne $L$ tracée sur $S$, pour tout point $p
\in L$, pour tout voisinage $N(p)$ de $p$ dans $V$ et pour tout $n \in \N$, il 
existe une isotopie $(\phi_t)_{t \in [0,1]}$ de $V$ à support dans $N(p)$ telle 
que $\phi_1(L)$ soit une courbe legendrienne et $\tb (\phi_1(L), \phi_1(S)) = 
\tb (L,S) - n$. 

\medskip

On suppose désormais que $S$ est $\xi$-convexe. Le choix d'un
épaississement homogène $U = S \times \R$ détermine sur~$S$ une
multi-courbe $\Gamma_{\!U}$ qui découpe $\xi S$~: elle est formée
des points de~$S$ où le vecteur $\partial_t$, $t \in \R$, appartient
au plan~$\xi$. De plus, toute multi-courbe qui découpe $\xi S$ est
le découpage associé à un certain épaisssissement homogène. Le lemme
ci-dessous~\cite {Gi1} montre que la multi-courbe $\Gamma_{\!U}$
concentre en fait toute la géométrie de contact de~$U$\,:

\begin{lemme}{\rm [Lemme de réalisation de feuilletages]}
\label{l:feuilletages} Soit $S$ une surface $\xi$-convexe, $U$ un
épaississement homogène de~$S$ et $\Gamma$ le découpage de $\xi S$
associé. Les feuilletages singuliers de $S$ tangents à~$\partial S$
et découpés par~$\Gamma$ constituent un espace contractile. De
plus, pour tout feuilletage~$\sigma$ dans cet espace, il existe un
plongement $\phi \from S \to U$ ayant les propriétés
suivantes\up{~:}
\begin{enumerate}
\item
$\phi$ coïncide avec l'inclusion sur $\partial S$\up{\,;}
\item
$\phi(S) \subset U = S \times \R$ est une surface transversale au
champ de vecteurs $\partial_t$, $t \in \R$\up{\,;}
\item
$\phi_*\sigma$ est le feuilletage caractéristique de $\phi(S)$.
\end{enumerate}
\end{lemme}

On note que, vu la propriété~2, $U$ est un épaississement homogène
de $\phi(S)$ et que le découpage associé est l'intersection de
$\phi(S)$ avec le cylindre $\Gamma \times \R$. En outre, $\phi$ est
isotope à l'inclusion parmi les plongements $S \to U$ satisfaisant
les conditions 1) et 2).

Un corollaire direct du lemme ci-dessus est le résultat suivant
\cite[exemple II.3.7]{Gi1}~:

\begin{lemme}{\rm [Lemme de réalisation legendrienne]} \label{l:courbes}
Soit $S$ une surface $\xi$-convexe, $U$ un de ses épaississements
homogènes, $\Gamma$ la multi-courbe de découpage de $\xi S$ associé
et $C$ une multi-courbe sur~$S$ transversale à~$\Gamma$ et telle que
chaque composante de $S \setminus C$ rencontre $\Gamma$. Il existe
un plongement $\phi \from S \to U$ ayant les propriétés
suivantes\up{~:}
\begin{enumerate}
\item
$\phi$ coïncide avec l'inclusion sur $\partial S$\up{\,;}
\item
$\phi(S) \subset U = S \times \R$ est une surface transversale au
champ de vecteurs $\partial_t$, $t \in \R$\up{\,;}
\item
$\phi(C)$ est legendrienne et son nombre de Thurston-Bennequin $\tb
(\phi(C), \phi(S))$ vaut $-\frac12 \Card (C \cap \Gamma)$.
\end{enumerate}
\end{lemme}

\subsection{Rocades legendriennes}

Dans l'exploration d'une variété de contact $(V,\xi)$, les
\emph{rocades} (\emph {bypasses})~\cite{Ho1} servent à pousser
l'investigation géométrique au-delà des domaines homogènes et à
analyser les changements (chirurgies) élémentaires que subit le
découpage d'une surface $\xi$-convexe quand on déforme celle-ci par
une grande isotopie. Elles jouent ainsi, en topologie de contact, un
rôle comparable à celui des disques de Whitney en topologie
différentielle.

\smallskip

Une \emph{rocade \up(legendrienne\up)} est un demi-disque $H$ plongé
dans $(V,\xi)$ et ayant les propriétés suivantes~:
\begin{itemize}
\item
le bord de $H$ est legendrien et $H$ est $\xi$-convexe --~en
particulier, les plans $T_pH$ et $\xi_p$ coïncident en chaque coin
$p$ de~$H$\,;
\item
l'enroulement de $\xi$ autour de $TH$ vaut $-1$ le long d'un des
deux arcs de $\partial H$, arc qu'on note $\arcinf H$, et $0$ le
long de l'autre, qu'on note $\arcsup H$ --~l'invariant de
Thurston-Bennequin de $\partial H$ est donc égal à~$-1$ et le
découpage de $\xi H$ se réduit à un arc dont les deux extrémités se
trouvent sur $\arcinf H$.
\end{itemize}
Une rocade matérialise ainsi une isotopie entre deux arcs
legendriens reliant les mêmes points, $\arcinf H$ et $\arcsup H$, le
second arc étant \guil{plus court} que le premier au sens où la
valeur absolue de l'enroulement de~$\xi$ y est moindre.

\medskip

Le lemme ci-dessous est un cas particulier du lemme de réalisation
legendrienne \ref{l:courbes}. Il permet de dénicher une rocade à
chaque fois qu'une surface $\xi$-convexe $S$ présente, dans son
découpage $\Gamma$, un arc propre \emph {parallèle au bord} (on dit
aussi \emph{$\partial$-parallèle}), c'est-à-dire bordant avec un arc
de $\partial S$ un demi-disque $H_0$ dont l'intérieur est disjoint
de~$\Gamma$.

\begin{lemme} \label{l:rocade1}
Soit $S$ une surface $\xi$-convexe à bord legendrien lisse, $U$ un
épaississement homogène de~$S$ et $\Gamma$ le découpage de~$\xi S$
associé. On suppose d'une part que, si $S$ est un disque, $\Gamma$
n'est pas connexe et d'autre part que l'adhérence d'une des
composantes connexes de $S \setminus \Gamma$ est un demi-disque $D$
dont un des arcs du bord est dans $\partial S$ et l'autre dans
$\Gamma$. On note $H$ un demi-disque voisinage de $D$ dans $S$ dont
un arc du bord est dans $\partial S$ et dont l'intersection avec
$\Gamma$ se réduit à l'arc $D \cap \Gamma$. Il existe alors un
plongement $\phi \from S \to U$ ayant les propriétés suivantes~:
\begin{enumerate}
\item
$\phi$ coïncide avec l'inclusion sur $\partial S$\up{\,;}
\item
$\phi(S) \subset U = S \times \R$ est une surface transversale au
champ de vecteurs $\partial_t$, $t \in \R$\up{\,;}
\item
$\phi(\Gamma)$ est l'intersection de $\phi(S)$ avec le cylindre
$\Gamma \times \R \subset U$\up{\,;}
\item
$\phi(H)$ est une rocade legendrienne.
\end{enumerate}
\end{lemme}

\medskip
On décrit maintenant l'attachement d'une rocade sur une surface
$\xi$-convexe et son effet sur le découpage de celle-ci~\cite{Ho1}.

Soit $S \subset V$ une surface $\xi$-convexe et $H$ une rocade
\emph{attachée} à $S$, c'est-à-dire vérifiant les propriétés
suivantes~:
\begin{itemize}
\item
$H$ s'appuie transversalement sur $S$ le long de $\arcinf H$ et ne
touche~$S$ nulle part ailleurs\,;
\item
il existe un découpage $\Gamma$ de $S$ qui passe par les extrémités
de $\arcinf H$ (et automatiquement par le milieu de $\arcinf H$).
\end{itemize}
On oriente $S$ de fa\c{c}on que $H$ soit du côté positif et on note
$N_\eps$ le $\eps$-voisinage de~$H$ privé de ce qui déborde du côté
négatif de $S$. Alors, pour $\eps>0$ assez petit, il existe sur $V$
un champ de vecteurs de contact qui est transversal à la surface
(lisse par morceaux)
$$ \Adh \bigl( (S \cup \partial N_\eps) \setminus (S \cap N_\eps) \bigr)  $$
et rentrant dans le quadrant extérieur à $N_\epsilon$
\cite{Ho1}. On peut donc lisser cette surface transversalement audit
champ de vecteurs pour obtenir une surface $\xi$-convexe $S'$. Il
existe clairement une isotopie de $S$ à~$S'$ et, si on l'utilise
pour identifier l'une à l'autre ces surfaces, on peut décrire comme
suit (à isotopie près) le découpage $\Gamma'$ de $S'$ en fonction du
découpage $\Gamma$ de $S$ et de la base $\arcinf H$ de $H$. On prend
sur $S$ un carré $R = [-1,1] \times [-1,1]$ (les coordonnées étant
compatibles avec l'orientation) qui contient $\arcinf H = [-1,1]
\times \{0\}$ et rencontre $\Gamma$ selon les trois segments $\{-1,
0, 1\} \times [-1,1]$. La multi-courbe $\Gamma'$ s'obtient à partir
de $\Gamma$ en retirant les segments
$$ \{-1\} \times [0,1], \quad \{0\} \times [-1,1] \quad \text{et} \quad
   \{1\} \times [-1,0] \phantom{.} $$
et en les remplaçant par
$$ [-1,0] \times \{1\}, \quad [-1,1] \times \{0\} \quad \text{et} \quad
   [0,1] \times \{-1\} . $$
On nomme cette opération \emph{chirurgie de $\Gamma$ le long de
$\arcinf H$}. Dit autrement, si on regarde plutôt le complémentaire
du découpage comme un coloriage de la surface en deux couleurs, le
passage de $S$ à $S'$ consiste à faire tourner, dans le carré $R$,
le coloriage d'un angle $\pi/2$ sans le changer en dehors.

\subsection{Modification de Lutz}

La modification de Lutz, définie ci-après, est une opération qui
consiste à changer une structure de contact au voisinage d'un tore
transversal.

\smallskip

Soit $\xi$ une structure de contact sur~$V$ et $T$ un tore plongé
dans~$V$ transversalement à~$\xi$. On note
$$ W = T \times [-1,1] \supset T = T \times \{0\} $$
un voisinage tubulaire compact de $T$ dont les fibres sont des arcs
legendriens dans $(V, \xi)$ et on paramètre~$T$ par $\R^2 \!/ \Z^2$.
Ainsi, la structure de contact $\xi$ admet dans~$W$ une équation du
type
$$ \cos \theta(x,t) \, dx_1 - \sin \theta(x,t) \, dx_2 = 0,
   \qquad (x,t) \in \R^2 \!/ \Z^2 \times [-1,1], $$
où la fonction $\theta \from \R^2 \!/ \Z^2 \times [-1,1] \to \S^1$
vérifie $\partial_t \theta > 0$ en tout point.

Pour tout entier $n \ge 0$, soit $\rho_n \from [-1,1] \to \S^1 = \R
/ 2\pi\Z$  la projection d'une fonction (faiblement) croissante
$[-1,1] \to \R$ qui vaut~$0$ près de~$-1$ et~$n\pi$ près de~$1$. La
structure de contact définie sur $W$ par l'équation
$$ \cos \bigl( \theta(x,t) + \rho_n(t) \bigr) \, dx_1
 - \sin \bigl( \theta(x,t) + \rho_n(t) \bigr) \, dx_2 = 0,
   \qquad (x,t) \in \R^2 \!/ \Z^2 \times [-1,1], $$
coïncide avec $\xi$ près de $T \times \{-1\}$ et avec $(-1)^n \xi$
près de $T \times \{1\}$. On peut donc la prolonger par $\xi$ en
dehors de~$W$ pour obtenir une structure de contact $\xi_n$ sur~$V$
mais celle-ci n'est orientable que si $n$ est pair ou si $T$
sépare~$V$. Par convention, on prend alors sur $\xi_n$ l'orientation
qui épouse celle de~$\xi$ près de $T \times \{1\}$. On dit que
$\xi_n$ est obtenue à partir de~$\xi$ par une modification de Lutz
de coefficient~$\frac{n}{2}$ le long de~$T$.

Des calculs simples montrent que la classe d'isotopie de~$\xi_n$ ne
dépend que de~$n$ (et pas du choix de $W$ ou de $\rho_n$) tandis que
sa classe d'homotopie dans l'espace des champs de plans tangents est
seulement fonction de la parité de~$n$ et est notamment égale à
celle de~$\xi$ quand $n$ est pair.

Les deux exemples suivants permettent d'éclairer le rôle joué
par les modifications de Lutz.

\begin{itemize}
\item Une modification de Lutz le long d'un tore compressible dont le
feuilletage caractéristique est une suspension (ou plus
généralement pour lequel  le champ de droites tracé par $\xi$
sur $T$ est homotope à un champ constant dans la trivialisation
habituelle de $T \cong \R^2 \!/ \Z^2$) donne toujours une structure
de contact vrillée.
\item Soit $T$ un tore dont le feuilletage caractéristique $\xi T$
est la réunion de deux composantes de Reeb qui pointent dans la
même direction (en particulier le champ de ses tangentes n'est pas
homotope à un champ de droites constant). Une modification de Lutz
de coefficient $1$ le long de $T$ donne une structure de contact
conjuguée à la structure initiale $\xi$ par un twist de Dehn, et
donc isotope à $\xi$ quand le tore est compressible. Pour s'en
convaincre, il suffit de constater que faire pivoter les tangentes
à $\xi T$ d'un angle $\theta$ conduit à un champ de droites dont
le feuilletage intégral se déduit de $\xi T$ par une isotopie
qui déplace les orbites périodiques. Lorsque $\theta$ varie
entre $0$ et $2\pi$, les positions successives prises par une orbite
périodique $\gamma$ de $\xi T$ forment un feuilletage de  $T$.
\end{itemize}

\subsection{Surfaces branchées}

Soit $G_0, G_1, G_2$ les graphes respectifs dans $\R^3 = \R^2 \times \R$ des 
fonctions $f_0, f_1, f_2$ définies sur $\R^2$ par
$$ f_0(x,y) = 0, \quad
   f_1(x,y) = \begin{cases}
   0 & \text{si $x \ge 0$,} \\
   e^{1/x} & \text{si $x < 0$,} \end{cases} \quad \text{et} \quad
   f_2(x,y) = - f_1(y,x) . $$
On note $\X$ l'intersection de $G_0 \cup G_1 \cup G_2$ avec le demi-espace 
$\{ x \ge -1\}$. Un homéomorphisme d'un ouvert $U$ de $\X$ sur un autre est dit 
\emph{lisse} si sa restriction à chaque intersection $U \cap G_i$, $0 \le i \le 
2$, est lisse comme application à valeurs dans $\R^3$.

Une \emph{surface branchée} est un espace topologique~$X$ séparé, de type 
dénombrable et localement modelé sur $\X$ par des cartes --~c'est-à-dire dont 
chaque point possède un voisinage homéomorphe à un ouvert de~$\X$~-- avec des 
changements de cartes lisses. On nomme
\begin{itemize}
\item
\emph{bord} de $X$ l'ensemble $\partial X$ des points que les cartes appliquent 
dans l'intersection $\X \cap \{ x=-1 \}$ ---~cet ensemble est un graphe doté 
d'une tangente en tout point, autrement dit un \emph{réseau ferroviaire}\,;
\item
\emph{point triple} tout point de~$X$ que les cartes envoient sur $(0,0,0)$\,;
\item
\emph{point double} tout point de~$X$ que les cartes appliquent dans $\X \cap 
\{ xy=0 \}$ sur un point autre que $(0,0,0)$\,;
\item
\emph{point régulier} tout point de $X$ qui n'est ni double ni triple.
\end{itemize}
Cette définition engendre une notion naturelle d'applications lisses d'une 
surface branchée dans une autre ou entre une surface branchée et une variété, 
d'où (par exemple) une notion de courbes immergées dans une surface branchée et 
de surfaces branchées plongées dans une variété.

Les sous-ensembles $X_0$, $X_1$ et $X_2$ constitués des points respectivement 
réguliers, doubles et triples forment une stratification de~$X$. L'ensemble
$X_0$, qu'on note aussi $\Reg (X)$, est appelé \emph{partie régulière} de~$X$ et 
son complémentaire $\Theta = X_1 \cup X_2$ le \emph{lieu singulier}. Ce lieu est
une multi-courbe immergée à points doubles ordinaires (les points triples de
$X$). L'ensemble $X_1$ est de plus naturellement coorienté~: par convention, son
côté positif est celui où $X$ a un seul feuillet. La direction donnée par cette 
coorientation est nommée \emph{direction de branchement}. Enfin, les composantes 
connexes de l'ensemble $X_0 = \Reg (X)$ sont appelées strates régulières ou \emph
{secteurs} de~$\Reg (X)$.

\subsection{Domaines fibrés}

Les surfaces branchées apparaissent souvent, de façon naturelle,
comme des quotients de variétés --~à bord anguleux~-- par certains
feuilletages de codimension $2$.

\begin{figure} [ht]
    {\epsfysize=2in\centerline{\epsfbox{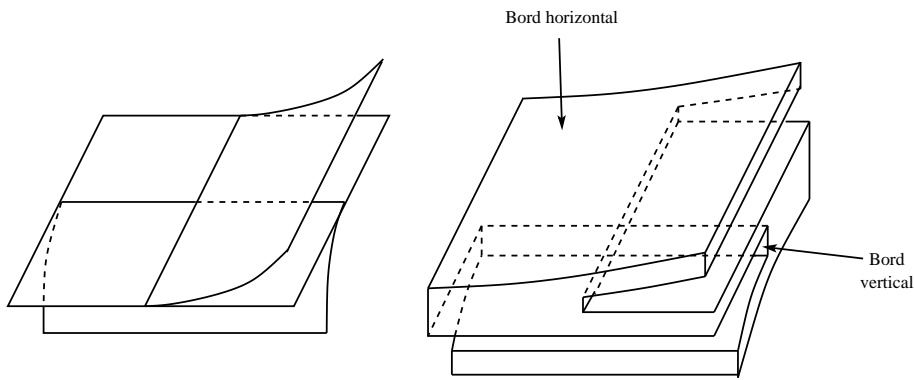}}}
    \caption{Surface branchée et domaine fibré.}
    \label{fig8}
\end{figure}

A titre de premier exemple, on considère trois copies $P_0$, $P_1$
et $P_2$ du plan $\R^2$ et on met sur leur union la relation
d'équivalence engendrée par les règles suivantes~:
\begin{itemize}
\item
un point $(x_0,y_0) \in P_0$ est équivalent à un point $(x_1,y_1)
\in P_1$ si et seulement si $x_0 = x_1$, $y_0 = y_1$ et $x_0 \le
0$\,;
\item
un point $(x_0,y_0) \in P_0$ est équivalent à un point $(x_2,y_2)
\in P_2$ si et seulement si $x_0 = x_2$, $y_0 = y_2$ et $y_0 \le 0$.
\end{itemize}
Le quotient de $P_0 \cup P_1 \cup P_2$ par cette relation est une
surface branchée difféomorphe à l'union $G_0 \cup G_1 \cup G_2$ des
trois graphes décrits au début de la partie~précédente (voir la
figure~\ref{fig8}).

Soit $M$ une variété compacte, orientable, à bord anguleux --~en
fait, ayant des arêtes lisses mais pas de coins~-- et munie d'un
feuilletage~$\tau$ en intervalles compacts. On dit que $M$ est un
\emph{domaine fibré} si $\partial M$ est l'union de deux surfaces
compactes lisses $\partial_h M$ et $\partial_v M$ satisfaisant aux
conditions suivantes~:
\begin{itemize}
\item
$\partial_h M$ est transversale au feuilletage $\tau$ tandis que
$\partial_v M$ est tangente à $\tau$\,;
\item
les surfaces $\partial_h M$ et $\partial_v M$ ont le même bord
--~qui est la partie anguleuse de $\partial M$~-- et chaque
composante $S$ de $\partial_v M$ a sur son bord des points où $M$
est concave (c'est-à-dire modelé non pas sur un quart mais sur trois
quarts d'espace)\,;
\item chaque feuille de $\tau$ traverse au plus deux fois $\partial_v M$.
\end{itemize}

\begin{proposition}
Le quotient $X = M / \tau$ est  naturellement une surface branchée
sans bord dont le lieu singulier est l'image des composantes $S$ de
$\partial_v M$ telles que $M$ soit concave près de chaque point de
$\partial S$.
\end{proposition}

\begin{proof}
Soit $\pi$ la projection $M \to X = M/\tau$. Pour tout point $p \in
X$ situé hors de $\pi (\partial_v M)$, on prend comme carte locale
une transversale à $\tau$. Si $p \in \pi (\partial_v M )$ est dans
l'image d'une seule composante, on a deux cartes locales naturelles
qui consistent à prendre des transversales à $\tau$ le long d'une
composante du bord ou le long de l'autre.
\end{proof}

Lorsque les deux premières conditions sont réalisées, on peut
toujours obtenir la troisième, quitte à perturber $\tau$ dans
l'intérieur de $M$.

On observe que chaque composante connexe de $\partial_v M$ est
feuilletée par des intervalles compacts transversaux au bord et est
donc un anneau.

Dès que $X$ a des points triples, il n'existe aucun plongement $\phi
\from X \to M$ qui soit une section de la projection $\pi \from M
\to X = M/\tau$. On peut cependant trouver des plongements $\phi_s
\from X \to M$, $s \in \mathopen] 0, 1]$, transversaux à~$\tau$ qui
sont des sections de $\pi$ en dehors du $s$-voisinage du lieu
singulier et ont la propriété que $\pi \circ \phi_s$ converge vers
l'identité uniformément quand $s$ tend vers~$0$. On dit alors
souvent, pour $s$ petit, que $M$ est un \guil{voisinage fibré} de
$\phi_s(X)$.

Une surface $S\subset (M,\tau )$ est \emph{portée} par $(M,\tau)$
si elle est transversale à $\tau$. Elle est \emph{pleinement
portée} lorsque de surcroît elle rencontre toutes les fibres
de $\tau$.

\section{Structures de contact et triangulations}

Le but de cette partie est la démonstration du résultat suivant~:

\begin{theoreme} \label{t:domaines}
Soit $V$ une variété close et orientée de dimension $3$. Il existe dans~$V$ un 
nombre fini de domaines fibrés $(M_i,\tau_i)$ flanqués chacun d'une structure de
contact~$\zeta_i$ sur $V \setminus \Int (M_i )$ tels que toute structure de contact
tendue sur~$V$, à isotopie près et pour un certain~$i$, soit égale à~$\zeta_i$ 
hors de~$M_i$ et tangente à~$\tau_i$ dans~$M_i$.
\end{theoreme}

\subsection{Triangulations de contact}

Soit $V$ une variété de dimension~$3$. Ce qu'on appelle \emph{triangulation} de
$V$, dans la suite, est une décomposition simpliciale possédant les propriétés
de régularité suivantes~:
\begin{itemize}
\item
chaque simplexe de dimension~$2$ (ou moins) est lisse\,;
\item
chaque simplexe de dimension~$3$ a, le long de ses arêtes mais en dehors des 
sommets, un angle d'ouverture strictement compris entre $0$ et~$\pi$.
\end{itemize}
Si $\Delta$ est une triangulation de $V$, on note $\Delta^i$, $0 \le i \le 3$,
son squelette de dimension~$i$. 
En outre, une \emph{isotopie de triangulations} est un chemin de triangulations.
Une telle isotopie $\Delta_t$ se prolonge en une isotopie ambiante $\psi_t$ de 
$V$ qui n'est bien sûr pas lisse en général mais se compose d'homéomorphismes 
dont les restrictions à chaque $2$-simplexe de~$\Delta_0$ sont lisses. On écrira
$\Delta_t = \psi_t(\Delta_0)$ pour dire que $\Delta_t$ est l'image de~$\Delta_0$
par~$\psi_t$.

Toute triangulation lisse --~c'est-à-dire dont chaque simplexe est lisse~\cite
{Wh}~-- est évidemment une triangulation au sens ci-dessus mais on aura besoin 
ici de triangulations qui ne sont pas lisses.

\begin{definition}[\cite{Gi6}] \label{d:triangulation}
Soit $(V,\xi)$ une variété de contact de dimension~$3$. Une \emph{triangulation 
de contact} de $(V,\xi)$ est une triangulation de~$V$ vérifiant les conditions 
suivantes (voir la définition \ref{d:convexiterelative} pour la condition de 
convexité relative)~:
\begin{enumerate}
\item
les $1$-simplexes sont des arcs legendriens\,;
\item
les $2$-simplexes sont $\xi$-convexes et relativement
$\xi$-convexes\,;
\item
les $3$-simplexes sont contenus dans des cartes de Darboux\footnote
{Cette condition est automatiquement remplie quand la structure de contact~$\xi$
est tendue puisque, d'après~\cite{El2}, toute boule munie d'une structure de 
contact tendue se plonge dans l'espace de contact standard.}.
\end{enumerate}
Une telle triangulation n'est pas lisse car les arêtes issues d'un sommet~$p$ 
quelconque sont toutes tangentes au même plan, à savoir le plan de contact
$\xi_p$.

\smallskip

On appelle \emph{nombre de Thurston-Bennequin} d'une triangulation de contact 
$\Delta$ de $(V,\xi)$ l'entier
$$ \TB(\Delta) = \TB(\Delta,\xi) = - \sum_F \tb(\partial F), $$
la sommation se faisant sur tous les $2$-simplexes~$F$. Comme $\tb (\partial F) 
\le -1$ pour tout~$F$ (d'après l'inégalité de Bennequin), l'entier $\TB(\Delta)$
est positif et au moins égal au nombre de $2$-simplexes de $\Delta$.

Une triangulation de contact $\Delta$ de $(V,\xi)$ sera dite \emph{minimale} si 
elle a le plus petit nombre de Thurston-Bennequin parmi les triangulations de 
contact obtenues en déformant $\Delta$ par une isotopie à support compact et 
relative à un voisinage des sommets.
\end{definition}

On explicite maintenant la condition technique de $\xi$-convexité relative qui, 
en pratique, est aussi facile à réaliser que la $\xi$-convexité ordinaire~:

\begin{definition} \label{d:convexiterelative}
Soit $(V,\xi)$ une variété de contact et $\Delta$ une triangulation de $V$ dont 
le $1$-squelette est legendrien. Soit $G$ un $3$-simplexe, $F_0, F_1$ deux faces
de~$G$ et $a$ leur arête commune. On oriente $F_1$ comme partie de $\partial G$ 
et $a$ comme partie de $\partial F_1$. On dira que $F_0$ et~$F_1$ sont \emph
{$\xi$-indisciplinées} le long d'un arc (orienté) $[p,q] \subset a$ si $\xi$~est
tangent à $F_0$ en $p$, à $F_1$ en $q$ et n'est un hyperplan d'appui de $G$ en
aucun point de $]p,q[$\footnote
{Un tel arc $]p,q[$ est, en quelque sorte, une orbite de $\xi \, \partial G$ qui
va d'une singularité négative à une singularité positive (voir la remarque \ref
{r:demitour}).}.
On dit que les $2$-simplexes de $\Delta$ sont \emph{relativement $\xi$-convexes}
si aucun $3$-simplexe ne possède de faces $\xi$-indisciplinées le long d'un arc 
de leur arête commune.
\end{definition}

\medskip

On rappelle (cf. section~\ref{subsection:tendues}) qu'un ensemble
$\XX$ de structures de contact tendues sur une variété $V$ close et
orientée est dit complet s'il représente toutes les classes
d'isotopie dans $\SCT(V)$.  Le point de départ de cette étude est la
proposition suivante, variante d'un résultat de \cite{Gi6} servant à
construire des livres ouverts adaptés aux structures de contact~:

\begin{proposition} \label{p:triangulation}
Sur toute variété close et orientée $V$ de dimension~$3$, il existe
un ensemble complet~$\XX$ de structures de contact tendues et une
triangulation~$\Delta$ ayant les propriétés suivantes\up{~:}
\begin{enumerate}
\item
toutes les structures de contact de $\XX$ coïncident --~comme champs de plans
non orientés~-- sur un voisinage des sommets de~$\Delta$\up{\,;}
\item
$\Delta$ est une triangulation de contact minimale de $(V,\xi)$ pour tout $\xi
\in \XX$.
\end{enumerate}
De plus, on peut choisir $\Delta$ isotope à n'importe quelle
triangulation lisse fixée de~$V$.
\end{proposition}

Avant de démontrer cette proposition, on traite le cas d'une seule
structure de contact~:

\begin{lemme} \label{l:triangulation}
Soit $(V,\xi)$ une variété de contact de dimension~$3$. Si la
structure $\xi$ est tendue, toute triangulation lisse de $V$ est
isotope à une triangulation de contact de $(V,\xi)$.
\end{lemme}

\begin{proof}
Soit $\Delta_0$ une triangulation lisse de~$V$. La possibilité de
déformer $\Delta_0$ par isotopie en une triangulation de
contact~$\Delta_1$ de $(V,\xi)$ tient avant tout au fait qu'on peut
déformer chaque arête de $\Delta_0$ en un arc legendrien par une
isotopie relative à ses extrémités --~et arbitrairement
$\classe0$-petite. Il faut toutefois un peu de soin aux sommets pour
garantir la différentiabilité du prolongement aux $2$-simplexes.

En chaque sommet~$s$, on se donne une projection $d\pi_s \from T_sV
\to \xi_s$ dont le noyau est transversal aux plans tangents des
$2$-simplexes en~$s$ et qui envoie les tangentes aux arêtes sur des
droites distinctes. Sur un voisinage compact convenable~$N(s)$, il
existe des coordonnées $(x,y,z) \in \D^2 \times [-1,1]$ centrées
en~$s$ dans lesquelles $d\pi_s$ est la dérivée de la projection
$\pi_s \from (x,y,z) \mapsto (x,y)$ et $\xi$ a pour équation $dz +
x\,dy - y\,dx = 0$. Par ailleurs, pour $N(s)$ assez petit, $\pi_s$
induit un plongement sur le $1$-squelette et sur chaque $2$-simplexe
de~$\Delta_0$. On déforme $\Delta_0$ dans~$N(s)$ comme suit (voir
l'exemple~\ref{x:modele1})~:
\begin{itemize}
\item
les arêtes de $\Delta_1$ sont les arcs legendriens sur lesquels les
arêtes de $\Delta_0$ se projettent\,;
\item
les faces de $\Delta_1$ sont des graphes au-dessus ou en-dessous de
celles de $\Delta_0$\,;
\item
les arêtes et les faces des triangulations intermédiaires
$\Delta_t$, pour tout $t \in [0,1]$, sont obtenues en prenant
l'isotopie barycentrique verticale au temps $t$ entre les arêtes et
les faces de~$\Delta_0$ et celles de~$\Delta_1$.
\end{itemize}
Moyennant quelques petites précautions dans la mise en place des
$2$-simplexes de $\Delta_1$, les angles diédraux de ses
$3$-simplexes dans~$N(s)$ sont compris strictement entre $0$
et~$\pi$ en dehors de~$s$ et il en va alors de même pour les angles
diédraux des $3$-simplexes de chaque $\Delta_t$, $t \in [0,1]$. On a
ainsi fabriqué les triangulations~$\Delta_t$ voulues près des
sommets.

On choisit ensuite une isotopie $\Delta_t^1$ du $1$-squelette de
$\Delta_0$ qui prolonge celle construite dans les cartes $N(s)$, est
lisse en dehors et aboutit à un graphe legendrien $\Delta_1^1$. Dans
$V' = V \setminus \bigcup \Int (N(s))$, cette isotopie s'étend en une
isotopie ambiante lisse $\psi_t$ engendrée par un champ qui, sur le
bord latéral $\partial_| N(s) = (\partial \D^2 )\times [-1,1]$ de
chacun des voisinages $N(s) = \D^2 \times [-1,1]$, est un multiple
du champ $\partial_z$. Les images $\psi_t (\Delta_0^2 \cap V')$ du
$2$-squelette ne se recollent pas \emph{a priori} aux $2$-simplexes
bâtis dans $N(s)$ mais les arcs qu'ils tracent sur les anneaux
$\partial_| N(s)$ sont des graphes au-dessus de ceux tracés par les
$2$-simplexes dans $N(s)$ (pour la projection $\pi_s$). Une
déformation de $\psi_t (\Delta_0^2 \cap V')$ près des anneaux
$\partial_| N(s)$ et laissant fixe le $1$-squelette $\Delta_t^1$
permet alors de le recoller aux $2$-simplexes de $\Delta_t$ dans
$N(s)$. On a ainsi modifié $\Delta_0$, par l'isotopie $\Delta_t$, $t
\in [0,1]$, en une triangulation $\Delta_1$ dont les $1$-simplexes
sont des arcs legendriens. De plus, comme $\xi$ est tendue, les
$3$-simplexes sont tous inclus dans des cartes de Darboux.

Quitte à décroître le nombre de Thurston-Bennequin du bord des
$2$-simplexes de $\Delta_1$ (par stabilisation), on peut, dans la
construction précédente, choisir l'isotopie de graphes $\Delta_t^1$
de sorte que, le long de chaque arête de~$\Delta_1$, l'enroulement
de~$\xi$ autour d'une face quelconque soit négatif ou nul. Ceci
permet de déformer $\Delta_1$, par une isotopie ambiante lisse
de~$V$ relative au $1$-squelette et à un voisinage des sommets, pour
que $\xi$ ne fasse aucun demi-tour inversé par rapport à une face le
long d'une arête de son bord. Par le même principe, on élimine tous
les arcs d'arêtes le long desquels deux faces adjacentes sont
$\xi$-indisciplinées. Les $2$-simplexes sont alors relativement
$\xi$-convexes et on les rend $\xi$-convexes par une ultime petite
perturbation de~$\Delta_1$ --~toujours stationnaire sur le
$1$-squelette et un voisinage des sommets. On obtient ainsi une
triangulation de contact $\Delta_2$ de $(V,\xi)$ isotope
à~$\Delta_0$. Enfin, puisque les triangulations de contact de
$(V,\xi)$ qui sont isotopes à $\Delta_2$ relativement à un voisinage
des sommets ont un nombre de Thurston-Bennequin minoré, l'une
d'elles est minimale.
\end{proof}

\begin{exemple} \label{x:modele1} Si une face $F_0$ de la triangulation
initiale $\Delta_0$ se projette sur le quadrant $Q = \{ x \ge 0,\; y
\ge 0 \}$ --~ce qui est toujours le cas quitte à prendre d'autres
coordonnées de Darboux~--, la face $F_1$ correspondante dans la
triangulation de contact $\Delta_1$ a une équation du type $z = xy
u(x,y)$ et son point complexe à l'origine est hyperbolique si $\abs{
u(0,0) } > 1$ et elliptique si $\abs{ u(0,0) } < 1$. On ne peut donc
pas choisir arbitrairement la nature des points complexes que les
faces de~$\Delta_1$ présentent en leurs sommets (si deux arêtes de
$\Delta_0$ issues de l'origine se projettent par exemple à
l'intérieur du quadrant $Q$ et se trouvent de part et d'autre de
$F_0$, le point complexe de $F_1$ ne peut être hyperbolique) mais on
peut faire en sorte qu'ils soient tous elliptiques.
\end{exemple}

\begin{proof}[Démonstration de la proposition~\ref{p:triangulation}]
Soit $\xi_0$ une structure de contact tendue et $\Delta_0$ une
triangulation de contact minimale de $(V,\xi_0)$
(lemme~\ref{l:triangulation}). Soit $N = \bigcup N(s)$ un voisinage
compact du $0$-squelette de $\Delta_0$ dont chaque composante
connexe est une carte de Darboux $N(s) = \D^2 \times [-1,1]$ autour
d'un sommet $s$. On prend les domaines $N(s)$ assez petits pour que
leur intersection avec chaque $2$-simplexe de $\Delta_0$ soit un
graphe au-dessus de sa projection sur $\D^2$ et évite $\D^2 \times
\{\pm1\}$. Toute structure de contact sur $V$ est isotope à une
structure~$\xi$ qui coïncide avec~$\xi_0$ sur~$N$. Pour peu que
$\xi$ soit tendue, le lemme \ref{l:triangulation} (ou plus
exactement sa preuve) montre que $(V,\xi)$ possède une triangulation
de contact (minimale) $\Delta$ isotope à $\Delta_0$ relativement à
$N$. Par une isotopie lisse $\psi_t$ de~$V$ stationnaire sur~$N$, on
déforme $\Delta_0$ en une triangulation $\Delta_1 = \psi_1
(\Delta_0)$ qui a le même $1$-squelette que $\Delta$ et des
$2$-simplexes tous $\xi$-convexes\footnote
{On utilise ici le fait que la $\xi$-convexité est une propriété à
la fois dense et ouverte.}
et relativement $\xi$-convexes (voir la preuve du lemme
\ref{l:triangulation}). Comme $\xi$ est tendue, $\Delta_1$ est une
triangulation de contact minimale de $(V,\xi)$. Par conséquent,
$\Delta_0$ est une triangulaton de contact minimale de $(V,
\psi_1^*\xi)$, ce qui établit la proposition.
\end{proof}

En pratique, les triangulations de contact minimales s'avèrent trop
rigides. Il est notamment difficile d'appliquer aux $2$-simplexes
les lemmes de réalisation de feuilletages caractéristiques sans
créer d'intersections indésirables entre eux près des sommets. Pour
cette raison, on introduit maintenant une notion de triangulations
de contact \guil{maniables} pour laquelle on donne une variante de
la proposition \ref{p:triangulation}.

\begin{definition} \label{d:maniable}
Soit $(V, \xi)$ une variété de contact de dimension $3$. Une
triangulation de contact $\Delta$ de $(V,\xi)$ est dite
\emph{maniable} si chaque simplexe $F$ de dimension $2$ écorné de
trois triangles aux sommets est un hexagone $H_F$ à bord legendrien
et $\xi$-convexe, et si, pour toute arête $a$, l'enroulement de
$\xi$ autour de $F$ le long de l'arc $r(a)=a \setminus \Int (\Lambda )
$, $\Lambda = \bigcup_F \overline{(F \setminus H_F)}$, est
strictement négatif. (Noter que cet arc n'est, en général, qu'un
morceau d'arête de $H_F$ mais pas une arête entière.) Dans le
$2$-squelette de $\Delta$, l'ensemble $\Lambda = \bigcup_F
\overline{(F \setminus H_F)}$ est un voisinage des sommets qu'on
nomme \emph{voisinage de sécurité}. On dit que $\Delta$ est {\it
$\Lambda$-minimale} si elle a le plus petit nombre de
Thurston-Bennequin parmi toutes les triangulations de contact
maniables de $(V,\xi)$ ayant $\Lambda$ pour voisinage de sécurité.
Dit autrement en faisant porter la déformation sur $\xi$, la
triangulation $\Delta$ est $\Lambda$-minimale si pour toute
structure $\xi'$ isotope à $\xi$ relativement à $\Lambda$ et
pour laquelle $\Delta$ est une triangulation de contact maniable,
$\TB (\Delta ,\xi' )\geq \TB (\Delta ,\xi )$. On note $\Lambda_F
(s)$ le triangle de $\Lambda$ inclus dans la face $F$ et contenant
le sommet $s$, et $\Lambda (s)$ la réunion des triangles de
$\Lambda$ contenant $s$.

\end{definition}

\begin{proposition} \label{p:maniable}
Sur toute variété close et orientée $V$ de dimension~$3$, il existe
un ensemble complet~$\XX$ de structures de contact tendues et une
triangulation~$\Delta$ ayant les propriétés suivantes\up{~:}
\begin{itemize}
\item
toutes les structures de contact de $\XX$ coïncident --~comme champs
de plans non orientés~-- sur un voisinage $U$ des sommets
de~$\Delta$\up{\,;}
\item
$\Delta$ est, pour tout $\xi \in \XX$, une triangulation de contact
maniable de $(V,\xi)$ dont le voisinage de sécurité $\Lambda$ est
fixe et contenu dans $U$ et qui est $\Lambda$-minimale.
\end{itemize}
\end{proposition}

\begin{proof}[Démonstration de la proposition~\ref{p:maniable}]
Soit $\xi_0$ une structure de contact tendue et $\Delta$ une
triangulation de contact de $(V,\xi_0)$ (le
lemme~\ref{l:triangulation} montre qu'il en existe). On reprend les
notations de la preuve du lemme~\ref{l:triangulation}. Dans chaque
petite boule $N(s)$ centrée en un sommet $s$ de $\Delta$, on
décroît localement, par une isotopie de stabilisation, le
nombre de Thurston-Bennequin de chaque arête $a$ issue de $s$ afin
que, dans $N(s)$, celui-ci soit strictement négatif et que $\xi_0$
imprime dans $N(s)$ un point singulier sur toute face adjacente à
$a$. Dans cette situation, on peut tracer sur chaque face $F$  de
sommet $s$ un arc $\gamma_{F,s}$ inclus dans $F \cap N(s)$, tangent
à $\xi$ en ses extrémités et délimitant un triangle
$\Lambda_F (s)$ dans $F$. Par approximation legendrienne, on est
alors en mesure d'effectuer une isotopie de $\xi_0$ à support dans
un voisinage de $\bigcup_{F,s} \gamma_{F,s}$ et relative au
$1$-squelette de $\Delta$ pour que $\gamma_{F,s}$ devienne un arc
legendrien d'invariant de Thurston-Bennequin relatif à $F$
inférieur à $-1$.  La réunion des arcs $\gamma_{F,s}$
découpe les arêtes de $\Delta$ en un certain nombre d'arcs
legendriens. Quitte à stabiliser chacun de ces arcs, on se
ramène au cas où, pour toute face $F$, ils ont tous un nombre de
Thurston-Bennequin relatif à $F$ inférieur à $-1$. Une
isotopie à support dans un voisinage du $1$-squelette permet
ensuite d'éviter tous les demi-tours inversés le long des
arêtes ainsi que des arcs $\gamma_{F,s}$. On rend alors les faces
$\xi_0$-convexes et relativement $\xi_0$-convexes à l'aide d'une
isotopie relative à $\Delta^1$ et aux arcs $\gamma_{F,s}$. La
triangulation $\Delta$ est de contact et maniable pour $(V,\xi_0)$
relativement au voisinage de sécurité $\Lambda =\bigcup_{F,s}
\Lambda_F (s)$.

La proposition~\ref{p:triangulation} permet maintenant d'isotoper
toute structure $\xi$ sur une structure $\xi'$ qui coïncide
avec $\xi_0$ le long de $\Lambda$ et qui a $\Delta$
pour triangulation de contact. Le procédé de stabilisation
fournit une isotopie
de $\xi'$ relative à $\Lambda$ qui donne aux trois côtés de
l'hexagone $H_F$ inclus dans $\partial F$ des nombres de
Thurston-Bennequin relatifs inférieurs à $-1$. On fait même en
sorte que les arcs $r(a) =\overline{a\setminus \Lambda} \subset a$ aient
un invariant de Thurston-Bennequin relatif à toute face adjacente
inférieur à $-\frac{1}{2}$. La triangulation $\Delta$ est donc
maniable pour $\xi'$ et a $\Lambda$ comme voisinage sûr.
Comme dans le lemme~\ref{l:triangulation},
on peut également la supposer $\Lambda$-minimale.
L'ensemble $\XX$ est constitué du choix d'une structure $\Lambda$-minimale
pour $\Delta$ dans chaque classe d'isotopie de structures tendues.
\end{proof}

Dans toute la suite du texte, on d\'esigne par
$\XX$ un syst\`eme complet de structures de contact tendues sur $V$ 
et par $\Delta$ une triangulation de $V$
tels que $\Delta$ soit une triangulation de contact maniable et
$\Lambda$-minimale pour toute structure $\xi \in \XX$, o\`u $\Lambda$
est un voisinage de s\'ecurit\'e fixe le long duquel toutes
les structures de $\XX$ co\"\i ncident.

\subsection{Propriétés des triangulations de contact minimales}

On commence par une propriété valable pour toutes les triangulations
de contact.

\begin{lemme}\label{lemme : cercle} Si $\Delta$ est une triangulation de contact pour
$(V,\xi )$, pour toute face $F\in \Delta^2$, aucune composante de la
courbe de découpage $\Gamma_F (\xi )$ n'est un cercle.
\end{lemme}
\begin{proof}[Démonstration.] Cette absence de courbes fermées vaut en fait
pour le découpage de toute surface $\xi$-convexe autre qu'une sphère
et au voisinage de laquelle $\xi$ est tendue. Or chaque $3$-simplexe
de $\Delta$ est inclus dans une carte de Darboux donc, \emph{a
fortiori}, $\xi$ est tendue au voisinage de $F$. (Si une composante
de $\Gamma_F$ était un cercle, le lemme de réalisation de
feuilletages permettrait de faire apparaître un disque
vrillé.)
\end{proof}

Soit $(V,\xi)$ une variété de contact tendue de dimension $3$ munie d'une
triangulation de contact $\Delta$ maniable et $\Lambda$-minimale, où $\Lambda$
est son voisinage de sécurité.
Pour une arête $a$ de $\Delta$, on rappelle que $r(a) =\overline{a\setminus \Lambda}$
et que $H_F$ est l'hexagone $F\setminus (\bigcup_{1\leq i\leq 3} \Int (\Lambda_F
(s_i )))$.

Soit $F$ une face de $\Delta$. On dira qu'une composante $\Gamma$ de
$\Gamma_{H_F} (\xi )$ est \emph{extrémale} si un des deux points
$p$ de $\partial \Gamma$ appartient à un arc $r(a)$, où $a$ est
une arête de $F$, ne peut pas être poussé en dehors de $r(a)$
par une isotopie de $\Gamma$ parmi les courbes transversales à
$\xi F$ (autrement dit il y a des singularités de $\xi F$ entre
$p$ et  les extrémités de $r(a)$) et est le plus proche d'une
des extrémités de  $r(a)$ parmi les points de $\Gamma_{H_F} (\xi
) \cap r(a)$ qui possèdent ces propriétés. En particulier,
pour un certain choix de $\Gamma_{H_F} (\xi )$, un point de $\Gamma
\cap r(a)$ est effectivement extrême. La multicourbe $\Gamma_{H_F}
(\xi )$ contient au plus six composantes extrémales.

\begin{lemme} \label{l:extremale}
Si $\xi$ est une structure de contact tendue sur $V$ et si $\Delta$ est une triangulation
de contact maniable et $\Lambda$-minimale pour $\xi$, 
alors toute composante du découpage $\Gamma_{H_F}(\xi)$ parallèle à 
$r(a)$ est extrémale.
\end{lemme}

\begin{proof}
On suppose qu'une composante de $\Gamma_{H_F} (\xi )$
parallèle à $r(a)$ n'est pas extrémale. Ceci signifie
en particulier que l'enroulement de $\xi$ le long de $r(a)$
relativement à $F$ est inférieur à $-2$.
Le lemme de réalisation~\ref{l:rocade1} permet,
par une isotopie de $H_F$ relative à un voisinage de
$\partial H_F \setminus r(a)$ et de support inclus dans
un petit voisinage homogène de $H_F$, de déformer
$H_F$ en une surface $H_F'$ qui contient une rocade $D$ s'appuyant
sur $r(a)$ (c'est-à-dire que $\partial_- D \subset r(a)$ et
$\Int (\partial_{\cap} D)\subset \Int (H_F' )$). La surface
$\overline{H_F \setminus D}$, une fois les coins de $D$ lissés,
est un hexagone $H_F'$ à bord legendrien qui est $\xi$-convexe
et isotope à $H_F$. Par construction,  $tb (H_F ')=tb (H_F )+1$.
L'isotopie entre $H_F$ et $H_F'$
se prolonge en une isotopie $(\phi_t )_{t\in [0,1]}$ de $V$
à support dans un voisinage de
$D$ et donc stationnaire sur $\Lambda$.
Le $1$-squelette
de $\Delta$ est legendrien pour $\xi' =\phi_1^* \xi$ et l'arc
$r(a)$ possède  un nombre de Thurston-Bennequin relatif aux faces
adjacentes inférieur ou égal à $-\frac{1}{2}$ (et même à $-1$
relativement à $F$) et strictement
supérieur à celui donné par $\xi$. Une isotopie de
$\xi'$ relative au $1$-squelette de $\Delta$ et à $\Lambda$
permet de rendre les faces $\xi'$-convexes et relativement
$\xi'$-convexes. Elle
fait  de $\Delta$ une
triangulation de contact maniable pour la nouvelle structure $\xi''$.
Ce faisant, on a strictement diminué le
nombre de Thurston-Bennequin et donc
$\Delta$ n'était pas $\Lambda$-minimale pour $\xi$.
\end{proof}

Si $F$ est un $2$-simplexe de $\Delta$, un \emph{quadrilatère fibré} dans $F$ 
est un quadrilatère $[0,1] \times [0,1] \subset F$ dont l'intersection avec
$\partial F$ est l'union des deux arêtes \emph{verticales} $\{0\} \times [0,1]$
et $\{1\} \times [0,1]$, celles-ci se trouvant à l'intérieur de deux arêtes 
distinctes de $F$.
Un arc simple et propre $A \subset F$ 
est \emph{porté} par un quadrilatère fibré
$Q$ s'il est inclus dans $Q \setminus \{ y=0,1\}$ et s'il est transversal au
champ de  vecteurs
$\partial_y$.

Soit $(V,\xi)$ une variété de contact tendue de dimension~$3$ et $\Delta$ une 
triangulation de contact de $(V,\xi)$. Pour toute face $F$ de $\Delta$,
on note $\Gamma_{\!F}$ le découpage de $\xi F$ associé à un épaississement 
homogène quelconque de~$F$. On appelle \emph{pièce} de~$F$ l'adhérence de toute 
composante connexe de $F \setminus \Gamma_{\!F}$ et on dit qu'une pièce est 
\emph{ordinaire} ou \emph{extraordinaire} selon que c'est ou non un quadrilatère
fibré. On insiste ici sur le fait qu'on ne fixe
pas l'épaississement homogène une fois pour toutes et qu'on s'autorise
donc des isotopies de  $\Gamma_F$ parmi les multi-courbes transversales
à $\xi F$.

\begin{corollaire} \label{c:quadrilatere}
Soit $\XX$ un syst\`eme complet de structures de contact tendues sur
$V$ et $\Delta$ une triangulation de $V$ qui est  maniable 
et $\Lambda$-minimale pour tout $\xi \in \XX$, avec  un voisinage de
s\'ecurit\'e $\Lambda$ fixe le long duquel toutes les structures 
de $\XX$ co\"\i ncident.
Il existe $C_0>0$ tel que, pour tout $\xi \in \XX$, toute face $F$ de $\Delta$ 
contienne trois quadrilatères fibrés $Q_1, Q_2, Q_3$ deux à deux disjoints, 
s'appuyant sur des paires d'arêtes de $F$ distinctes et inclus dans $H_F 
\setminus \Lambda$ \up(en particulier, pour $i \in \{ 1, 2, 3 \}$, les arêtes 
verticales de $Q_i$ sont incluses dans $\bigcup_{a \subset \Delta^1} r(a)$\up), 
avec les propriétés suivantes\up{~:}
\begin{itemize}
\item
chaque $Q_i$ est une union de pièces ordinaires\up{\,;}
\item
au plus $C_0$ pièces du découpage de $\xi F$ ne sont pas incluses dans $Q_1 \cup
Q_2 \cup Q_3$.
\end{itemize}
\end{corollaire}

\begin{figure} [ht]
\begin{center}
    \resizebox{8cm}{!}{\input{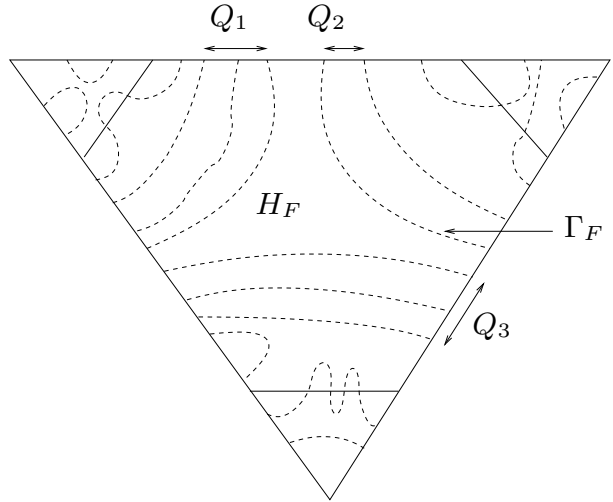}}
    \caption{Le découpage de $F$ et $H_F$.}
    \label{fig13}
\end{center}
\end{figure}

\begin{proof}
Il découle du lemme~\ref{l:extremale} que les pièces ordinaires qui 
n'intersectent pas $\Lambda$ forment trois
quadrilatères fibrés \guil{non parallèles}. Vu le
lemme~\ref{lemme : cercle}, restent comme pièces au plus six
demi-disques extrémaux, une pièce \guil{centrale}, et des
pièces qui rencontrent $\Lambda$ et dont le nombre est donc
borné indépendemment de $\xi \in \XX$ (voir la
figure~\ref{fig13}).
\end{proof}

\subsection{Prismes fibrés}

Soit $Y$ un triangle ou un quadrilatère. Un prisme fibré est un polyèdre $P = Y
\times [0,1]$ dont on ne retient de la structure produit que la projection sur
$Y$. Les faces de $P$ sont dites \emph{verticales} ou \emph{horizontales} selon
qu'elles se trouvent dans $\partial Y \times [0,1]$ ou dans $Y \times \{0,1\}$.

Soit $\Delta$ une triangulation de $V$. Un prisme fibré dans $(V,\Delta)$ est un
plongement d'un prisme fibré $P$ dans $V$ avec les propriétés suivantes~:
\begin{itemize}
\item
$P$ est contenu dans un $3$-simplexe $G$ de $\Delta$ et son intersection avec
$\partial G$ est l'union de ses faces verticales\,;
\item
chaque face verticale de $P$ est un quadrilatère fibré d'une face de $G$.
\end{itemize}
Un $3$-simplexe donné contient au plus cinq prismes fibrés deux à deux disjoints
et non isotopes parmi les prismes fibrés. Il y a, à isotopie près parmi les 
prismes fibrés, trois telles familles de cinq prismes fibrés~: chacune d'elle 
possède quatre prismes à base triangulaire, situés près des sommets, et un
prisme de base quadrilatérale, situé en diagonale.

Une \emph{configuration de prismes fibrés} dans $(V,\Delta)$ est la donnée d'une
collection de prismes fibrés $P = (P_i)_{1 \le i \le k}$ de $(V,\Delta)$ qui 
intersecte chaque $3$-simplexe de $\Delta$ en une sous-famille de l'une des 
trois familles maximales de cinq prismes précédentes et telle que l'intersection
entre deux prismes $P_i \ne P_j$ de $P$ donne soit un quadrilatère fibré (avec 
concordance des fibrations données de part et d'autre), soit un arc d'intérieur 
non vide inclus dans une arête de $\Delta$, soit l'ensemble vide. Dans le 
premier cas, les prismes $P_i$ et $P_j$ sont contenus dans des $3$-simplexes qui
ont une face en commun et dans le deuxième cas, dans des $3$-simplexes qui ont
une arête commune.

\begin{lemme} \label{l:fini}
Toute variété triangulée $(V,\Delta )$ possède un nombre
fini de configurations de prismes fibrés, à isotopie près
parmi les configurations de prismes fibrés.
\end{lemme}

Une configuration de prismes fibrés $P$  est dite \emph{admissible}
pour $\xi \in \XX$ si les faces verticales des prismes de la famille $P$
sont des unions de pièces ordinaires qui ne rencontrent pas $\Lambda$.
 Pour $\xi \in \XX$, on note $\P_{\Delta ,\xi}$ l'ensemble des
configurations de prismes fibrés admissibles pour $\xi$. On munit
les classes d'isotopies de prismes fibrés dans $\P_{\Delta ,\xi}$
d'une relation d'ordre partiel~: si $P,Q \in \P_{\Delta ,\xi}$
représentent des classes d'isotopies $[P]$ et $[Q]$, on dit que
$[P]\preceq [Q]$ s'il existe une isotopie de $P$ dans $\P_{\Delta
,\xi}$ en $P'\subset Q$.

\begin{lemme}
Pour $\xi \in \XX$ fixée, les classes d'isotopie de
$\P_{\Delta,\xi}$ sont en nombre fini.
\end{lemme}

\begin{proof}
Il suffit d'observer qu'un élément de $\P_{\Delta ,\xi}$
est déterminé, à isotopie près dans $\P_{\Delta ,\xi}$, par
la classe d'isotopie de ses arêtes horizontales parmi les arcs
transversaux aux feuilletages caractéristiques des faces. Or toute
arête horizontale incluse dans une face $F$ est isotope à une
composante de $\Gamma_F (\xi )$, ce qui donne un nombre fini de
classes d'isotopie possibles pour chaque arête horizontale, mais
également pour leur collection.
\end{proof}

Les sous-sections suivantes sont consacrées à la démonstration
du lemme fondamental ci-dessous~:

\begin{lemme}\label{l:prismes}
Soit $\XX$ un ensemble complet de structures de contact tendues sur
$V$ et $\Delta$ une triangulation de $V$ qui est  maniable 
et $\Lambda$-minimale pour tout $\xi \in \XX$, avec  un voisinage de
s\'ecurit\'e $\Lambda$ fixe le long duquel toutes les structures 
de $\XX$ co\"\i ncident.
Il existe $C_1 >0$ tel que pour tout $\xi \in \XX$, et pour tout
élément  $P=(P_i )_{1\leq i\leq n}$ de $\P_{\Delta ,\xi }$ dont
la classe d'isotopie est  maximale, au plus $C_1$ pièces des faces
de $\Delta$ ne sont pas incluses dans $\Int (\bigcup_{1\leq i\leq n}
P_i )$.
\end{lemme}

\begin{remarque}\label{remarque : minoration}
Quitte à considérer des sous-familles de prismes $P$ des
éléments maximaux de $\P_{\Delta ,\xi}$ et à augmenter la
valeur de $C_1$, on peut également assurer que chaque face
verticale des prismes de $P$ contient au moins $20$ pièces
(ordinaires).
\end{remarque}

\subsubsection{Une première normalisation des faces}

Pour toute arête $a$ de $\Delta$, on note $s(a) =
\bigcup_{F\supset a} (a\cap H_F )$.

\begin{lemme}\label{lemme : normalisation}
Quitte, pour tout $\xi \in \XX$, à effectuer une isotopie de $\xi$,
on peut supposer que, si $\xi \in \XX$~:
\begin{enumerate}
\item $\Delta$ est une triangulation de contact maniable et
$\Lambda$-minimale pour $\xi$\,;
\item pour toute arête $a\subset \Delta^1$,
il existe un germe de feuilletage $\F_a$ au voisinage de $s(a)$ par
des arcs parallèles à $s(a)$, $s(a)$ étant une feuille, qui
est tangent au germe des faces adjacentes à $s(a)$ et qui est
legendrien pour tout $\xi \in \XX$\,;
\item pour toute face $F$ et pour toute pièce ordinaire $R$
incluse dans $F\setminus \Lambda$,
$$\xi \vert_R =\{ \sin \theta dx+\cos \theta dy =0 \}$$
si $R\simeq \{ y=0\} \subset \{ (x,y,\theta )\in [0,1]\times
[-1,1]\times [-\pi /2, \pi /2] \}$. En particulier, le feuilletage
caractéristique de $R$ possède une courbe de singularités $\{
y=0, \theta =0\}$ portée par $R$.
\end{enumerate}
\end{lemme}

\begin{proof}
La deuxième condition n'est \emph{a priori} pas réalisée dans
le voisinage de sécurité $\Lambda$. Pour démontrer ce lemme,
on reprend la démonstration de la proposition~\ref{p:maniable},
dans laquelle on exige une propriété supplémentaire pour la
structure $\xi_0$~: on  déforme $\xi_0$ par isotopie  à l'aide
du lemme de Darboux pour obtenir le point 2. C'est possible sans
stabilisation supplémentaire puisque tous les enroulements
relatifs aux faces sont négatifs. On déroule la preuve de la
proposition~\ref{p:maniable} pour cette structure $\xi_0$, dont on
déduit la construction d'un système complet $\XX$ convenable.
Pour une structure $\xi \in \XX$ quelconque, on a alors 2 au
voisinage de  $s(a) \cap \Lambda$ où $\xi =\xi_0$. Le lemme de
Darboux, appliqué à  $\xi \in \XX$ au voisinage d'une arête $a$
relativement à ce qui a déjà été fait au voisinage de
$\Lambda$, fournit une isotopie de $\xi$ stationnaire sur  $\Lambda$
qui conduit à une structure, notée à nouveau $\xi$, qui
possède également un feuilletage legendrien par des arcs
parallèles à $s(a)$ et tangent aux faces près de $s(a)$.

Pour obtenir un germe de feuilletage  legendrien indépendant de
$\xi \in \XX$, on se sert du fait que le nombre de Thurston-Bennequin
relatif à toute face $F$ contenant $a$ le long de $r(a)$ est
inférieur à $-\frac{1}{2}$, ce pour tout $\xi \in \XX$.
Précisément, sur un voisinage $K$ de $s(a)$, $\xi$ est solution
d'une équation $\cos f_0(x,y,\theta ) dx-\sin f_0(x,y,\theta
)dy=0$ dans des coordonnées $(x,y,\theta )\in D^2 \times [0,1]$,
$s(a)=\{ x=y=0\}$, données par un plongement $\phi_0 : D^2 \times
[0,1] \to V$ et pour lesquelles les faces sont
$\partial_\theta$-invariantes. On déforme le feuilletage
legendrien dirigé par $\partial_\theta$ sur un petit voisinage de
$s(a)$ dans $K$ relativement à $\partial K$ tout en préservant
les faces, pour le faire coïncider avec le feuilletage  donné
par $\xi_0$ tout près de $s(a)$. On obtient ainsi une isotopie
$(\phi_t )_{t\in [0,1]}$ de $\phi_0$ pour laquelle $\phi_{1*}
\partial_\theta$ est tangent à $\xi_0$ près de $s(a)$ et $\phi_t
\vert_{\partial K} =\phi_0 \vert_{\partial K}$ pour tout $t\in
[0,1]$. À chaque instant $t\in [0,1]$, la structure $\phi_t ^*
\xi$ est donnée près de $D^2 \times \{ 0,1\}$ par une équation
$\cos f_t (x,y,i) dx -\sin f_t (x,y,i)dy =0$, $i=0,1$, où $f_t$
dépend continûment de $t\in [0,1]$. Le fait que l'enroulement de
$\xi$ le long de $s(a)$ par rapport à une face quelconque
d'arête $a$ soit  inférieur à $-\frac{1}{2}$ implique que pour
tout $t\in [0,1]$, $f_t (x,y,1)>f_t (x,y,0)$. On peut donc étendre
$f_t$ en une famille continue  de fonctions sur $D^2 \times [0,1]$
constamment égale à $f$ près du bord avec la propriété
$\partial_\theta f_t >0$. L'équation $\cos f_t (x,y,\theta )dx
-\sin f_t (x,y,\theta )dy=0$ définit une structure de contact sur
$D^2 \times [0,1]$, dont l'image par $\phi_t$ se recolle à $\xi$
hors de $K$ pour former un chemin de structures de contact entre
$\xi$ et une structure $\xi'$ tangente au même feuilletage
legendrien que $\xi_0$ près des arêtes. Le théorème de Gray
convertit ce chemin en une isotopie de $V$ stationnaire sur
$(V\setminus K) \cup \Delta^1$.

On obtient 3 relativement à la déformation déjà effectuée
près des arêtes à l'aide du lemme de réalisation de feuilletages.
\end{proof}

Dans la suite, le syst\`eme  $\XX$ qu'on consid\`ere
poss\`ede, en plus des pr\'ec\'edentes, les propri\'et\'es explicit\'ees dans
les conclusions du lemme~\ref{lemme : normalisation}.

\subsubsection{Holonomie}

Soit $G$ un $3$-simplexe dans une variété de contact $(V,\xi )$.
Un arc legendrien de classe $\classe1$ par morceaux $\gamma_1 :[0,1]
\to \partial G$ qui évite les sommets de $G$ est dit \emph{étal}
si le champ de plans $\xi$ est un plan d'appui à $\partial G$ le
long de $\gamma_1$. Précisément, $\gamma_1$ intersecte
l'intérieur de chaque face le long d'une ligne singulière de son
feuilletage caractéristique et, le long d'une arête, $\xi$ est
à l'\guil{extérieur} de $G$.

\begin{figure} [ht]
\begin{center}
    \resizebox{8cm}{!}{\input{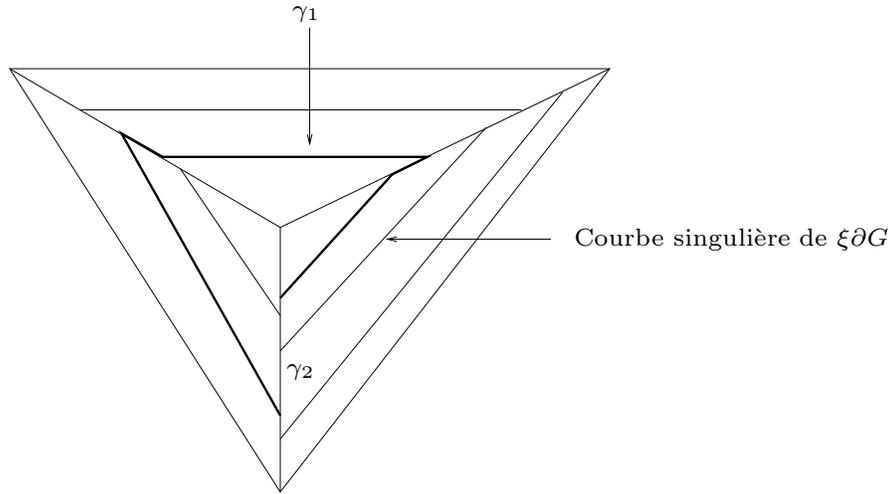}}
    \caption{Courbe d'holonomie $-1$.}
    \label{fig14}
\end{center}
\end{figure}

Une \emph{courbe d'holonomie} est une courbe legendrienne $\gamma$
dans $\partial G$ constituée de la concaténation de deux arcs
$\gamma_1$ et $\gamma_2$, où $\gamma_1$ est étal et  $\gamma_2$
inclus dans une arête de $G$. On appelle champ de plans {\it
médian} de $G$ le long de $\gamma_2$, tout champ de plans tangent
à $\gamma_2$ qui rencontre l'intérieur du secteur de
$T_{\gamma_2} V$ délimité par $G$, c'est-à-dire qui n'est un
plan d'appui en aucun point de $\gamma_2$. L'\emph{holonomie} $Hol
(\gamma )$ d'une courbe d'holonomie $\gamma$ est un entier dont la
valeur absolue est égale à celle du nombre de Thurston-Bennequin
de $\xi$ le long de $\gamma_2$, calculé relativement à un  champ
de plans médian quelconque de $G$ donné le long de $\gamma_2$.
La valeur absolue de l'holonomie est donc la moitié du nombre de
points de $\gamma_2$, comptés algébriquement, où $\xi$ est
égal au plan médian. Il est positif si $\gamma_2$ est orienté
comme le bord de la face qui contient $\gamma_1$ près de $\gamma_1
(0)$ et négatif sinon. La figure~\ref{fig14} montre un exemple de
courbe d'holonomie $-1$. Sur cette figure, lorsqu'on change le sens
de parcours de $\gamma$, on change également la face servant à
étalonner l'orientation de $\gamma_2$ et le signe de $Hol (\gamma
)$ ne varie pas. C'est toujours le cas lorsque $\gamma_1$ aborde
$\gamma_2$ par deux faces différentes.

\begin{lemme}\label{lemme : holonomie}
Soit $\xi \in \XX$ et $G$ un $3$-simplexe de $\Delta$. Toute courbe
d'holonomie $\gamma$ pour $\xi$ incluse dans $\partial G \setminus
\Lambda$ a une holonomie égale à $-1$.
\end{lemme}

\begin{remarque}\label{remarque : holonomie}
Si $\gamma =\gamma_1 \cup \gamma_2$ est une courbe d'holonomie pour
laquelle $\gamma_1$ aborde $\gamma_2$ des deux bouts par la même
face, et si $\overline{\gamma}$ est $\gamma$ orientée en sens
inverse, alors $Hol(\overline{\gamma} )= -Hol (\gamma )$. En
particulier ce cas de figure est exclus par le lemme~\ref{lemme :
holonomie}.
\end{remarque}

\begin{proof}
Comme $\xi$ est tendue, l'inégalité de Bennequin appliquée à
$\gamma$ dit que $Hol (\gamma )\neq 0$. On suppose que $Hol (\gamma
)\neq -1$. La courbe $\gamma$ se décompose par définition en la
réunion de deux arcs $\xi$-legendriens $\gamma_1$ et $\gamma_2$,
où $\gamma_1$ est étal et $\gamma_2$ est inclus dans un certain
$r(a)$, $a\in \Delta^1$. Comme $\gamma_1$ est étal, on peut le
pousser sur un arc legendrien lisse $\gamma_1'$ situé à
l'intérieur de $G$ et s'appuyant sur $\partial \gamma_2$, {\it
via} une famille d'arcs legendriens. En particulier, la courbe
fermée $\gamma_1' \cup \gamma_2$ borde un demi-disque $D$ avec
$\Int (D)\subset \Int (G)$. Le nombre de Thurston-Bennequin de
$\gamma_2$ relatif à  $D$ vaut $- \vert Hol (\gamma )\vert$ si
$\gamma_1$ aborde les deux bouts de $\gamma_2$ par le même
côté, et sinon $- \vert Hol (\gamma )\vert +\frac{1}{2} $
lorsque  $Hol (\gamma )<0$ et $-\vert Hol (\gamma ) \vert
-\frac{1}{2} $ lorsque $Hol (\gamma )>0$. En effet, cet enroulement
est le même que celui de $\xi$ le long de $\gamma_2$, calculé
par rapport à un champ de plans qui est médian de $G$ le long de
$\Int (\gamma_2 )$ et égal à $\xi$ (i.e. tangent à une face de
$G$) aux points de  $\partial \gamma_2$. Selon les cas, un tel champ
de plans a un enroulement $0$, $-\frac{1}{2}$, où $\frac{1}{2}$
par rapport à un champ de plans médian. De la même façon,
comme $\gamma_1$ est étal, l'invariant de Thurston-Bennequin de
$\gamma_1'$ relatif à $D$ vaut respectivement $0$, $-\frac{1}{2}$
et $\frac{1}{2}$ dans les trois cas énumérés ci-dessus.
 Dans toutes les situations, puisque $Hol (\gamma )\neq -1$,  on a
$tb(\gamma_1' ,D)\geq tb(\gamma_2 ,D)+1$. On stabilise $\gamma_1'$
pour obtenir un arc $\gamma_1''$, tel que $tb(\gamma_1'',D'
)=tb(\gamma_2 ,D')+1$ (si $D' \subset G$ désigne l'image de $D$
par l'isotopie de stabilisation). Le disque $D'$ va jouer le rôle
de rocade.

L'isotopie qui consiste à pousser $\gamma_2$ sur $\gamma_1''$ le
long de $D'$ peut être réalisée  relativement à $\Lambda$.
Appliquée à $\xi$ elle fait de $\Delta$ une triangulation de
contact pour son image $\xi'$ (après déformations habituelles
près des faces), avec un  nombre de Thurston-Bennequin strictement
moindre que celui de $\xi$. On constate de plus facilement que comme
$Hol (\gamma )\neq -1$, pour toute face $F$ de $\Delta$ et toute
arête $a$ de $F$, le nombre de Thurston-Bennequin le long de
$r(a)$ relatif à $F$ reste inférieur à $-\frac{1}{2}$. Cette
triangulation est donc maniable pour $\xi'$ et on obtient une
contradiction avec la $\Lambda$-minimalité de $\Delta$.
\end{proof}

\subsubsection{Démonstration du lemme~\ref{l:prismes}}

Soit $G$ un $3$-simplexe de $\Delta$. On note $s_1$, $s_2$, $s_3$ et
$s_4$ les sommets de $G$. La notation $[s_i s_j]$ désigne
l'arête de $G$ qui joint les sommets $s_i$ et $s_j$. On note $(s_i
s_j s_k )$ la face de $G$ qui contient les sommets $s_i$, $s_j$ et
$s_k$.

Un \emph{paquet} est une union de pièces ordinaires contenues dans
une même face qui forme un quadrilatère fibré (connexe).
L'\emph{épaisseur} d'un paquet est le nombre de pièces qui le
constituent.

Soit $\xi \in \XX$ et $P=( P_i )_{1\leq i\leq n} \in \P_{\Delta
,\xi}$ une configuration de prismes admissible maximale. Les
fonctions $f_i :\R \rightarrow \R$ qui apparaissent dans la suite
sont affines, de dérivées strictement positives et
indépendantes de $(\xi , P)\in \XX \times \P_{\Delta ,\xi}$. On
raisonne par l'absurde en supposant que $N\gg 1$ pièces  ne sont
pas incluses dans $\Int (\bigcup_{1\leq i\leq n} P_i )$. Parmi
celles-ci  au moins $N_1 =f_1 (N)$ pièces se trouvent dans une
même face $F$ de $\Delta$. Le corollaire~\ref{c:quadrilatere}
nous dit alors qu'au moins $N_2 =f_2 (N_1 )$ pièces ordinaires
sont situées dans un même hexagone $H_F \setminus \Lambda$\,;
leur union formant un paquet $Q_1$ (d'épaisseur $N_2$), lui-même
contenu dans $H_F \setminus \Int_G((\bigcup_{1\leq i\leq n} P_i)\cap
G)$ pour un des deux simplexes $G$ de $\Delta$ contenant $F$.
C'est-à-dire que $Q_1$ ne rencontre pas la réunion $P_G$ des
prismes de $P$ inclus dans $G$. Soit $s_1$, $s_2$ et  $s_3$ les
sommets de $F$ et $s_4$ le quatrième sommet de $G$. Pour fixer les
idées, disons que le paquet $Q_1$ joint $[s_1 s_2 ]$ à $[s_1 s_3
]$. Au moins  $N_2 -1$ composantes de $\Gamma_{(s_1 s_3 s_4 )} (\xi
)$ partent de $[s_1 s_3 ] \cap Q_1$, dont au moins $N_3 =f_3 (N_2 )$
restent dans $H_{(s_1 s_3 s_4 )} \setminus \Lambda$ et  ne sont pas
parallèles à une arête (corollaire~\ref{c:quadrilatere}).

\begin{figure} [ht]
\begin{center}
    \resizebox{16cm}{!}{\input{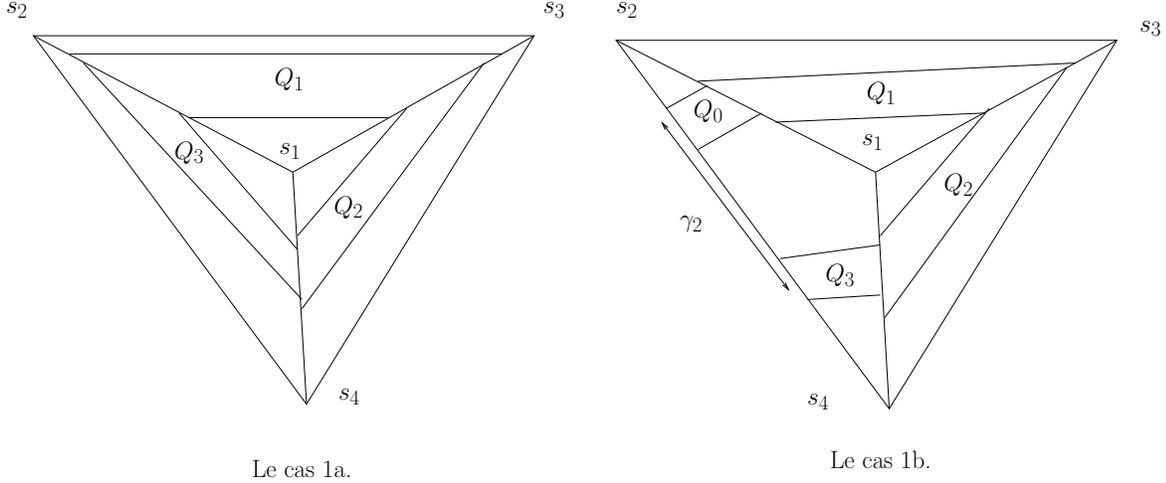}}
    \caption{Les cas 1a et 1b.}
    \label{fig10}
\end{center}
\end{figure}

{\bf Cas 1.} Parmi celles-ci, au moins la moitié $N_4 =f_4 (N_3 )$
vont vers $[s_1 s_4 ]$ et délimitent un paquet $Q_2$ d'épaisseur
$N_4 -1$. Comme toutes ont une extrémité dans $[s_1 s_3 ] \cap
Q_1$, $Q_2$ ne rencontre pas $P_G$. On distingue à présent deux
cas (figure~\ref{fig10})~:

{\bf Cas 1a.}  Au moins $11$ (pour être tranquille) composantes de
$\Gamma_{H_{(s_1 s_2 s_4 )}} (\xi )$ issues de $Q_2 \cap [s_1 s_4 ]$
reviennent vers $r([s_1 s_2])$. On en déduit l'existence d'un
paquet $Q_3$ d'épaisseur au moins $10$ entre $[s_1 s_4 ] \cap Q_2$
et $[s_1 s_2]$ qui ne rencontre pas $P_G \cup \Lambda$.

\begin{lemme}
Il existe un prisme fibré $P'$ dont le bord vertical est inclus
dans  $Q_1 \cup Q_2 \cup Q_3$ et est une union de pièces
ordinaires. De plus, ce prisme ne rencontre pas $P_G \cup \Lambda$.
\end{lemme}

\begin{proof}[Démonstration.]
Il s'agit de construire les faces verticales de $P'$. On ordonne les
pièces contenues dans $Q_3$ de un à dix. Chaque pièce porte
une  courbe legendrienne singulière, numérotée comme la
pièce dont elle fait partie. On part de la cinquième courbe
singulière de  $\xi Q_3$, notée $c_3$. L'extrémité de
celle-ci dans $[s_1 s_4 ]$ est encadrée par les extrémités de
$2$ courbes singulières de $\xi Q_2$. On note $c_2$, celle dont
l'extrémité est située juste après  $\partial c_3 \cap [s_1
s_4]$, pour l'orientation de $[s_1 s_4]$ donnée comme bord de la
face $(s_1 s_2 s_4 )$. En particulier, la portion d'arête de $[s_1
s_4 ]$ située entre $c_3$ et $c_2$ n'est \emph{pas} un arc étal.
L'extrémité de $c_2$ dans $[s_1 s_3 ]$ est encadrée  par les
extrémités de deux courbes singulières de $\xi Q_1$. À
nouveau, on note $c_1$ celle dont l'extrémité est située
après $\partial c_2 \cap [s_1 s_3 ]$ pour l'orientation de $[s_1
s_3]$ induite par celle de la face $(s_1 s_3 s_4 )$. Comme
l'holonomie de toute courbe d'holonomie vaut $-1$, l'extrémité
$c_1 \cap [s_1 s_2]$ est  voisine de celle de  $c_3$, comme sur la
figure~\ref{fig12} (ce qui ne serait pas le cas si on avait choisit
$c_1$, $c_2$ et $c_3$ pour construire un arc étal). Soit $U_1$,
$U_2$ et $U_3$ des petits voisinages tubulaires de $c_1$, $c_2$ et
$c_3$ dans respectivement $Q_1$, $Q_2$ et $Q_3$. On note de plus
$V_1$, $V_2$ et $V_3$ des voisinages, dans $\partial G$, des arcs
délimités dans les arêtes $[s_1 s_3]$, $[s_1 s_4]$ et $[s_1
s_2]$ par  $\bigcup_{1\leq i\leq 3} (\partial c_i )$. Les faces
verticales recherchées sont obtenues par un lissage de
$$U_1 \cup U_2 \cup U_3 \cup V_1 \cup V_2 \cup V_3.$$

\begin{figure} [ht]
\begin{center}
    \resizebox{13cm}{!}{\input{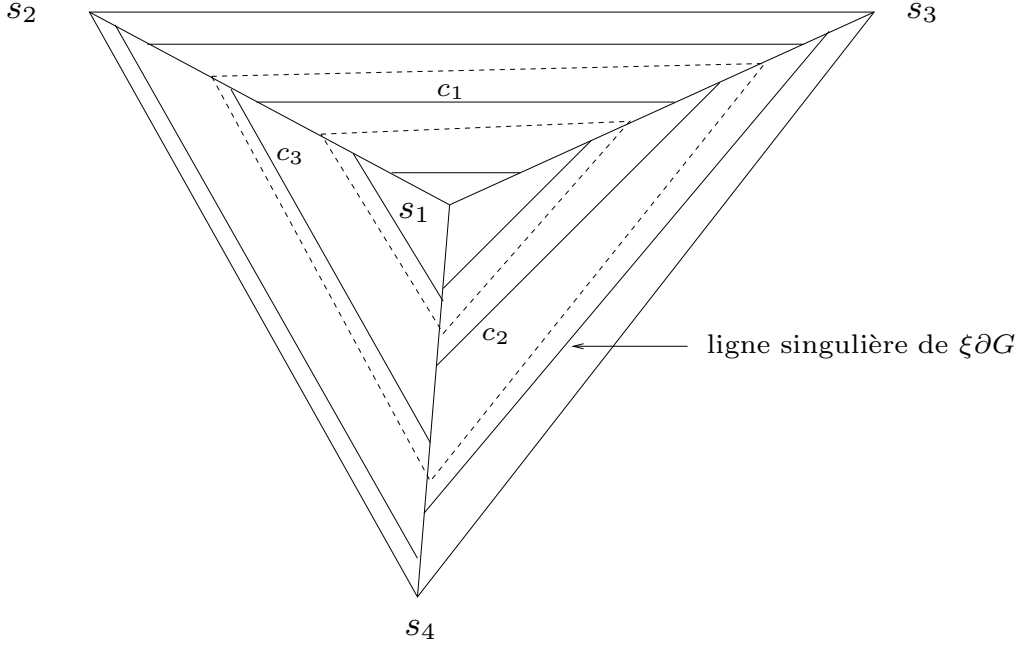}}
    \caption{Les faces verticales de $P'$.}
    \label{fig12}
\end{center}
\end{figure}

Elles sont par construction incluses dans $Q_1 \cup Q_2 \cup Q_3$ et
elles ne rencontrent donc pas $P_G \cup \Lambda$. Ce sont les faces
verticales d'un prisme fibré $P'$ inclus dans $G$ qui ne rencontre
pas $P_G \cup \Lambda$. Sur la figure~\ref{fig12}, les faces
verticales de $P'$ sont délimitées par les traits pointillés\,;
les traits pleins représentent les lignes singulières du
feuilletage des faces de $G$.
\end{proof}

On déduit facilement de ce lemme que la famille $P$ n'était pas
maximale~: si $P'$ n'est parallèle à aucun prisme de $P_G$, on
ajoute $P'$ à $P$\,; si $P'$ est parallèle à un prisme $P_1$ de
$P_G$, on remplace $P_1$ par un prisme qui contient $P'$ et $P_1$ et
dont le bord horizontal est contenu dans $\partial P' \cup \partial
P_1$.

{\bf Cas 1b.} On peut construire un paquet $Q_3 \subset (s_1 s_2
s_4) \setminus \Lambda$ entre $[s_1 s_4]\cap Q_2$ et $[s_2 s_4]$
d'épaisseur $N_5 =N_4-10-C_0 =f_5 (N_3 )$.

Dans ce cas, on considère le sous-paquet de $Q_1$ d'épaisseur
$N_4-1$ contenant toutes les pièces issues de $Q_2 \cap [s_1
s_3]$. On rebaptise $Q_1$ ce nouveau paquet. On peut alors
construire un paquet $Q_0 \subset (s_1 s_2 s_4 )\setminus \Lambda$
d'épaisseur au moins $N_6 =N_4 -1-10-C_0 =f_6 (N_4 )$ entre $[s_1
s_2] \cap Q_1$ et $[s_2 s_4]$, sinon on se retrouve dans le cas
$1a)$ avec un prisme fibré dont le bord vertical tourne autour de
$s_1$.

On construit une courbe d'holonomie $\gamma =\gamma_1 \cup \gamma_2$
pour laquelle $\gamma_1$ possède un voisinage de ses deux bouts
dans la même face, ce qui donne une contradiction avec le
lemme~\ref{lemme : holonomie} au vu de la remarque~\ref{remarque :
holonomie}. L'arc étal $\gamma_1$ est la concaténation de $4$
arcs singuliers $c_0 \subset Q_0$, $c_1 \subset Q_1$, $c_2 \subset
Q_2$, $c_3 \subset Q_3$, du feuilletage caractéristique des faces
et de $3$ portions d'arêtes qui relient leurs extrémités.
L'arc $\gamma_2$ qui mesure l'holonomie sera pris dans $[s_2 s_4]$.
Pour cela, pour $N$ assez grand, on prend  pour $c_3$ un arc
singulier du feuilletage de $Q_3$ qui  découpe $Q_3$ en deux
paquets d'épaisseur au moins $20+C_0$ (pour être tranquille).
Les choix des arcs singuliers $c_0$, $c_1$ et $c_2$ de même que
celui des trois arêtes sont alors imposés par la condition
d'être étal. Par la position de $c_3$, on est assuré que $c_2
\subset Q_2$ et $c_1 \subset Q_1$. De plus, $c_1$ découpe $Q_1$ en
deux paquets d'épaisseur au moins $19+C_0$. Comme au plus $14
+C_0$ composantes de $\Gamma_{(s_1 s_2 s_4 )} (\xi )$ (estimation
lâche) ayant leur extrémités dans $Q_1$ ne sont pas dans
$Q_0$, on a bien $c_0 \subset Q_0$. C'est l'arc $\gamma_1$
recherché.

{\bf Cas 2.} On peut construire un paquet $Q_2 \subset (s_1 s_3 s_4
)\setminus \Lambda$ d'épaisseur $N'_4 =f'_4 (N_3 )$ entre $[s_1
s_3 ] \cap Q_1$ et $[s_3 s_4]$.

\begin{figure} [htpb]
\begin{center}
    \resizebox{10cm}{!}{\input{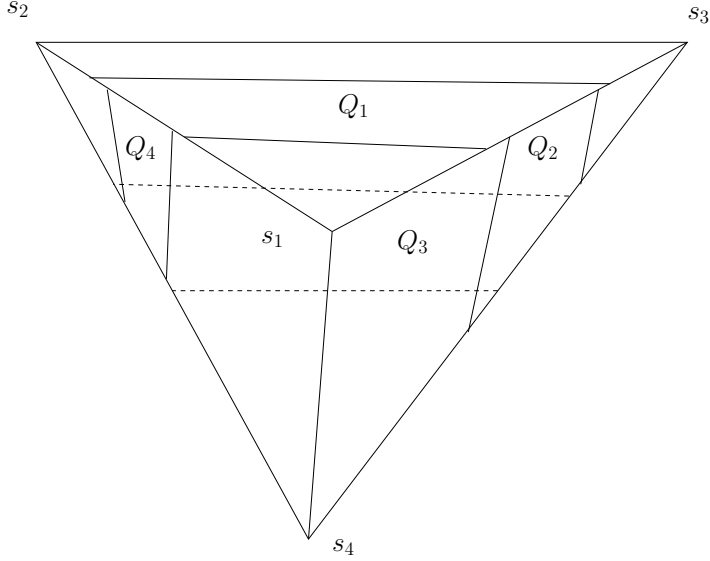}}
    \caption{Le cas 2.}
    \label{fig11}
\end{center}
\end{figure}

On obtient ensuite un paquet $Q_3 \subset (s_2 s_3 s_4 )\setminus
\Lambda$ d'épaisseur au moins $N'_5 =f'_5 (N'_4 )$ entre $[s_3 s_4
]\cap Q_2$ et $[s_2 s_3]$, ce qui, pour $N$ assez grand, ramène au
cas $1)$, ou entre $[s_1 s_4 ]\cap Q_2$ et $[s_2 s_4]$. À nouveau,
le seul cas non traité est celui où on peut construire un paquet
$Q_4 \subset (s_1 s_2 s_4 )\setminus \Lambda$ d'épaisseur au moins
$11$ entre $[s_2 s_4 ]\cap Q_3$ et $[s_1 s_2 ]$ (voir
figure~\ref{fig11}).
Dans ce dernier cas, on construit, comme dans le cas $1a$, un
prisme fibré dont le bord vertical est dans $Q_1 \cup Q_2 \cup Q_3
\cup Q_4$. Comme précédemment, la famille $(P_i)_i$ n'était
donc pas maximale.

\subsection{Construction des domaines fibrés}

On démontre à présent le théorème~\ref{t:domaines}.

\subsubsection{Normalisation des faces verticales}

Les lemmes ci-dessous permettent, par des isotopies et des partitions
successives de $\XX$, d'affiner la normalisation des structures de $\XX$.
On rappelle que  $\XX$ est un ensemble complet de structures de contact tendues sur
$V$ et que $\Delta$ est  une triangulation de $V$,  maniable 
et $\Lambda$-minimale pour tout $\xi \in \XX$, toutes
les structures de $\XX$ \'etant \'egales le long de leur voisinage de
s\'ecurit\'e commun $\Lambda$. En outre, les structures de $\XX$ ont d\'ej\`a subi
une premi\`ere normalisation et satisfont toutes aux 
conclusions du lemme~\ref{lemme : normalisation}.

\begin{lemme} \label{l:partition1}
Il existe $C_2>0$ et un nombre fini de configurations de prismes fibrés $P^1,
\dots, P^k$ tels que, pour  toute structure $\xi \in \XX$, on puisse isotoper 
$\xi$, par une isotopie de $V$ stationnaire sur $\Lambda$ et préservant 
$\Delta$, en $\xi'$ admettant un des $P^j$ pour configuration admissible et pour
laquelle au plus $C_2$ composantes de $\Gamma_{\Delta^2} (\xi ')$ ne soient pas 
incluses dans $\Int (\bigcup_{P \in P^j} P)$.
\end{lemme}

Une structure de contact $\xi \in \XX$ qui vérifie ces propriétés pour une 
configuration $P^j$ sera dite \emph{portée} par $P^j$.

\begin{proof}
\`A toute structure $\xi \in \XX$, on peut associer un élément
maximal $P(\xi )$ dans $\P_{\Delta ,\xi}$. Le lemme~\ref{l:fini} 
affirme qu'il n'y a qu'un nombre fini de classes d'isotopie
de configurations de prismes fibrés. Une isotopie de $V$ fixant
$\Delta$ et relative à $\Lambda$ permet d'envoyer chaque $P(\xi )$
sur un élément d'une sous-famille finie $P^1$,...,$P^k$ de la
famille $P(\xi )_{\xi \in \XX}$. Cette isotopie transporte $\xi$ sur
une structure $\xi'$. Le lemme~\ref{l:prismes} affirme
l'existence de la borne universelle $C_2 =C_1 >O$.
\end{proof}

Dans la suite, on remplace $\xi$ par $\xi'$ comme représentant de
la classe d'isotopie de $\xi$ dans $\XX$. On démontre le
théorème~\ref{t:domaines} pour l'ensemble $\XX_1$ des structures
$\xi \in \XX$ portées par $P^j =P=(P_i)_{1\leq i\leq n}$.

\begin{lemme}\label{l:partition2}
Quitte à déformer chaque structure $\xi \in \XX_1$ par une
isotopie relative à $\Lambda$ qui préserve $\Delta$, il existe
une partition de $\XX_1$ en un nombre fini de sous ensembles $\XX_1^1
,...,\XX_1^k$ tels que, pour $j=1,...,k$~:
\begin{enumerate}
\item  toutes les structures de $\XX_1^j$ possèdent une  même courbe de
découpage sur les faces de $\Delta^2 \setminus \Int (\bigcup_{1\leq
i\leq n} P_i )$\,;
\item toutes les structures $\xi \in \XX_1^j$ sont tangentes aux fibres des
faces des prismes de $P$\,;
\item pour toute arête $a$, toutes les structures $\xi \in \XX_1^k$ sont
égales sur un voisinage de $a\setminus \Int (\bigcup_{1\leq i\leq n}
P_i )$.
\end{enumerate}
\end{lemme}

\begin{proof}[Démonstration.]
Pour toute structure $\xi \in \XX_1$, le nombre de composantes du
découpage $\Gamma_{\Delta^2 \setminus \Int (\bigcup_{1\leq i\leq n}
P_i )} (\xi )$ est inférieur à $C_2$. Le nombre de classes
d'isotopie de multi-arcs
$$(\Gamma_{\Delta^2 \setminus \Int (\bigcup_{1\leq i\leq n} P_i )}
(\xi ))_{\xi \in \XX_1}$$
dans $\Delta^2 \setminus \Int (\bigcup_{1\leq i\leq n} P_i )$ est donc
fini. Quitte à déformer chaque structure $\xi \in \XX_1$ par une
isotopie relative à $\Lambda$ qui préserve $\Delta$ et à
partitionner $\XX_1$ en $X_1^1 ,...,X_1^k$, pour $j=1,...,k$, on peut
supposer que toutes les structures de $\XX_1^j$ possèdent une
même courbe de découpage sur $\Delta^2 \setminus \Int
(\bigcup_{1\leq i\leq n} P_i )$. Les propriétés 2 et 3
découlent du lemme~\ref{lemme : normalisation} (pour obtenir dans
3 l'égalité à partir du lemme~\ref{lemme : normalisation}, et
pas seulement une coïncidence des feuilletages legendriens, il
faut encore utiliser la contractibilité de $\mbox{Diff}^+([0,1])$
comme il est expliqué en détail dans la preuve du
lemme~\ref{lemme : ajuste}).
\end{proof}

On note $\XX_2$ l'un des $\XX_1^j$, $j=1,...,k$, et on fixe $\zeta \in
\XX_2$. On est ramené à démontrer le
théorème~\ref{t:domaines} pour $\XX_2$.

\begin{lemme}\label{l:normalisation-hexagone}
Toute structure $\xi \in \XX_2$ est isotope à une structure
maniable pour $\Delta$ et $\Lambda$-minimale $\xi'$ qui coïncide
avec $\zeta$ le long de $\Delta^2 \setminus \Int (\bigcup_{1\leq i\leq
n} P_i )$ et tangente aux fibres le long des faces verticales des
polyèdres de $P$.
\end{lemme}

On note $\XX_3$ l'ensemble obtenu par déformation de $\XX_2$. Le
long des portions de faces verticales qui ne sont pas incluses dans
$\Int (\bigcup_{1\leq i\leq n} P_i )$, toutes les structures $\xi \in
\XX_3$ sont à la fois égales et tangentes aux fibres.

\begin{proof}[Démonstration.]
Soit $R$ une composante de $H_F \setminus \Gamma_{H_F}$ incluse dans
$H_F \setminus \Int (\bigcup_{1\leq i\leq n} P_i)$. Les structures
$\xi$ et $\zeta$ sont égales sur un voisinage $K$ de $\partial
R\cap \partial F$ dans $R$, lequel $K$ est en outre feuilleté par
des arcs legendriens. On rétrécit $R$ en  $R_0$ en poussant
chaque arc de $\partial R \cap \partial F$ sur un arc legendrien qui
lui est parallèle dans $K$. Par application du lemme de
réalisation de feuilletages~\ref{l:feuilletages} à $R_0$, on
trouve une isotopie à support dans un voisinage de $R_0$, dont le
temps $1$ envoie $\xi$ sur $\zeta$ le long de $R_0$. Le point
important est que le support de cette isotopie ne rencontre ni
$\Lambda$ ni les faces autres que $F$, si bien $\Delta$ demeure
maniable et $\Lambda$-minimale pour l'image de $\xi$.
\end{proof}

\subsubsection{Normalisation des faces horizontales}

Pour chaque prisme $P_i$, $1\leq i\leq n$, on fixe un feuilletage
non singulier $\FF_i$ sur ses faces horizontales $Y_i \times \{
0,1\}$, égal au germe de feuilletage tracé par $\zeta$ au bord
$(\partial Y_i )\times \{ 0,1\}$ et transversal à une direction
fixée de $Y_i$ (on munit $Y_i$ d'une structure affine). Un tel
feuilletage non singulier existe grâce au fait que l'holonomie
vaut $-1$ et donc que l'indice du germe de  feuilletage
caractéristique $\zeta Y_i \times \{ j\}$ le long de $(\partial
Y_i) \times \{ j\}$, $j=0,1$ vaut $0$. Il détermine un germe de
structure $\eta_0$ le long des faces horizontales.

Grâce à la remarque~\ref{remarque : minoration}, on se place
dans la situation où, pour tout $\xi \in \XX_3$, toute face
verticale des prismes de la configuration $P$ contient au moins $20$
pièces.

\begin{lemme}
Il existe une  rétraction compacte $K= \bigcup_{1\leq i\leq n} P_i
\setminus N( \partial (\bigcup_{1\leq i\leq n} P_i ))$  de $\bigcup_{1\leq
i\leq n} P_i $, une structure de contact $\zeta_0$ sur $V\setminus
\Int (K)$ et pour tout $\xi \in \XX_3$ une isotopie de $\xi$ en
$\xi'$ qui vérifie~:
\begin{enumerate}
\item $\xi' =\zeta_0$ sur $V\setminus \Int (K)$\,;
\item $\xi'$ est tangente aux fibres verticales de chaque $P_i$.
\end{enumerate}
\end{lemme}
\begin{proof}[Démonstration.]
On commence par démontrer le lemme suivant~:

\begin{lemme}\label{l:int}
Soit $f : Y\times [-2,2] \rightarrow (V, \xi )$ un plongement dans
une variété de contact tendue. On suppose que~:
\begin{itemize}
\item chaque arc $\lbrace x\rbrace \times [-2,2]$, pour $x$
dans un voisinage de $\partial Y$, est legendrien pour $f^* \xi$\,;
\item sur chaque face verticale, les courbes
de découpage vont d'une arête verticale à l'autre\,;
\item l'holonomie de toute courbe d'holonomie tournant  autour
du bord vertical vaut $-1$\,; les arêtes de  $f( (\partial Y)
\times \{ \pm 1\} )$ sont transversales à $\xi$ \,;
\item il y a au moins quatre  courbes de singularités
sur chaque face verticale aux altitudes supérieures à $1$ et
inférieures à $-1$\,;
\item les arêtes horizontales de $Y\times [-2,2]$
sont transversales à $f^* \xi$.
\end{itemize}
Pour tout feuilletage non singulier non singulier $\FF$ tracé sur
les faces horizontales de $Y\times [-1,1]$ qui est tangent à $f^*
\xi$ au bord et transversal à une direction donnée de $Y$ (muni
de sa structure affine), il existe une isotopie de $\xi$ relative au
bord de l'image de $f$ en une structure $\xi'$, telle que $f^* \xi'$
trace le feuilletage $\FF$ sur $Y\times \lbrace \pm 1\rbrace$ et
que chaque arc $\lbrace x\rbrace \times [-1,1]$ soit legendrien.
\end{lemme}

\begin{proof}[Démonstration.]
On munit $Y\times [-2,2]$ de la structure $f^* \xi$. Grâce à la
présence de quatre arcs singuliers sur chaque face en dessous et
au-dessus des altitudes $-1$ et $1$, on construit sur $(\partial Y
)\times [1,2]$ et $(\partial Y )\times [-2,-1]$ une courbe
d'holonomie legendrienne $\gamma_{\pm 1}$ par concaté\-nation de
portions d'arêtes verticales et d'arcs  de singularités sur les
faces. On obtient que  $tb(\gamma_{\pm 1} )=-1$ car l'holonomie vaut
$-1$. On prend alors un disque convexe (de bord lisse par morceaux)
$D_{\pm 1} \subset Y\times [-2,2]$ qui s'appuie sur $\gamma_{\pm
1}$ transversalement au bord vertical de $Y \times [-2,2]$. La
portion de bord vertical comprise entre $\gamma_{-1}$ et $\gamma_1$
complétée par $D_{-1}$ et $D_1$ forme une sphère $S$ plongée
dans $Y \times [-2,2]$.

On va construire un modèle de cette sphère dans $\R^3 =\{
(x,y,t)\}$. Pour cela, on munit $\R^3$ de la structure $\eta$
d'équation $\cos tdx -\sin tdy=0$. On fixe un triangle (ou un
quadrilatère) $Y'$ dans le plan des $(x,y)$, et on considère le
polyèdre $Q=\lbrace (x,y)\in Y' , t\in [0 ,2k\pi ] \rbrace$.

Comme l'holonomie  de $\xi$ autour de $(\partial Y) \times [-2,2]$
est $-1$, il existe $k\in \R$ et $Y'$ pour que le bord vertical de
ce polyèdre $Q$, avec le germe de structure $\eta$, soit
conjugué à $(\partial Y) \times [-2,2]$ avec le germe de
structure $f^* \xi$ par un germe de difféomorphisme de contact
fibré (qui préserve la direction verticale, celle de $Q$ est
dirigée par $\partial_t$) $$\phi:(\partial
Y)\times[-2,2]\rightarrow (\partial Y')\times[0,2k\pi].$$

Il existe alors deux disques $D$ et $D'$ qui s'appuient
respectivement sur les courbes  $\phi ((\partial Y) \times \{ \pm
1\} )$ du bord vertical du polyèdre $Q$, transversaux à la
direction verticale $\partial_t$ et tels que $\phi$ s'étende en un
difféomorphisme $\Phi : Y \times [-1,1] \rightarrow Q'$, où
$Q'$ est le polyèdre de $\R^3$, dont le bord est union de $D$ et
$D'$ avec la portion de bord vertical de $Q$ comprise entre les
courbes $\phi ((\partial Y) \times \{ \pm 1\} )$, avec les
propriétés suivantes~:
\begin{itemize}
\item les fibres verticales de $Y\times [-1,1]$
sont envoyées par $\Phi$ sur les fibres verticales de $Q'$
(dirigées par $\partial_t$ qui sont legendriennes)\,;
\item $\Phi_* \FF$ est le feuilletage $\zeta (D\cup D')$.
\end{itemize}
Pour cela, on étend d'abord $\phi$ en un plongement fibré $\psi
:Y \times [-1,1] \rightarrow Q$. Il ne reste plus qu'à relever
les feuilletages $\psi_* \FF$ de manière legendrienne
transversalement à $\partial_t$ pour obtenir $D$ et $D'$, puis à
composer $\psi$ au but par l'extémité d'une isotopie qui
amène, en glissant le long des fibres, l'image par $\psi$ des
faces horizontales sur $D$ et $D'$. Le fait que $\FF$ soit
transversal à une même direction sur les faces supérieure et
inférieure assure que $D$ et $D'$ ne se rencontrent pas.

Les courbes $\phi (\gamma_{\pm 1} )$ sont legendriennes et
d'invariant de Thurston-Bennequin $-1$. On leur fait border dans
$(\R^3 ,\eta )$ deux disques convexes qui ne rencontrent pas $Q'$.
On note $S'$ la sphère construite comme précédemment à
partir de ces disques convexes et de la portion de bord vertical
qu'ils délimitent dans $Q$. D'après le lemme de réalisation de
feuilletages, on peut choisir ces deux disques de sorte que $\phi$
s'étende en un difféomorphisme
$\phi_1$ de $(S,f^* \xi \vert_S )\rightarrow (S',\eta \vert_{S'}
)$.

D'après le théorème~\ref{t:boule} d'Eliashberg, $\phi_1$
s'étend alors en un difféomorphisme de contact $\phi_2$ entre
les boules bordées par ces sphères.

Le plongement $\phi_2^{-1} \circ \Phi : Y \times [-1,1] \rightarrow
Y \times [-2,2]$ envoie $\FF$ sur le feuilletage caractéristique
de l'image des faces horizontales, et les fibres sur des courbes
legendriennes. Il vaut l'identité sur le bord vertical et
préserve l'orientation~: il est isotope à l'identité
relativement au bord de $Y \times [-2,2]$ par la restriction d'une
isotopie de $Y \times [-2,2]$.

On considère l'image de $f^* \xi$ par cette isotopie que l'on
propage à l'aide de $f$ dans $V$. Cette nouvelle structure donne
le résultat souhaité.
\end{proof}

Pour pousuivre, on applique le lemme~\ref{l:int} à chaque
structure $\xi \in \XX_3$ et au germe de structure $\eta_0$ le long
des faces horizontales~: il existe une isotopie de $\xi$ relative
à $\Delta$ qui envoie $\xi$ sur une structure $\xi'$ tangente aux
fibres de $P$ et égale au germe $\eta_0$ sur les faces
horizontales. On remplace  chaque structure $\xi \in \XX_3$ par la
structure, à nouveau notée $\xi$, obtenue après cette
première isotopie et on pioche une structure $\zeta_0$ dans
$\XX_3$.

Le long d'une fibre $I$ de $P_i$, toute structure $\xi$ tangente aux
fibres est repérée par une fonction angle $\theta_I : I
\rightarrow \R /\Z$ dans le fibré normal aux fibres, dont la
dérivée est strictement positive. Dans un petit voisinage  $N(
\partial (\bigcup_{1\leq i\leq n} P_i )) \cap (\bigcup_{1\leq i\leq n} P_i
)$ de $\partial (\bigcup_{1\leq i\leq n} P_i )$ dans $\bigcup_{1\leq i\leq
n} P_i$, on peut après isotopie (toujours par contractibilité de
$\mbox{Diff}^+ (I)$) rendre toutes ces fonctions angle égales à
celle de $\zeta_0 \in \XX_3$. On réalise une telle isotopie sur
chaque $\xi \in \XX_3$ pour obtenir $\xi' \in \XX'_3$.

Les prismes $P_i$ de la famille $P$ découpent chaque simplexe $G$
de $\Delta$ en polyèdres, homéomorphes à la boule, au bord
desquels toutes les structures  $\xi' \in \XX'_3$ coïncident.
Comme les structures de $\XX'_3$ sont tendues, le théorème
d'unicité~\ref{t:boule} d'Eliashberg donne, sur chaque polyèdre,
une isotopie stationnaire au bord entre une structure quelconque
$\xi' \in \XX'_3$ et $\zeta_0$.
\end{proof}

\begin{corollaire}\label{cor : lissage}
Il existe un domaine fibré à bord $(M,\tau )$ avec $K \subset M
\subset \bigcup_{1\leq i\leq n} P_i $, dont la fibration en intervalles
$\tau$  est la restriction de la fibration de $\bigcup_{1\leq i\leq n}
P_i $. En particulier, selon la terminologie de la sous-section
suivante, toutes les structures de $\XX_3$ sont ajustées à
$(M,\tau ,\zeta )$.
\end{corollaire}

\begin{proof}[Démonstration.]
On retire à $\bigcup_{1\leq i\leq n} P_i$ un petit voisinage fibré
$U$, $U\cap K =\emptyset$, des portions d'arêtes qui ne sont pas
dans $\Int (\bigcup_{1\leq i\leq n} P_i )$. On note $\tau'$ la
fibration en intervalles de $\bigcup_{1\leq i\leq n} P_i $. L'espace
quotient $\Sigma =(\bigcup_{1\leq i\leq n} P_i \setminus U)/\tau'$ est
une surface branchée dont le lieu singulier est non générique.
Un lissage  du bord de $\bigcup_{1\leq i\leq n} P_i \setminus U$
transversal à la fibration au-dessus des points réguliers de la
projection et  une petite modification de $\bigcup_{1\leq i\leq n} P_i
\setminus U$ correspondant à une perturbation générique de
$\Sigma$ pour obtenir un lieu de branchement générique permet
d'obtenir le  domaine fibré recherché.
\end{proof}

Pour démontrer le théorème~\ref{t:domaines}, il reste une
difficulté~: le domaine  fibré du corollaire~\ref{cor : lissage}
est \guil{à bord}. Dans la suite, on indique  comment on se
ramène à des voisinages fibrés de surfaces branchées sans
bord.

\subsubsection{Structures de contact ajustées à un domaine fibré}

Soit $(M,\tau )$ un domaine fibré, voisinage d'une surface
branchée à bord, et  $X =M/\tau $ la surface branchée
quotient. On note $\pi :M \rightarrow X$ et $\zeta$ une structure de
contact sur $V\setminus \Int (M)$. Une structure de contact est {\it
ajustée} à $(M,\tau, \zeta )$ si elle est égale à $\zeta$
hors de $\Int (M)$ et tangente à $\tau$ dans $M$.

\begin{lemme}\label{lemme : ajuste}
Toute structure de contact $\xi$  ajustée à $(M,\tau ,\zeta )$
est déterminée, à isotopie près parmi les structures
ajustées, par la fonction $$a_\xi :\partial_h M:\rightarrow
]0,\infty[$$ qui est continue sur chaque secteur et  qui associe à
chaque point $p$ l'angle de rotation total de $\xi$ le long de la
feuille de $\tau$ partant de $p$. (On mesure cet angle avec une
métrique auxiliaire et l'holonomie de $\tau$.)
\end{lemme}

\begin{proof}[Démonstration.]
Ce lemme est une conséquence de la contractibilité de l'espace
$\mbox{Diff}^+(I)$ des difféomorphimes de l'intervalle qui
préservent l'orientation. Soit $\xi_0$ et $\xi_1$ deux
structures de contact sur $D^2 \times I$ muni des coordonnées
$((x,y),t)$, égales le long de $D^2 \times \{ 0,1\}$, tangentes
aux fibres $\{ p\} \times I$ et pour lesquelles $a_{\xi_0}
=a_{\xi_1}$ sur $D^2$. Elles sont alors données par $\alpha_0
=\cos f_0 (x,y,t)dx -\sin f_0 (x,y,t)dy$ et $\alpha_1 =\cos f_1
(x,y,t)dx -\sin f_1 (x,y,t)dy$ avec, pour $i=0,1$, $\frac{\partial
f_i}{\partial t} >0$ et $f_0 (x,y,i)=f_1 (x,y,i)$. Par la
contractibilité de $\mbox{Diff}^+ (I)$, il existe une famille à
un paramètre de fonctions $f_s : D^2 \times I \rightarrow \R$,
$s\in [0,1]$, qui vérifient $\frac{\partial f_s}{\partial t} >0$
et qui sont indépendantes de $s\in [0,1]$ sur $D^2 \times \{
0,1\}$. Le chemin de formes de contact $\alpha_s =\cos f_s (x,y,t)dx
-\sin f_s (x,y,t)dy$ donne une isotopie relative à $D^2 \times \{
0,1\}$ entre $\xi_0$ et $\xi_1$ parmi les structures tangentes à
$\{ *\} \times I$.

Pour démontrer le lemme, on donne une version relative de la
discussion précédente. Soit $\xi$ et $\xi'$  deux structures
ajustées à $(M,\tau ,\zeta )$, qui ont les mêmes fonctions
angles $a_\xi =a_{\xi'}$. Si $B$ est un secteur de $X$, alors
$\partial B$ est un polygone $\partial B=\delta_1 \cup \delta_2\cup
...\cup \delta_m$, où les $\delta_i$ sont les côtés successifs
de $B$ qui se rencontrent le long des points triples de $X$. Comme
$(B\times I) \cap \partial_v M$  consiste en la réunion de
produits $\delta_i \times [a_i ,b_i ]$ où $[a_i ,b_i] \subset
[0,1]$, les structures $\xi$ et $\xi'$ sont égales sur $\delta_i
\times [a_i ,b_i ]$. On déforme d'abord $\xi$ en $\xi'$ le long de
$\delta_i \times ([0,1] \setminus [a_i ,b_i ])$ en utilisant la
contractibilité de $\mbox{Diff}^+ (I)$, puis on effectue
l'isotopie entre $\xi$ et $\xi'$ sur $B\times I$ relativement à
$(\partial B)\times I$.
\end{proof}

Soit $\xi_0$ une structure  de contact
ajustée à $(M,\tau ,\zeta )$ fixée.

\begin{proposition}
Les classes d'isotopie des structures de contact $\xi$  ajustées à
$(M,\tau ,\zeta )$  sont en bijection avec les fonctions
$$w_\xi :\pi_0 (Reg (X))\rightarrow \Z,$$  dites fonctions \emph{poids},
vérifiant la condition
$$w_\xi (R)>-\frac{1}{2\pi} \inf(a_{\xi_0})$$  et les relations d'adjacence
$$w_\xi (R_1 )+w_\xi (R_2 )=w_\xi (R_3 )$$
pour toute paire de feuillets réguliers $R_1$, $R_2$ de $X$ qui se
joignent pour donner $R_3$.
\end{proposition}

\begin{proof}[Démonstration.]
Pour toute structure de contact $\xi$ ajustée à
$(M,\tau,\zeta)$, on définit
$$w_\xi =\frac{1}{2\pi} (a_\xi -a_{\xi_0} ).$$
Le poids $w_\xi (p)$ est un entier qui varie continûment sur les
composantes de $Reg (X)$. Il est donc constant sur $Reg (X)$. La
relation d'adjacence est donnée en regardant une fibre située
au-dessus d'un point $p$ du lieu singulier de $X$ où $R_1$ et
$R_2$ se rejoignent pour donner $R_3$. La condition de contact donne
que $a_\xi >0$, ce qui fournit l'inégalité  $w_\xi
(R)>-\frac{1}{2\pi} \inf (a_{\xi_0})$.

Réciproquement, si on dispose d'une telle fonction poids, la
fonction $a_\xi =w_\xi + a_{\xi_0}$ détermine une structure de
contact qui convient.
\end{proof}

\subsubsection{Le lemme d'élagage}

\begin{lemme}{\rm [Lemme d'élagage]} \label{l:elagage}
Soit $(M,\tau,\zeta)$ un domaine fibré et $\XX$ un ensemble de structures de 
contact ajustées à $(M,\tau,\zeta)$. On suppose qu'il existe un réel $C$ et un 
point $p \in \partial_h M$ tels que $a_\xi(p) < C$ pour tout $\xi \in \XX$ et on
note $X_1, \dots, X_k$ les adhérences des strates régulieres de $X = M/\tau$ qui 
contiennent $\pi(p)$. On peut alors trouver des structures de contact $\zeta_1, 
\dots, \zeta_l$ sur le complémentaire du domaine fibré
$$ (M',\tau') = \bigl( M \setminus \Int (\bigcup_{j=1}^k \pi^{-1}(X_j)),
   \tau \vert_{M'} \bigr) $$
telles que toute structure $\xi \in \XX$ soit isotope à une structure ajustée à 
l'un des $(M',\tau',\zeta_i)$.
\end{lemme}

\begin{proof}
La borne sur $a_\xi$ donne une borne sur tous les poids $w_\xi (X_i
)$, $1\leq i\leq k$. Ceux-ci peuvent donc prendre seulement un
nombre fini de valeurs. On peut partitionner $\XX$ en
$\XX_1$,...,$\XX_l$, de sorte que toutes  les $\xi \in \XX_j$ donnent
le même poids. Si $\zeta_i$ est une structure de $\XX_i$, on peut
alors, comme dans le lemme~\ref{lemme : ajuste}, trouver une
isotopie, parmi les structures ajustées à $(M,\tau ,\zeta )$,
entre toute structure $\xi \in \XX_i$ et une structure $\xi'$ égale
à $\zeta_i$ sur $\pi^{-1} (\bigcup_{1\leq j\leq k} X_j )$.
\end{proof}

L'opération d'élagage diminue strictement le nombre de secteurs de $X$. En la 
répétant un nombre fini de fois, on obtient donc l'ensemble vide.

\begin{corollaire} \label{c:bord}
Si $\XX$ est un ensemble de structures de structures de contact
ajustées à $(M,\tau ,\zeta )$, alors il existe $(M_1 ,\tau_1
,\zeta_1 )$,..., $(M_l ,\tau_l ,\zeta_l )$ obtenues par élagage de
$(M,\tau ,\zeta )$ telles que, pour $1\leq i\leq l$, $X_i =M_i
/\tau_i$ soit une surface branchée sans bord, et que toute
structure $\xi \in \XX$ soit isotope  à une structure $\xi'$
ajustée à l'un des $(M_i,\tau_i ,\zeta_i )$.
\end{corollaire}

\begin{proof}[Démonstration.]
Tout point $p \in \partial X$ vérifie les hypothèses du lemme
d'élagage. Par ailleurs, un élagage réduit strictement le
nombre de secteurs de $X$. Un nombre fini d'applications du lemme
d'élagage~\ref{l:elagage} conduit au résultat.
\end{proof}

Le corollaire~\ref{c:bord} termine la démonstration du théorème~\ref{t:domaines}.

\section{Finitude homotopique}

On démontre ici le théorème~\ref{t:generation}. Le théorème~\ref{t:homotopie} en
résulte car une modification de Lutz à coefficient entier le long d'un tore ne 
change pas la classe d'homotopie dans les champs de plans.

Soit $V$ une variété close. Vu le théorème~\ref{t:domaines}, il suffit de 
démontrer le théorème~\ref{t:generation} pour un ensemble $\SS$ de structures de
contact tendues toutes ajustées à un même domaine fibré $(M,\tau,\zeta)$. Soit
$X$ la surface branchée $M/\tau$ avec son lieu singulier $\Theta$ et ses strates
régulières $X_0, \dots, X_d$. Étant donné $\xi_0 \in \SS$, chaque structure $\xi
\in \SS$ est déterminée (à déformation près parmi les structures ajustées) par 
son poids $w_\xi = (w_\xi (X_1), \dots, w_\xi(X_d)) \in \Z^d$ relatif à $\xi_0$.
Par application du lemme d'élagage et quitte à effectuer une partition de $\SS$,
on se ramène au cas où, pour tout $\xi \in \SS$, le poids $w_\xi$ appartient à
$\N^d$. À chaque arête lisse $C$ de $\Theta$, on associe l'équation linéaire 
$x_i = x_j+x_k$ sur $\R^d$ où $X_j$ et $X_k$ sont les secteurs de $X$ qui se 
joignent le long de $C$ pour former $X_i$. Soit $W$ le sous-espace vectoriel de 
$\R^d$ des solutions de ce système. Chaque poids $w_\xi$ appartient à $W \cap 
\N^d$.  

On note $\preceq$ l'ordre partiel défini sur $\Z^d$ par
$$ (x_1, \dots, x_d) \preceq (y_1, \dots, y_d) \quad \text{si} \quad
   x_i \leq y_i \text{\ \ pour\ \ } 1 \le i \le d. $$

\begin{lemme}\label{l:minimaux} Soit $W$ un sous-espace vectoriel de $\R^d$.
Les éléments minimaux de $W \cap \N^d$ pour l'ordre partiel $\preceq$ 
sont en nombre fini et engendrent $W \cap
\N^d$.
\end{lemme} 
\begin{proof}[D\'emonstration.]  
Pour montrer que les \'el\'ements minimaux de $W\cap \N^d$ sont en nombre
fini, on raisonne par r\'ecurrence sur $d$.
Lorsque $d=1$, on a un unique \'el\'ement minimal.
On suppose le r\'esultat d\'emontr\'e pour tout sous-espace de 
$\R^{d-1}$, $d-1 \geq 1$.
Soit maintenant $W$ un sous-espace de $\R^d$, $d\geq 2$. On raisonne par l'absurde en
supposant que $W\cap \N^d$ poss\`ede une infinit\'e d'\'el\'ements 
minimaux. Quitte a permuter l'ordre des coordonn\'ees,
on peut alors trouver une suite  d'\'el\'ements 
minimaux $(v^i )_{i\in \N}$, $v^i =(v^i_1 ,\dots ,v^i_d )$, de $W\cap \N^d$
pour laquelle la suite des derni\`eres coordonn\'ees  $(v^i_d )_{i\in \N}$ tend vers
l'infini. On note $\pi :\R^d \rightarrow \R^{d-1}$
la projection sur les $d-1$ premi\`eres coordonn\'ees.
En  appliquant l'hypoth\`ese de r\'ecurrence
\`a $\pi (W )$, on obtient un nombre fini
d'\'el\'ements $u_1 , \dots ,u_n $ de $\N^d$ dont les images
par $\pi$ sont les \'el\'ements minimaux de $\pi (W) \cap \N^{d-1}$.
Si $i\in \N$ est assez grand, $v^i_d$ est sup\'erieur \`a chacune
des derni\`eres coordonn\'ees des vecteurs $u_1, \dots ,u_n$ et en particulier,
$v^i$ est sup\'erieur \`a un des $u_j$. C'est une contradiction.

Si $u \in W\cap \N^d$ n'est pas un \'el\'ement minimal, il est sup\'erieur 
\`a l'un d'entre eux $u_i$. Soit le vecteur $u-u_i$ est minimal, soit on peut \`a nouveau
lui retirer un des $u_j$, $j=1,\dots n$. En un nombre fini d'\'etape, on \'ecrit $u$
comme une somme d'\'el\'ements minimaux.
\end{proof}

On note $u_1, \dots, u_k$ les \'el\'ements minimaux de $W\cap \N^d$
donn\'es par le lemme~\ref{l:minimaux}. Par ailleurs,

\begin{lemme} \label{porte}
Les classes d'isotopie des surfaces compactes plongées dans $M$ transversalement
à $\tau$ sont en bijection avec les éléments non nuls de $W \cap \N^d$.
\end{lemme}

Si $\TT$ est une surface portée par $(M,\tau)$, on note $w_\TT = (w_\TT (X_1), \dots, 
w_\TT (X_d))$ son poids, où $w_\TT (X_i)$ est le nombre de composantes de $\TT \cap 
\pi^{-1}(X_i)$ et $\pi$ la projection $M \to X = M/\tau$.

\begin{remarque}
Si $X$ est compact sans bord et $\SS$ non vide, les surfaces compactes portées par $X$ 
--~c'est-à-dire plongées dans $M$ transversalement à $\tau$~-- sont closes et 
transversales aux structures de contact de $\SS$. Leurs composantes connexes 
sont donc des tores et/ou des bouteilles de Klein.
\end{remarque}

Soit $T_1, \dots, T_k$ les surfaces portées par $(M,\tau)$ et correspondant
aux éléments $u_1 ,\dots ,u_k$ de $W \cap \N^d$. On suppose que $T_1,\dots ,T_l$
sont des tores et que $T_{l+1} ,\dots ,T_k$ sont des bouteilles de Klein.
Lorsque $T_i$ est une bouteille de Klein, on note 
$T_i'$ le tore, toujours port\'e par $(M,\tau )$, bordant un de ses petits
voisinages tubulaires $N(T_i)$. Le poids $u_i'$ de $T_i'$ vaut $2u_i$.

\begin{lemme} \label{Lutz}
Toute structure de contact $\xi \in \SS$ est obtenue à partir de $\xi_0$ par 
modification de Lutz sur les tores $T_i$, $1\leq i\leq l$ et 
$T_i'$, $l+1\leq i\leq k$.
\end{lemme}

\begin{proof}
Soit $\xi \in \SS$. On a $w_\xi = \sum_{i=1}^k n_i(\xi) u_i$, avec $n_i (\xi ) \in
\N$. La structure obtenue à partir de $\xi_0$ par chirurgie de Lutz de 
coefficient $n_i (\xi )$ sur $T_i$ lorsque $T_i$ est un tore
et de coefficient $\frac{1}{2} n_i (\xi )$ sur  $T_i' =\partial N(T_i )$
lorsque $T_i$ est une bouteille de Klein est ajustée à $(M,\tau,\zeta)$ et a le même poids 
que $\xi$. Elle lui est donc isotope relativement à $\partial M$.
\end{proof}

M\^eme si les coefficients de chirurgie apparaissant dans le lemme~\ref{Lutz}
sont des demi-entiers, celui-ci implique  
le th\'eor\`eme~\ref{t:generation}~: on ajoute comme \guil{structures de base}
\`a $\xi_0$ les $2^{k-l-1}$ structures  obtenues en faisant ou non 
des modifications de Lutz de coefficient
$\frac{1}{2}$ le long des tores $T_i'$, $l+1\leq i\leq k$. Toutes
les structures de $\SS$ sont alors obtenues 
\`a partir de l'une d'entre elles par modifications de
Lutz \`a coefficients entiers sur les collections de tores 
$(T_i )_{1\leq i\leq l}$ et 
$(T_i' )_{l+1\leq i\leq k}$.

\section{Finitude géométrique}

On démontre maintenant le théorème~\ref{t:torsion}. Comme précédemment, grâce au 
théorème~\ref{t:domaines}, il suffit d'établir le résultat pour un ensemble 
$\SS$ de structures tendues toutes ajustées à un même moule $(M,\tau,\zeta)$. On
fixe une structure de référence $\xi_0 \in \SS$, toute autre $\xi \in \SS$ étant
alors déterminée par son poids $w_\xi$ dont on peut supposer, sauf pour $\xi_0$,
qu'il est dans $(\N \setminus \{0\})^d$ 
(voir le lemme d'\'elagage~\ref{l:elagage}). Pour une telle structure  $\xi \in \SS 
\setminus \{\xi_0\}$, le poids $w_\xi$ correspond à une surface close $T$ 
pleinement portée par $(M,\tau)$. Chaque composante de $T$ est un tore ou une 
bouteille de Klein.

\begin{lemme} \label{l:porte}
Il existe une surface orientée $\TT$ pleinement portée par $(M,\tau)$ qui contient
le bord horizontal de $M$ comme sous-surface orientée.
\end{lemme}

\begin{proof}
Le domaine $M$ contient une surface $T$ pleinement portée, par exemple celle que
détermine le poids d'une structure $\xi \in \SS$. En doublant $T$ si nécéssaire
et en remplaçant chaque bouteille de Klein par le bord d'un de ses voisinages
tubulaires, on s'assure que les composantes de $T$ sont des tores et que $T$ 
intersecte au moins deux fois chaque fibre de $M$. Parmi ces intersections, deux
sont les plus proches de $\partial_h M$. On pousse les composantes connexes de 
$T$ contenant ces intersections extriêmes jusque dans $\partial M$. Quitte à
doubler $T$ à nouveau pour pouvoir l'orienter comme voulu, on obtient la surface
cherchée.
\end{proof}

\subsection{Le lemme du degré}

Soit $A$ une composante connexe de $\partial_v M$ et $C$ une composante du bord 
de $A$. On définit le degré $\deg(A)$ de $A$ comme la valeur absolue du degré de
l'image de $\zeta_x$ par rapport à $T_x A$ dans le quotient $T_xV / T_x\tau$ 
lorsque $x$ parcourt $C$.

\begin{assertion} \label{a:degre}
Quitte à modifier $(M,\tau,\zeta)$ par élagage et isotopie, on peut supposer que
pour chaque composante $A$ de $\partial_vM$ et chaque composante $C$ de 
$\partial A$, il y a exactement $2\deg(A)$ points le long de $C$ où $T_xA$ 
coïncide avec $\zeta_x$. En particulier, si $\deg(A) = 0$, le feuilletage 
$\zeta A$ contient une courbe de singularités isotope à l'âme de $A$.
\end{assertion}

\begin{proof}
On étend $A$ le long des fibres legendriennes de $M$ jusqu'à obtenir un anneau 
$A'$, où $A' \setminus A \subset M$, pour lequel les structures $\xi \vert_{A'}$
sont égales pour tout $\xi \in \SS$ et la condition de l'assertion est vérifiée 
le long de $A'$. Cette extension est rendue possible par le fait qu'on peut 
supposer que toutes les structures $\xi \in \SS$ ont une fonction poids minorée
par une constante $c>0$ que l'on peut choisir arbitrairement grande (quitte à
procéder à des partitions de $\SS$ et à élaguer $M$). On considère alors un 
épaississement $A' \times [0,1]$ de $A'$ dans $M$, où $A' \times \{1\} = A'$ et 
où chaque $\xi \in \SS$ est $I$-invariante (en particulier chaque $A' \times 
\{t\}$ est feuilleté par des intervalles legendriens). On prend pour nouveau
domaine fibré $M' = M \setminus (A' \times I)$ dont on lisse les coins de sorte 
que $A' \times \{0\} \subset M'$ soit une composante du bord vertical de $M'$.
\end{proof}

Dorénavant, toutes les composantes de $\partial_v M$ sont réputées satisfaire 
aux conclusions de l'assertion~\ref{a:degre}.

\begin{lemme} \label{lemme : degre}{\rm [Lemme du degré]}
Il existe des domaines fibrés $(M_i ,\tau_i ,\zeta_i )_{1\leq
i\leq k}$ obtenus par élagage de $(M,\tau ,\zeta )$ avec les
propriétés suivantes~:
\begin{itemize}
\item  toute structure
de contact $\xi \in \SS$ est conjuguée, par un produit de twists de
Dehn et d'isotopies,  à une structure ajustée à l'un des $(M_i
,\tau_i ,\zeta_i )$\,;
\item tout domaine fibré $(M_i ,\tau_i )$ porte pleinement une surface
$\TT_i$ comme dans le lemme~\ref{l:porte}\,;
\item toutes les composantes de $\partial_v M_i$ satisfont aux conclusions
de l'assertion~\ref{a:degre}\;
\item toute composante de $\partial_v M_i$ de degré
non nul  intersecte $\TT_i$ le long de courbes contractibles
dans $\TT_i$.
\end{itemize}
\end{lemme}

\begin{proof}[Démonstration.]
Soit $T$ une composante de  $\mathcal{T}$ qui intersecte une
composante $A$ de  $\partial_v M$ avec  $\deg(A)\not=0$ le long de
$c$. On suppose que $c$ est non contractible sur  $T$. Pour tout
$\xi\in \SS$ avec $w_\xi$ suffisamment grande, il existe un
plongement $\phi:T\times[0,1]\rightarrow M$, où  $\phi(T,0)=T$ et
$\phi^* \xi$ est donnée par  $\cos (f(x,y)+2\pi t)dx -\sin
(f(x,y)+2\pi t) dy=0$. Ici, les coordonnées sur
$T\times[0,1]=\R^2/\Z^2\times[0,1]$ sont  $(x,y,t)$, et  $f$ est une
fonction à valeurs dans le cercle  $T\rightarrow \R/2\pi\Z$. On a
alors la propriété suivante~:

\begin{assertion}\label{perturb}
Il existe un tore  $T'\subset \phi(T\times[0,1])$ isotope à $T$ et
transversal aux fibres legendriennes tel que $T'$ soit convexe et
$\#\Gamma_{T'} (\xi )\leq 2 \deg(A)$.
\end{assertion}

\begin{proof}[Démonstration de l'assertion~\ref{perturb}]
Après une perturbation  $\classe\infty$-petite de $T$, on peut
supposer que $T$ est convexe. Comme  $T\pitchfork \xi$ pour tout
$\xi\in \SS$, le feuilletage caractéristique  $\xi T$ est non
singulier, et donc $\#\Gamma_T (\xi )$ est égal au nombre
d'orbites fermées $\gamma_i$ de  $\xi T$.  Les $T\setminus \bigcup_i
\gamma_i$ sont des composantes annulaires qui sont soit de Reeb (ne
possèdent pas d'arc transversal s'appuyant au bord et qui
intersecte toutes les feuilles) ou tendues (il existe un tel arc
transversal). On peut supposer que  $c$ est transversal à
$\bigcup_i\gamma_i$. En analysant les composantes de $c\setminus
\bigcup_i\gamma_i$, tout arc (séparant ou non séparant) situé
dans une composante de Reeb donne au moins une tangence, tandis que
les arcs situés dans une composante tendue ne contribuent pas
nécéssairement. Donc le nombre de composantes de Reeb est
majoré  par $2\deg (A)$ ($=$ le nombre de tangences de  $c$), si
les orbites $\gamma_i$ ont une intersection non triviale avec  $c$.
Pour voir que les orbites  $\gamma_i$ ont une intersection
géométrique non triviale avec  $c$, on observe  que $2\deg(A)$,
qui est le décompte algébrique du nombre de tangences entre  $c$
et  $\xi T$, est invariant par isotopie. Si le nombre d'intersection
géométrique est nul, alors le degré doit aussi être nul.
Pour conclure, toutes les composantes tendues peuvent être
supprimées en isotopant $T$ à une distance bornée dans
$\phi(T\times[0,1])$. On peut également remarquer que
l'assertion~\ref{a:degre} implique que les composantes de Reeb
pointent dans la même direction.
\end{proof}

Le point clé pour le tore convexe $T$ modifié comme dans le
lemme précédent est que  $\#\Gamma_T (\xi )$ est borné
indépendemment du choix de $\xi \in \SS$.  On suppose que  $\xi\in
\SS$ satisfait $w_\xi \gg nw_T$, où  $n={1\over 2} \#\Gamma_T (\xi
)$. Alors il existe un plongement $\psi: T\times[0,n]\rightarrow
N(\mathcal{B})$, où $T\times\{0\}=T'$ et $\psi^*\xi$ est donnée
par  $\cos (g(x,y)+2\pi t)dx -\sin (g(x,y) +2\pi t) dy=0$.  Si on
excise $\psi(T\times[0,n])$ et on recolle $\psi(T\times\{0\})$ avec
$\psi(T\times\{n\})$ via l'identification naturelle donnée par la
fibration legendrienne, on obtient une structure de contact  $\xi'$
correspondant au poids $w_\xi - nw_T$. Maintenant, $\xi$ et  $\xi'$
sont isomorphes car elles diffèrent par des twists de Dehn le long
du tore $T$.   Pour cette raison, on peut réduire inductivement
$w_\xi \rightarrow w_\xi - nw_T$ jusqu'à ce qu'un secteur de $X$
ait un petit poids. Un tel secteur peut alors être élagué.
\end{proof}

En appliquant le lemme du degré à $(M,\tau ,\zeta )$, on se
ramène au cas où $(M,\tau ,\zeta )$ est l'un des $(M_i ,\tau_i
,\zeta_i )$. On peut encore améliorer $M$ par élagage grâce
à  la proposition suivante~:

\begin{proposition}\label{improved degree lemma}
Après des élagages successifs, on peut supposer, quitte à
conjuguer les structures de $\SS$ par des produits de twists de Dehn,
que, pour toute composante $A$ de $\partial_v M$,
\begin{enumerate}
\item $\deg (A)=0$ si et seulement si les deux composantes de $\partial A$
sont non contractibles dans  $\mathcal{T}$.
\item $\deg (A)=1$ si et seulement si les deux composantes de $\partial A$
bordent des disques dans $\TT$.
\end{enumerate}
\end{proposition}

\begin{proof}[Démonstration.]
Par le lemme du degré, si une composante de $\partial A$ est non
contractible, alors $\deg(A)=0$. À l'inverse, si une composante
$c$ de $\partial A$ borde un disque $D$ dans $\TT$, alors par
l'assertion~\ref{a:degre} et la non singularité du feuilletage
caractéristique de $D$, il ne peut y avoir que deux points le long
de $c$ pour lesquels $T_xA=\zeta_x$.  Donc  $\deg (A)=1$. Ainsi,
soit les deux composantes de $\partial A$ sont non contractibles
dans $\TT$, soit elles bordent toutes deux un disque.
\end{proof}

\begin{remarque}
Si les deux composantes de $\partial A$ bordent un disque, alors ces
disques doivent être du même côté de $A$, sinon la sphère
constituée de l'union de ces deux disques avec $A$ pourrait être
lissée en une sphère transversale à $\xi_0 \in \SS$.
\end{remarque}

\subsection{\'Elimination des disques de contact}

Dans cette partie, on simplifie le domaine fibré $M$ en
éliminant les disques de contact. Un \emph{disque de contact} est
un disque proprement plongé $D\subset M$, transversal aux fibres
et dont le bord est dans $\partial_v M$.

\begin{lemme}\label{replace}
Soit $A$ une composante de $\partial_v M$. S'il existe un disque de
contact $D$ dont le bord est dans $A$, alors les composantes $c_1$,
$c_2$ de $\partial A$ bordent des disques  $D_1$, $D_2\subset
\mathcal{T}$ de sorte que $D_i$ soit dans l'intérieur de $M$
près de  $\partial D_i$.
\end{lemme}

\begin{proof}[Démonstration.]
Comme  $A$ admet un disque de contact $D$ et que le feuilletage
caractéris\-tique de $D$ est non singulier, $\deg(A)$ doit être
égal à un. Par le lemme du degré, $c_i$ doit border un disque
$D_i$ dans  $\TT$. On note que $D_i$ ne peut pas être dans la
\guil{direction opposée} à $D$, c'est-à-dire que $D_i$ ne peut
pas contenir la composante de $\partial_h M$ adjacente à $c_i$.
Dans le cas contraire, $D\cup D_i$ (augmenté de parties de $A$ et
après lissage) formerait une sphère immergée transversale à
la fibration legendrienne, ce qui est impossible.
\end{proof}

\begin{remarque}
Il est possible que  $D_1\subset D_2$ où  vice-versa.
\end{remarque}

\begin{proposition}\label{no disks of contact}
Soit  $V$ une variété close et irréductible de dimension
trois. Il existe un nombre fini de paires  $(N_i,\zeta_i)$,
$i=1,\dots,k$, satisfaisant aux conditions suivantes~:

\begin{enumerate}
\item $N_i\subset V$ est une union finie de tores épaissis
$T^2\times [0,1]$, d'anneaux épaissis $A\times [0,1]$ et de
voisinages de bouteilles de Klein $N(K)$, où chaque
$A\times\{j\}$, $j=0,1$, est collé de façon incompressible sur
un $\partial(T^2\times [0,1])$ ou sur un $\partial N(K)$, et où
certaines composantes de bord de  $T^2\times I$ peuvent être
identifiées entre elles ou avec un $\partial N(K)$\,;

\item $\zeta_i$ est une structure de contact tendue sur $V\setminus
\Int (N_i )$\,;

\item toute structure de contact tendue sur $V$ est conjuguée,
par un produit de twists de Dehn et d'isotopies à une structure
égale à l'une des $\zeta_i$ sur $V\setminus \Int (N_i )$.
\end{enumerate}
\end{proposition}

\begin{proof}
On élimine d'abord les disques de contact de $M$, tout en
préservant la condition que $M$ porte pleinement une union de
tores $\mathcal{T}$. (Voir la  remarque~\ref{loss} ci-dessous.) S'il
y a un disque de contact pour $A$, alors en utilisant le
lemme~\ref{replace}, on peut le remplacer par des disques de contact
(pour $A$) $D_1$ et $D_2$  dans $\TT$. Sans perte de
généralité, on suppose que $D_1$ est un disque de contact le
plus intérieur pour $\TT$. Alors, soit
 $D_1$ et $D_2$ sont disjoints, soit $D_1\subset D_2$.
Comme $D_1$ peut contenir des disques de $\partial_h M$, on le
pousse légèrement le long des fibres pour l'éloigner de
$\partial_h M\cap \Int (D_1 )$. On appelle $D_1'$ cette déformation
de $D_1$. On change alors $M$ en $M\setminus D_1'$, $\TT$ en $(\TT
\setminus D_1 ) \cup D_1'$, et $D_2$ en  $ (D_2\setminus D_1)\cup
D_1'$ si  $D_1\subset D_2$ (puis, après le paragraphe ci-dessous
on les rebaptise $M$, $\mathcal{T}$ et  $D_2$).

On explique maintenant comment transformer  $D_1'$ en une surface
convexe (de bord legendrien), de sorte que la structure de contact
$\zeta$ s'étende de manière unique en une structure de contact
tendue sur $M\cup N(D_1')$, c'est-à-dire de sorte qu'on puisse
isotoper toutes les structures $\xi \in \SS$ sur $M$ relativement à
$\partial M$ pour les faire coïncider sur (un voisinage de)
$D_1'$. Après un lissage des coins de $\partial M$ et une
perturbation générique, $\partial M$ devient convexe. Le fait
que   $\deg(A)=1$ se traduit par le fait qu'on peut trouver une
isotopie de $\zeta$ près de $\partial M$ (et donc une isotopie
concomitante des structures $\xi \in \SS$) qui rende $\partial D_1'$
legendrien avec $tb(\partial D_1')=-1$.  Dans ce cas, si $D_1'$ est
perturbé en une surface convexe de bord legendrien, il n'y a
qu'une seule possibilité pour  $\Gamma_{D_1'} (\xi )$ à isotopie
près.  Dès lors, en appliquant le lemme de réalisation de
feuilletages on peut supposer que  toutes les structures $\xi \in
\SS$ donnent le même germe le long de $D_1'$.

Comme (le nouveau) $\mathcal{T}$ n'est pas pleinement porté par le
(nouveau) $M$ (voir la remarque~\ref{loss}~: la fibration persiste
topologiquement, de même que la notion de surface portée), on
modifie $\TT$ de la façon suivante~: soit $T$ la composante de
$\TT$ qui contient $D_2$, et soit $T'$ une déformation parallèle
de $T$. On remplace $T$ par un lissage de $(T\setminus D_2)\cup
A\cup D_1'$. Si on dédouble ce nouveau tore, on obtient une
nouvelle collection de surfaces $\TT$ qui contient $\partial_h M$.
Comme $\partial_v M$ est constitué d'un nombre fini de
composantes, on élimine ainsi tous les disques de contact en un
nombre fini d'opérations. On observe que les composantes de
$\partial_v M \cap \TT$ qui étaient non  homotope à zéro
(resp. homotope) dans $\TT$ le demeurent.

Une fois éliminés tous les disques de contact, on étudie les
composantes de $\partial_h M$. Elles sont de trois types~: (i)
disques, (ii) anneaux incompressibles dans $\TT$ et (iii) tores.
Tous les disques de $\partial_h M$ peuvent être éliminés de la
façon suivante. Soit $D$ une composante de $\partial_h M$
difféomorphe à un disque et $A$ un anneau qui partage une
composante de bord avec $D$. Le degr\'e $\deg(A)$ doit être non nul et
si $S$ est une composante de $\partial_h M$ qui intersecte l'autre
composante de bord de $A$, alors, par le lemme du degré, $S$ ne
peut pas être un anneau incompressible dans $\TT$. Donc $S$ est
aussi un disque. Mais alors  $D\cup A\cup S$ est une sphère qui
borde une boule $B^3$ d'un côté ou de l'autre. Dans un cas, on
remplace $M$ par $M\cup B^3$ et dans l'autre par $M\setminus B^3$
(où on est sûr qu'à isotopie près toutes les structures de $\SS$
coïncident d'après le théorème de
classification~\ref{t:boule} d'Eliashberg). À la fin du processus,
tous les disques du bord horizontal sont éliminés. Ceci implique
que toutes les composantes de $M\setminus \TT$ sont des tores
épaissis, des voisinages de bouteilles de Klein ou des anneaux
épaissis qui sont collés de façon incompressible les uns aux
autres.
\end{proof}

\begin{remarque}\label{loss}
En éliminant les disques de contact, on perd le contrôle sur la
fibration legendrienne, bien qu'elle persiste topologiquement.
Ainsi, au lieu de considérer les classes de conjugaison de
structures de contact tendues \emph{ajustées} à un triplet $(M,\tau
,\zeta )$, on doit considérer les classes de conjugaison de
structures de contact sur $V$ qui coïncident avec $\zeta$ sur
$V\setminus \Int (N)$.  Du fait de cette perte
d'information, on doit s'appuyer sur un résultat de finitude pour
une classe simple de variétés $N$, précisément
lorsque $N$ est un fibré en cercles au-dessus
d'une surface.
\end{remarque}

\subsection{Réduction aux fibrés en cercles}

Soit  $(N,\zeta)=(N_i,\zeta_i)$ une paire comme dans la
proposition~\ref{no disks of contact} et $\TT (N,\zeta )$ un
ensemble de structures de contact tendues sur $V$ égales à
$\zeta$ sur $V\setminus \Int (N)$ et de torsions bornées par $n$.
L'objectif de cette partie est de démontrer la proposition
suivante.

\begin{proposition}\label{p:fibre}
Il existe une famille finie de paires $(N_1 ,\zeta_1 ),...,(N_l
,\zeta_l )$ satisfaisant aux conditions suivantes~:
\begin{itemize}
\item les $N_i$, $i=1,...,l$, sont des fibrés en cercles
au-dessus de surfaces compactes à bord, éventuellement non
orientables, dans $V$\,;
\item $\zeta_i$, $i=1,...,l$, est une structure de contact tendue
sur $V\setminus \Int (N_i )$\,;
\item toutes les composantes de $\Gamma_{\partial N_i}$ sont isotopes
aux fibres de $N_i$\,;
\item toute structure de $\TT (N,\zeta )$ est conjuguée,
par un produit de twists de Dehn et d'isotopies, à une structure
égale à $\zeta_i$ sur $V\setminus \Int (N_i )$ pour un certain
$i\in [1,l]$.
\end{itemize}
\end{proposition}

On observe déjà que la sous-variété $N$ n'est pas quelconque.

\begin{lemme}
La variété $N$ est une variété graphée constituée par
collage, le long d'une famille de tores incompressibles, de blocs de
l'un des trois types suivants~:
\begin{enumerate}
\item un tore épais\,;
\item un voisinage de bouteille de Klein\,;
\item un fibré en cercles au-dessus d'une surface compacte de
caractéristique inférieure ou égale à $-1$ et non
nécéssairement orientable. Toute composante fibrée $M$ de ce
type a un bord qui intersecte non trivialement $\partial N$, et les
courbes de $\Gamma_{(\partial N)\cap M} (\zeta )$ sont isotopes aux
fibres de $M$.
\end{enumerate}
\end{lemme}

\begin{remarque}
La variété $N$ n'est pas forcément connexe.
\end{remarque}

\begin{proof}[Démonstration.]
À chaque fois que deux composantes $T^2\times I$ font bord commun,
elles peuvent être regroupées pour former un unique  $T^2\times
I$. Si on coupe $N$ le long de l'union des $T^2\times \{{1\over
2}\}$ donnés par la proposition~\ref{no disks of contact}, alors
les composantes connexes sont de l'un des trois types recherchés.
Le cas 3 apparaît pour les composantes connexes qui contiennent
un produit $A\times I$. Sur chaque composante $M$ de ce découpage
qui est de la forme 3, les fibres sont isotopes aux courbes
$\partial A \times \{ *\}$ des produits $A\times I$ qui la composent
et qui sont de degré nul. Comme les anneaux $(\partial A) \times
I$ contiennent une courbe de singularités isotope à leur âme
(et donc aux fibres de $M$ d'après l'assertion~\ref{a:degre}), les
courbes de $\Gamma_{\partial N} (\zeta ) \cap M$  sont isotopes aux
fibres.
\end{proof}

\begin{lemme}\label{lemme : tore}
Soit $T \subset N$  un tore de recollement entre deux composantes
fibrées $M$ et $M'$ de type 3 le long duquel les fibrations ne
correspondent pas. Il existe un voisinage tubulaire $T\times [-1,1]$
de $T =T\times \{ 0\}$ inclus dans $M\cup M'$ et des structures
$\zeta_1 ,...,\zeta_k$ sur $T\times [-1,1]$ pour lesquelles les
tores $T\times \{ \pm1 \}$ sont convexes, les composantes des
courbes de découpage $\Gamma_{T\times \{ \pm 1\}} (\zeta_i )$ sont
respectivement deux fibres de $M'$ et pour lesquelles, pour tout $\xi
\in \TT (N,\zeta )$, il existe une isotopie de $\xi$ en $\xi'$ dans
$V$ relative à $V\setminus \Int (N)$  qui fait co\"incider $\xi'
\vert_{T\times [-1,1]}$ avec une des structures $\zeta_i$.
\end{lemme}

\begin{proof}[Démonstration.]
On commence par prouver le fait suivant.

\begin{lemme}\label{l:epaississement}
Il existe un épaississement $T\times [-2,2]$ de  $T$ dans $M \cup
M'$ tel que toute composante de $\Gamma_{T\times \{ -2\}} (\xi)$
soit isotope à une fibre de $M$ et toute composante de
$\Gamma_{T\times \{ 2\}} (\xi )$ soit isotope à une fibre de $M'$.
\end{lemme}

\begin{proof}[Démonstration.]
Les intersections $\partial M \cap \partial N$ et $\partial M' \cap
\partial N$ sont non vides et contiennent une courbe de
singularités isotope aux fibres. On prend des copies parallèles
de ces courbes singulières dans $\Int (M)$ et $\Int (M')$. Ce sont
des courbes legendriennes isotopes aux fibres dont l'enroulement
(calculé par rapport à la fibration) est nul. On choisit alors
un épaississement $T\times [-2,2]$ de $T=T\times \{ 0\}$ qui
incorpore ces courbes legendriennes dans son bord. La condition
d'enroulement nul fait que celles-ci peuvent être réalisées,
après isotopie éventuelle, comme des orbites non singulières
de $\xi T\times \{ \pm2\}$. Elles y imposent donc la condition
voulue sur $\Gamma_{T\times \{ \pm 2\}} (\xi)$.
\end{proof}

D'après la classification des structures tendues sur le tore
épais~\cite{Gi4,Ho1}, si les pentes des courbes de découpage au
bord de $T\times [-2,2]$ ne coïncident pas, on peut réduire
à deux  le nombre de leurs composantes. On peut donc, sous les
hypothèses du lemme~\ref{lemme : tore} et les conclusions du
lemme~\ref{l:epaississement}, isotoper toute structure  $\xi
\vert_{T\times [-2,2]}$ relativement au bord de $T\times [-2,2]$ en
$\xi'$ de sorte que les découpages $\Gamma_{T\times \{ \pm 1\}
}(\xi' )$ aient deux composantes  isotopes aux fibres de $M$ et de
$M'$. On peut également obtenir que la torsion de $\xi'
\vert_{T\times [-1,1]}$ soit nulle.  Il n'y a qu'un nombre fini de
structures de contact tendues sur $T\times [-1,1]$ qui satisfassent
ces conditions. D'où les structures $\zeta_1 ,...,\zeta_k$
recherchées.
\end{proof}

On a un résultat dans le même esprit pour les composantes
$N(K)$. On rappelle que la bouteille de Klein possède seulement
quatre classes d'isotopie distinctes de courbes fermées simples
essentielles et non orientées. Précisément, dans un système
de coordonnées $(y,t) \in (\R /\Z ) \times [0,1] $ où $K\simeq
(\R /\Z )\times [0,1] /(y,0)\sim (-y,1)$, on obtient quatre
représentants avec les courbes $\gamma_1 =\{ t=0 \}$, $\gamma_2
=\{ y=0\}$, $\gamma_3 =\{ y=\frac{1}{2}\}$ et $\gamma_4 =\{
y=\frac{1}{4} , y=\frac{3}{4}\}$. Similairement, $N(K)$ possède,
à isotopie près, deux fibrations de Seifert différentes~: une
fibration en cercles sur la bande de M\"obius dont les fibres sont
isotopes à $\gamma_1$ et une fibration de Seifert sur le disque
dont les fibres régulières sont isotopes à $\gamma_4$ et les
deux fibres singulières à $\gamma_2$ et $\gamma_3$.

\begin{lemme}\label{l:N(K)}
Soit $T$ un tore de recollement entre une composante fibrée $M$ de
type 3 et un voisinage $N(K)$. Si aucune des deux fibrations
possibles de $N(K)$ ne prolonge celle de $M$, alors il existe un
voisinage $T\times [-1,1]$ de $T=T\times \{ 0\}$ dans $M\cup N(K)$
et une famille $\zeta_1 ,...,\zeta_k$ de structures de contact sur
$T\times [-1,1]$ avec les propriétés suivantes~:
\begin{itemize}
\item les tores $T\times \{ \pm1\}$ sont $\zeta_i$-convexes et
les courbes de découpages $\Gamma_{T\times \{ \pm 1\}} (\zeta_i )$
sont respectivement deux fibres de $M$ et de $N(K)$ (pour une des
deux fibrations de $N(K)$)\,;
\item toute structure $\xi \in \TT (N,\zeta )$ est isotope à une
structure $\xi'$ égale à l'une des structures $\zeta_i$ sur
$T\times [-1,1]$ par une isotopie qui est stationnaire en dehors de
$N$.
\end{itemize}
\end{lemme}

\begin{proof}[Démonstration.]
On commence par normaliser toute structure $\xi \in \TT (N,\zeta )$
au voisinage de $K$. Pour cela, on rend, par isotopie de $\xi$, la
courbe $\gamma_1$ legendrienne, puis on maximise son invariant de
Thurston-Bennequin relatif à $K$, parmi ceux qui sont négatifs
ou nuls. On distingue deux cas.

Si ce nombre de Bennequin maximal est nul, alors $\gamma_1$ est,
après isotopie de $\xi$, une courbe non singulière du
feuilletage caractéristique de $K$. On peut alors trouver un
voisinage $T\times [-2,2]$ de $T=T\times \{ 0\}$ dans $M\cup N(K)$
pour lequel le feuilletage caractéristique $\xi T\times \{ 2\}$
contient une courbe non singulière isotope à $\gamma_1$,
c'est-à-dire aux fibres de la fibration de $N(K)$ sur la bande de
M\"obius, et $\xi T\times \{ -2\}$ contient une courbe non
singulière isotope aux fibres de $M$. Comme on l'a déjà vu
précédemment, la classification des structures tendues sur le
tore épais permet, dans le cas où ces deux directions ne sont
pas isotopes dans $T^2 \times [-2,2]$, d'isotoper $\xi$ dans
$T\times [-2,2]$ pour l'amener sur un nombre fini de modèles
fixés au préalable dans $T\times [-1,1]$, tous possédant les
propriétés requises dans les conclusions du lemme~\ref{l:N(K)}.

Lorsque le nombre de Bennequin maximal de $\gamma_1$ est strictement
négatif, on isotope $K$ pour que $\xi$ ne fasse pas de demi-tour
inversé le long de $\gamma_1$. On découpe alors $K$ selon
$\gamma_1$. Il apparaît un anneau $A$ à bord legendrien que
l'on rend convexe par isotopie $\classe\infty$-petite relative à
$\partial A$. Le découpage de $A$ est alors constitué de
traverses. Ces traverses donnent dans $K$ une multi-courbe $\Gamma$
qui joue le rôle de \guil{courbe de découpage} de $K$. Ici, on
observe que chacune des composantes de $\Gamma$ est isotope à
$\gamma_2$, $\gamma_3$ ou $\gamma_4$. De plus, le lemme de
réalisation appliqué sur $A$ permet de dessiner sur $K$ des
courbes legendriennes non singulières parallèles aux composantes
de $\Gamma$. Celles-ci sont isotopes aux fibres régulières de la
fibration de $N(K)$ sur le disque et peuvent être incorporées
dans le bord d'un voisinage de $T$. On conclut alors comme
précédemment.
\end{proof}

\begin{proof}[Démonstration de la proposition~\ref{p:fibre}.]
Si $\partial N =\emptyset$, alors soit $N(=V)$ est un fibré en
tores sur le cercle, soit $N$ est obtenu en recollant deux
voisinages $N(K)$ le long de leur bord. Le cas des fibrés en
cercles est traité dans~\cite{Gi4} et~\cite{Ho2}. On y obtient le
théorème~\ref{t:torsion} et donc on peut prendre $N_i
=\emptyset$ dans la proposition~\ref{p:fibre}. Celui où
$V=N=N(K)\cup N(K)$ s'étudie de manière similaire.

\begin{proposition} \label{prop: double bouteille}
Soit $V=N(K_1)\cup N(K_2)$ où $K_1$ et $K_2$ sont des bouteilles
de Klein. L'espace des structures de contact tendues de torsion
inférieure à $n\in \N$ possède un nombre fini d'orbites sous
l'action du groupe engendré par les isotopies et les twists de
Dehn.
\end{proposition}

\begin{proof}[Esquisse de démonstration.]
Soit $\xi$ une structure de contact tendue sur $V$. On normalise le
\guil{découpage} de $K_1$ puis celui de $K_2$ en minimisant le
nombre de leurs composantes. On choisit des petits voisinages
\guil{homogènes} $N_\varepsilon(K_i)$ de $K_i$ et on identifie
l'adhérence de leur complémentaire à $T\times[-2,2]$. Les
découpages $\Gamma_{\partial N_\varepsilon (K_i)}$ sont
constitués de fibres de l'une des fibrations des $N(K_i)$. Si
aucune des deux fibrations de $N(K_1)$ ne s'étend à $N(K_2)$,
alors on peut trouver une isotopie de $\xi$ en $\xi'$ pour laquelle
les découpages $\Gamma_{T\times \{ \pm 1\}}(\xi')$ sont
constitués de deux fibres de la composante $N(K_i)$
correspondante. Si le nombre de composantes de $\Gamma_{\partial
N_\varepsilon (K_i)}$ n'était pas deux, on trouverait une rocade
s'appuyant sur $\partial N_\varepsilon (K_i)$ et par suite une
rocade s'appuyant sur $K_i$ qui permettrait de réduire son
découpage. D'où une contradiction.  On peut ainsi se ramener au
cas où toutes les structures considérées coïncident au
voisinage des bouteilles $K_i$. La proposition découle alors de la
classification des structures tendues sur le tore épais.

Le cas restant est celui où $V$ est un fibré de Seifert sur la
bouteille de Klein, l'espace projectif où la sphère, avec
respectivement $0$, $2$ et $4$ fibres singulières. Pour toute
structure tendue $\xi$, la variété $V$ possède de plus une
fibre régulière d'enroulement nul. On excise un voisinage
standard d'une telle fibre régulière pour ce ramener au cas d'un
fibré en cercles sur une surface à bord, lequel est traité
dans la section~\ref{section: cercle}.
\end{proof}

On suppose maintenant que $\partial N\not= \emptyset$.  Si deux
composantes de type 3 de $N$ sont adjacentes le long d'un tore $T$,
soit leurs fibrations coïncident au bord, et on les recolle pour
former une composante fibrée élargie, soit elles sont
différentes et, à l'aide du lemme~\ref{lemme : tore}, on se
ramène à l'étude d'une famille de structures tendues sur le
découpage de $N$ par $T\times [-1,1]$ qui coïncident au bord.
Comme $\partial N\neq \emptyset$, toute composante $N(K)$, voisinage
d'une bouteille de Klein $K$, est alors isolée ou adjacente à
une composante $M$ de type 3. Dans le deuxième cas, le
lemme~\ref{l:N(K)} permet, quitte à effectuer un découpage le
long du tore de recollement $T$ (et après isotopie des
structures), de se ramener aux cas de $M$ et d'une composante $N(K)$
isolée ou de supposer que la fibration de $M$ s'étend en une
fibration (de Seifert) sur $M \cup N(K )$. Pour obtenir la
proposition~\ref{p:fibre}, les composantes $N(K)$ isolées peuvent
être traitées à part, de la même façon que dans la
proposition~\ref{prop: double bouteille}.

Après ces opérations de réduction, la sous-variété $N$ est
 un fibré de Seifert au-dessus d'une surface compacte à bord
non vide.  Pour toute structure $\xi \in \TT (N,\zeta )$, le fibré
$N$ possède une fibre régulière d'enroulement nul. On se
ramène au cas où $N$ est un fibré en cercles en retirant à
$N$ des voisinages normalisés $N(\gamma_1 ),...,N(\gamma_l )$ des
fibres singulières (rendues legendriennes au préalable)
$\gamma_1$,...,$\gamma_l$. Pour cela, on incorpore dans le bord de
chaque $N(\gamma_i )$ une fibre régulière d'enroulement nul qui
forme une orbite périodique de $\xi \partial N(\gamma_i )$.  Ceci
impose que $\Gamma_{\partial N(\gamma_i)}$ est constitu\' e de
courbes isotopes aux fibres. On peut également imposer que $\#
(\Gamma_{\partial N(\gamma_i )} ) =2$ et d'après le théorème
de classification des structures tendues sur le tore
solide~\cite{Gi4,Ho1}, il y a un nombre fini de structures tendues
sur $N(\gamma_i)$ qui ont un tel découpage au bord.

Enfin, maintenant que l'on s'est ramené à un fibré en cercles,
on remarque que, à isotopie près, toutes les composantes des
courbes de découpage présentes au bord sont $\S^1$-invariantes.
En conclusion, on obtient la proposition~\ref{p:fibre}.
\end{proof}

\subsection{Le cas des fibrés en cercles} \label{section: cercle}

Pour conclure la preuve du théorème~\ref{t:torsion}, on utilise
des résultats de~\cite{Gi4, Gi5} et~\cite{Ho2} qui classifient les
structures de contact tendues sur les fibrés en cercles~:
l'ensemble $\TT (N,\zeta )$ contient un nombre fini de structures de
contact tendues de torsion donnée, à conjuguaison près par des
produits de twists de Dehn le long de tores incompressibles dans $N$
et d'isotopies. La classification est même plus précise.

\begin{theoreme}\cite{Gi4,Gi5,Ho2}\label{t:GH}
Soit $S$ une surface compacte orientée de bord non vide et $E$ un
ensemble formé d'un nombre pair (et non nul) de points sur chaque
composante de $\partial S$. Soit également $\zeta$ un germe de
structure de contact le long de $\partial M$. On suppose que
$\partial M$ est $\zeta$-convexe et que sa courbe de découpage est
$S^1 \times E$. Les classes d'isotopies de structures de contact
tendues sur $M$ qui co\"incident avec $\zeta$ sur le bord sont en
correspondance bijective avec les sous-variétés incompressibles
(mais pas nécéssairement $\partial$-incompressibles) de
dimension un $\Gamma$  de $S$ avec $\partial \Gamma =E$.
\end{theoreme}

Dans ce théorème, il faut voir $\Gamma$ comme la courbe de
découpage \guil{minimale} d'une section $\{ 0 \} \times S$ rendue
convexe au préalable. Lorsque $\zeta$ est $S^1$-invariante, alors
toutes les structures tendues sur $M$ qui coïncident avec
$\zeta$ au bord le sont également (à isotopie près).

Dans le cas qui nous intéresse cependant, la base $S$ de $\pi
:N\to S$ n'est pas toujours orientable. On se donne une multi-courbe
compacte plongée dans $S$ dont le complémentaire est orientable
et dont la préimage par $\pi$ est une union disjointe de
bouteilles de Klein. Pour toute structure $\xi \in \TT (N,\zeta )$,
on peut réaliser une isotopie de $\xi$ à support dans $N$ pour
que chacune de ces bouteilles $K$ contienne une fibre de $\pi$
d'enroulement nul qui est une orbite non singulière $\gamma$ de
son feuilletage caratéristique. La surface obtenue en découpant
$K$ selon $\gamma$ est un anneau que l'on peut rendre convexe par
perturbation $\classe\infty$-petite de $\xi$ relative à $\gamma$.
Par extension, on considère que la bouteille $K$ est
\guil{convexe} et que son découpage s'identifie à celui de $A$.
Celui-ci est consitué, après isotopie, d'une union de fibres de
$\pi$. En appliquant le lemme de réalisation de feuilletages à
$A$, on arrive sans difficulté à conjuguer la structure $\xi$
sur  un voisinage $N(K)$ de $K$ à la structure d'équation $\sin
((2n+1)\pi t)dx+\cos ((2n+1)\pi t)dy=0$ pour l'identification
donnée par
$$N(K)\simeq \{ (x,y,t)\in [-1,1]\times \R /\Z  \times [0,1 ]\}
/(x,y,0)\sim (-x,-y,1 )$$
où $K=\{ x=0\}$.

Sur la section $B_0=\{ y=0\} $ de $\pi : N(K) \to B$, on peut
matérialiser une  courbe de découpage par les $2n+1$ traverses
$\{ y=0,\; t=\frac{2k+1}{4n+2} \}$, $k\in [0,2n+1[$. Même si $B_0$
n'est pas orientable, ce découpage prend un sens lorsqu'on
considère le disque  à bord legendrien $B_1$ obtenu en
découpant $B_0$ le long de $\{ t=0\}$. La structure $\xi$ est
$S^1$-invariante au-dessus de $B_1$.

Dans le complémentaire de l'union de ces voisinages de bouteilles
de Klein normalisés, qui est le  produit $M_0$ d'une surface
compacte orientable $S_0$ par le cercle, on  minimise, comme dans le
théorème~\ref{t:GH}, le découpage d'une section $\{ 0\} \times
S_0$ dont le bord est égal à celui des bandes  de M\"obius $B_0$
trouvées dans chaque $N(K)$. Comme elle l'est sur les composantes
de $\partial M_0$ (après utilisation du lemme de réalisation de
feuilletages), la structure $\xi$ est, d'après le
théorème~\ref{t:GH}, (isotope à une structure)
$S^1$-invariante sur $S^1  \times S_0$. On recolle la courbe de
découpage obtenue sur $\{ 0\} \times S_0$ avec les arcs tracés
sur les $B_0$. On obtient ainsi une \guil{courbe de découpage}
$\Gamma_S (\xi )$ de $S$. Soit $k(\xi )$ le nombre de composantes de
$\Gamma_S (\xi )$.

Si $k (\xi )$ est borné indépendemment de $\xi \in \TT (N,\zeta
)$, alors toutes les configurations possibles pour la famille
$(\Gamma_S (\xi ))_{\xi \in \TT (N,\zeta )}$ se déduisent d'un
nombre fini d'entre elles par des twists de Dehn et des isotopies de
$S$. Comme $\xi$ est déterminée par $\Gamma_S (\xi )$ et que
tout twist de Dehn sur $S$ provient d'un twist le long d'un tore de
$N$, on obtient la conclusion.

Dans le cas contraire, pour tout $p\geq 0$, il existe une structure
$\xi_p \in \TT (N,\zeta )$ pour laquelle la courbe de découpage
$\Gamma_S (\xi_p )$ contient une famille de $p$ courbes
parallèles. Celle-ci est incluse dans un anneau $A \subset S$. La
fibration est triviale au-dessus de $A$, et sur $\pi^{-1} (A)$ la
structure $\xi_p$ est $S^1$-invariante. On vérifie facilement que
le tore épais $(\pi^{-1} (A), \xi \vert_{\pi^{-1} (A)})$ contient
un tore de torsion $[\frac{p}{2}] -1$, d'où une contradiction.

\section{Finitude pour les n\oe uds legendriens}

Pour finir, on démontre ici le théorème~\ref{t:noeuds}. Soit $\KK$ une classe 
topologique de n\oe uds dans $\S^3$ et $K$ un de ses représentants. On note 
$\zeta_0$ la structure de contact standard sur $\S^3$, $V$
la variété à bord $\S^3 \setminus \Int (N(K))$ et $\KK_{\zeta_0,n}$ l'ensemble des
n\oe uds legendriens de $(\S^3, \zeta_0)$ dans la classe $\KK$ dont l'invariant
de Thurston-Bennequin vaut $n$. Tout n\oe ud $L \in \KK_{\zeta_0,n}$ possède un 
voisinage tubulaire $N(L) = \D^2 \times \S^1$ contenant $L$ comme $\{0\} \times
\S^1$ et \guil{standard} au sens où $\zeta_0$ y est définie par
$$ \sin \theta \, dx + \cos \theta \, dy = 0, \qquad 
   (x,y,\theta) \in \D^2 \times \S^1. $$ 
Pour tout $L \in \KK_{\zeta_0,n}$, la restriction de $\zeta_0$ à $\S^3 \setminus
\Int (N(L))$ donne une structure de contact $\xi_L$ sur $V$. La feuilletage
caractéristique $\xi_L\, \partial V$ est déterminé, à isotopie près, par la 
pente de ses lignes singulières. Cette pente elle-même est déterminée par 
l'invariant de Thurston-Bennequin $n$ de $L$. On peut donc supposer que toutes 
les structures $\xi_L$, $L \in \KK_{\zeta_0,n}$, coïncident le long de $\partial
V$. Or la variété $V$ est irréductible et les structures $\xi_L$ sont tendues et
de torsion nulle~: d'après le théorème~\ref{t:torsion-bord}, 
elles vivent dans un nombre fini de classes, modulo isotopie et conjugaison par 
des twists de Dehn sur des tores. Vus dans $\S^3$, ces twists sont isotopes à 
l'identité. Si $\xi_L$ et $\xi_{L'}$ sont dans une même classe, il existe donc 
une isotopie de $\S^3$ issue de l'identité dont le temps $1$ envoie $L$ sur $L'$
et préserve $\zeta_0$. Comme l'espace des structures de contact tendues 
positives sur $\S^3$ est simplement connexe \cite{El2}, cette isotopie peut être
déformée en une isotopie de contact qui amène $L$ sur $L'$. C'est une isotopie 
legendrienne entre $L$ et $L'$.

\begin{question}
Dans une classe d'isotopie de n\oe uds donnée, existe-t-il un nombre fini de
n\oe uds legendriens dont tout n\oe ud legendrien se déduise par isotopie
legendrienne et stabilisation~?
\end{question}

$ $\\
Université de Nantes, UMR 6629 du CNRS, 44322 Nantes, France\\
Vincent.Colin@math.univ-nantes.fr\\
$ $\\
École normale supérieure de Lyon, 69000 Lyon, France\\
Emmanuel.Giroux@umpa.ens-lyon.fr\\
$ $\\
University of Southern California, Los Angeles, CA 90089\\
khonda@math.usc.edu\\
http://rcf.usc.edu/\char126 khonda


\begin{thebibliography}{convex}

\bibitem[Be]{Be}
D.~Bennequin, 
\textit{Entrelacements et équations de Pfaff}. \ 
Astérisque \textbf{107--108} (1983), 87--161.

\bibitem[Br]{Br}
M. Brittenham, 
\textit{Essential laminations and Haken normal forms}. \ 
Pac. J. Math. \textbf{168} (1995), 217--234.

\bibitem[Co1]{Co1}
V.~Colin, 
\textit{Chirurgies d'indice un et isotopies de sphères dans les variétés de 
contact tendues}. \ 
C.~R. Acad. Sci. Paris Sér. I Math. \textbf{324} (1997), 659--663.

\bibitem[Co2]{Co2}
V.~Colin,
\textit{Sur la torsion des structures de contact tendues}. \ 
Ann. Sci. Éc. Norm. Sup. \textbf{34} (2001), 267--286.

\bibitem[Co3]{Co3}
V.~Colin, 
\textit{Une infinité de structures de contact tendues sur les variétés 
toroïdales}. \ 
Comment. Math. Helv. \textbf{76} (2001), 353--372.

\bibitem[Co4]{Co4}
V.~Colin,
\textit{Structures de contact tendues sur les variétés toroïdales et 
approximation de feuilletages sans composante de Reeb}. \ 
Topology \textbf{41} (2002), 1017--1029.

\bibitem[CGH1]{CGH1}
V.~Colin, E.~Giroux et  K.~Honda, 
\textit{On the coarse classification of tight contact structures}. \ 
Proceedings of Symposia in Pure Mathematics \textbf{71} (2003), 109--120.

\bibitem[El1]{El1}
Y.~Eliashberg, 
\textit{Classification of over-twisted contact structures on $3$-manifolds}. \ 
Invent. Math. \textbf{98} (1989), 623--637.

\bibitem[El2]{El2}
Y.~Eliashberg,
\textit{Contact $3$-manifolds twenty years since J.~Martinet's work}. \ 
Ann. Inst. Fourier \textbf{42} (1992), 165--192.

\bibitem[El3]{El3}
Y.~Eliashberg, 
\textit{Filling by holomorphic discs and its applications}. \ 
London Math. Soc. Lect. Notes Ser. \textbf{151} (1991), 45--67.

\bibitem[ET]{ET}
Y.~Eliashberg et W.~Thurston,
\textit{Confoliations}. \ 
University Lecture Series \textbf{13}, Amer. Math. Soc. 998).

\bibitem[Et]{Et}
J.~Etnyre,
\textit{Tight contact structures on lens spaces}. \ 
Commun. Contemp. Math. \textbf{2} (2000), 559--577.

\bibitem[EH1]{EH1}
J.~Etnyre et K.~Honda,
\textit{On the non-existence of tight contact structures}. \ 
Ann. of Math. \textbf{153} (2001), 749--766.

\bibitem[EH2]{EH2}
J.~Etnyre et K.~Honda,
\textit{Tight contact structures with no symplectic fillings}. \ 
Invent. Math. \textbf{148} (2002), 609--626.

\bibitem[FO]{FO}
W.~Floyd et U.~Oertel, 
\textit{Incompressible surfaces via branched surfaces}. \ 
Topology \textbf{23} (1984), 117--125.

\bibitem[Ga]{Ga}
D.~Gabai, 
\textit{Essential laminations and Kneser normal form}. \ 
J.~Diff. Geom. \textbf{53} (1999), 517--574.

\bibitem[Gi1]{Gi1}
E.~Giroux, 
\textit{Convexité en topologie de contact}. \ 
Comment. Math. Helv. \textbf{66} (1991), 637--677.

\bibitem[Gi2]{Gi2}
E.~Giroux, 
\textit{Une structure de contact, même tendue, est plus ou moins tordue}. \  
Ann. Sci. Éc. Norm. Sup. \textbf{27} (1994), 697--705.

\bibitem[Gi3]{Gi3}
E.~Giroux, 
\textit{Une infinité de structures de contact tendues sur une infinité de 
variétés}. \ 
Invent. Math. \textbf{135} (1999), 789--802.

\bibitem[Gi4]{Gi4}
E.~Giroux, 
\textit{Structures de contact en dimension trois et bifurcations des 
feuilletages de surfaces}. \ 
Invent.Math.\ \textbf{141} (2000), 615--689.

\bibitem[Gi5]{Gi5}
E.~Giroux, 
\textit{Structures de contact sur les variétés fibrées en cercles au-dessus 
d'une surface}. \ 
Comment. Math. Helv. \textbf{76} (2001), 218--262.

\bibitem[Gi6]{Gi6}
E.~Giroux, 
\textit{Structures de contact, livres ouverts et tresses fermées}. \ 
En préparation.

\bibitem[Gr]{Gr}
M.~Gromov, 
\textit{Pseudoholomorphic curves in symplectic manifolds}. \ 
Invent. Math. \textbf{82} (1985), 307--347.

\bibitem[Ha]{Ha}
W.~Haken, 
\textit{Theorie der Normalflächen. Ein Isotopiekriterium für den Kreisknoten}. \ 
Acta Math. \textbf{105} (1961), 245--375.

\bibitem[Ho1]{Ho1}
K.~Honda, 
\textit{On the classification of tight contact structures I}. \ 
Geom. Topol. \textbf4 (2000), 309--368.

\bibitem[Ho2]{Ho2}
K.~Honda, 
\textit{On the classification of tight contact structures II}. \ 
J.~Diff. Geom. \textbf{55} (2000), 83--143.

\bibitem[Ho3]{Ho3}
K.~Honda, 
\textit{Gluing tight contact structures}. \ 
Duke Math. J. \textbf{115} (2002), 435--478.

\bibitem[HKM]{HKM}
K.~Honda, W.~Kazez et  G.~\textsc{Mati\'c}, 
\textit{Convex decomposition theory}. \  
Internat. Math. Res. Notices \textbf{2002}, 55--88.

\bibitem[Ka]{Ka}
Y.~Kanda, 
\textit{The classification of tight contact structures on the $3$-torus}. \ 
Comm. Anal. Geom. \textbf{5} (1997), 413--438.

\bibitem[Kn]{Kn}
H.~Kneser, 
\textit{Geschlossene Flächen in dreidimensionalen Mannigfaltigkeiten}. \  
Jahres. der Deut. Math.-Verein. \textbf{38} (1929), 248--260.

\bibitem[KM]{KM}
P.~Kronheimer et  T.~Mrowka, 
\textit{Monopoles and contact structures}. \ 
Invent. Math. \textbf{130} (1997), 209--255.

\bibitem[Li]{Li}
T.~Li, 
\textit{Laminar branched surfaces in $3$-manifolds}. \ 
Geom. Topol. \textbf{6} (2002), 153--194.

\bibitem[LS]{LS}
P.~Lisca et A.~Stipsicz, 
\textit{An infinite family of tight, not semi-fillable contact three-manifolds}. \ 
Geom. Topol. \textbf7 (2003), 1055--1073.

\bibitem[Ma]{Ma}
S.~\textsc{Makar-Limanov}, 
\textit{Morse surgeries of index $0$ on tight manifolds}. \ 
Prépublication (1997).

\bibitem[Th]{Th}
W.~Thurston,
\textit{A norm for the homology of $3$-manifolds}. \ 
Mem. Amer. Math. Soc. \textbf{59} (1986), 99--130.

\bibitem[Wa]{Wa}
F.~Waldhausen, 
\textit{On irreducible $3$-manifolds which are sufficiently large}. \ 
Ann. of Math. \textbf{87} (1968), 56--88.

\bibitem[Wh]{Wh}
J.~H.~C.~Whitehead, 
\textit{On $\classe1$ complexes}. \ 
Ann. of Math. \textbf{41} (1940), 809--824.

\end{thebibliography}
\end{document}